\documentclass[twoside,11pt]{article}

\RequirePackage{amsthm}
\usepackage{jmlr2e}

\usepackage{graphicx}

\RequirePackage[OT1]{fontenc}

\usepackage{hyperref}
\usepackage{multicol}
\usepackage{graphicx}
\usepackage{multirow}
\usepackage{booktabs}
\usepackage{mathrsfs,amsmath,amssymb,color}
\usepackage{tcolorbox}
\usepackage{bm}
\usepackage{bbm, dsfont}
\usepackage{subcaption}
\theoremstyle{plain}
\newtheorem{thm}{Theorem}
\theoremstyle{plain}
\newtheorem{lem}[theorem]{Lemma}
\theoremstyle{plain}
\newtheorem{ass}[theorem]{Assumption}
\theoremstyle{plain}
\newtheorem{prop}[theorem]{Proposition}

\newcommand{\be} {\begin{eqnarray*}}
\newcommand{\ee} {\end{eqnarray*}}

\def\boxit#1{\vbox{\hrule\hbox{\vrule\kern6pt  \vbox{\kern6pt#1\kern6pt}\kern6pt\vrule}\hrule}}

\def\floorx{\left \lfloor{x}\right \rfloor}
\def\floorbeta{\left \lfloor{\beta}\right \rfloor}

\newcommand{\Var}{\mathop{\rm var}}
\newcommand{\Cov}{\mathop{\rm cov}}

\definecolor{newred}{RGB}{156, 34, 54}
\def\tcr{ \color{black} }

\newcommand{\norm}[1]{\left\Vert#1\right\Vert}

\newcommand\wt[1]{{ \widetilde{#1} }}

\newcommand \bbE{\mathbb{E }}

\def\T{{ \mathrm{\scriptscriptstyle T} }}

\newcommand{\ind}{\mathbbm{1}}
\def\bsm{\boldsymbol}

\def\m{\mathcal}
\def\mb{\mathbb}

\def\boxit#1{\vbox{\hrule\hbox{\vrule\kern6pt  \vbox{\kern6pt#1\kern6pt}\kern6pt\vrule}\hrule}}

\jmlrheading{1}{2023}{1-59}{12/21; Revised 1/23}{10/00}{}{Shuang Zhou, Debdeep Pati, Tianying Wang, Yun Yang, and Raymond J. Carroll}

\ShortHeadings{Gaussian Processes with Errors in Variables}{Zhou, Pati, Wang, Yang and Carroll}
\firstpageno{1}

\begin{document}

\title{Gaussian Processes with Errors in Variables: Theory and Computation}

\author{\name Shuang Zhou \email szhou98@asu.edu \\
       \addr School of Mathematical and Statistical Sciences\\
       Arizona State University\\
       Tempe, AZ 85287-1804, USA
       \AND
       \name  Debdeep Pati  \email debdeep@stat.tamu.edu \\
       \addr Department of Statistics\\
       Texas A\&M University\\
       College Station, TX 77843-3143, USA
       \AND
       \name  Tianying Wang \email tianyingw@tsinghua.edu.cn\\
       \addr Center for Statistical Science\\
       Department of Industrial Engineering\\
       Tsinghua University\\
      Beijing 100084, China
              \AND
       \name  Yun Yang  \email yy84@illinois.edu \\
       \addr Department of Statistics\\
       University of Illinois at Urbana-Champaign\\
       Champaign, IL 61820-3633, USA 
              \AND
       \name   Raymond J.\ Carroll  \email carroll@stat.tamu.edu \\
       \addr Department of Statistics\\
       Texas A\&M University\\
       College Station, TX 77843-3143, USA\\
       School of Mathematical and Physical Sciences\\
       University of Technology Sydney\\
       Ultimo NSW 2007, Australia}

\editor{}

\maketitle

\begin{abstract}
Covariate measurement error in nonparametric regression is a common problem in nutritional epidemiology and geostatistics, and other fields.  Over the last two decades, this problem has received substantial attention in the frequentist literature. Bayesian approaches for handling measurement error have only been explored recently and are surprisingly successful, although there still is a lack of a proper theoretical justification regarding the asymptotic performance of the estimators.    By specifying a Gaussian process prior on the regression function and a Dirichlet process Gaussian mixture prior on the unknown distribution of the unobserved covariates, we show that the posterior distribution of the regression function and the unknown covariate density attain optimal rates of contraction adaptively over a range of H\"{o}lder classes, up to logarithmic terms.  
We also develop a novel surrogate prior for approximating the Gaussian process prior that leads to efficient computation and preserves the covariance structure, thereby facilitating easy prior elicitation.  We demonstrate the empirical performance of our approach and compare it with competitors in a wide range of simulation experiments and a real data example.
\end{abstract}

\begin{keywords}
Approximated Gaussian processes, measurement error model, nonparametric Bayes, smoothing and nonparametric regression, supersmooth errors 
\end{keywords}

\section{Introduction}
The general formulation of a deconvolution problem assumes that the observations are the true underlying variables contaminated with measurement error.   In an errors-in-variables regression problem, responses $Y_i$'s are observed corresponding to evaluations of an unknown regression function $f_0$ on noise-contaminated covariates $W_i$'s as
\begin{equation}\label{eq:model}
\begin{aligned}
Y_i &= f_0(X_i) + \epsilon_i, \quad \epsilon_i \overset{i.i.d.}{\sim} \mbox{N}(0, \sigma^2), \\
W_i &=X_i + u_i, \quad u_i \overset{i.i.d.}{\sim} g_u, \quad X_i\overset{i.i.d.}{\sim} p_0,\quad i=1, \ldots, n,
\end{aligned}
\end{equation}
where $X_i$'s are the unknown true covariates and we denote by $p_0$ the marginal distribution of the true covariate, and we write ``i.i.d." short for ``identically and independently distributed". In model \eqref{eq:model}, we consider the centered Gaussian error $\epsilon_i$ with unknown standard deviation $\sigma$, and denote by $g_u$ the known measurement error distribution. The goal is to recover the true regression function $f_0$ and the true density function $p_0$.

From a frequentist perspective, there is a rich literature addressing these problems.  Historically, the density deconvolution problem was first addressed in \cite{carroll1988optimal,fan1991optimal,stefanski1990deconvolving}, where it was noted that the fundamental difficulty in recovering the true density lies in the nature of the distribution of the measurement errors, and a class of deconvolution kernel density estimators was proposed. In a nonparametric regression setting \cite{fan1993nonparametric} developed a globally consistent deconvolution kernel type estimator. 
Later on, \cite{ioannides1997nonparametric} generalized the estimator while \cite{delaigle2007nonparametric} extended the theory to the heteroscedastic case.  Refer to a review article \citep{delaigle2014nonparametric} for a detailed discussion on kernel-based deconvolution estimators. Other methods such as deconvolution estimators based on Fourier-techniques, local linear and polynomial estimators are also popular, see \cite{carroll1996asymptotics,carroll1999nonparametric,cook1994simulation,SIMEXdelaigle08,delaigle2006nonparametric,delaigle2009design,du2011nonparametric,stefanski1995simulation}.

It is well known that the optimal rate of convergence of deconvolution estimators can be quite slow compared to the classical minimax rate for estimating smooth densities or functions.  The rate of convergence is controlled by the tail behavior of the characteristic function of the measurement error density; faster decaying rate of the characteristic function leads to a slower convergence rate and vice versa. In particular, the optimal rate is only of the logarithmic order when the measurement error distribution is a ``supersmooth" distribution, whose characteristic function decays exponentially in the tails. This includes the Gaussian and the Cauchy densities.  This slow rate of convergence renders estimation practically infeasible unless the measurement error variance is allowed to be sufficiently small \citep{carroll1999nonparametric,delaigle2008alternative,fan1992deconvolution} with respect to the sample size.    In particular, it has been shown in \cite{delaigle2008alternative,fan1992deconvolution} that 
the optimal rate of convergence in the ``supersmooth" case is improved to $n^{-\beta/(2\beta+1)}$ for estimating a function in a H\"{o}lder class with regularity level $\beta$ if the error standard deviation  of a Gaussian error density decreases to zero at the rate of $n^{-1/(2\beta+1)}$.  
This requirement on the error standard deviation can be easily satisfied by generating replicates $n^{1/(2\beta+1)}$ times per data point. In many applications, such as nutritional epidemiology, it is customary to collect multiple recalls of dietary intake from the respondents which serve as the replicated proxies and can boost the rate of convergence.  

Another critical point regarding the performance of classical deconvolution estimators is the choice of an appropriate kernel and associated bandwidth. 
Many effective bandwidth selection procedures have been developed for practical purposes, refer to \cite{delaigle2004bootstrap,delaigle2004practical,SIMEXdelaigle08}. 
 In absence of the knowledge of the true regularity level, data-driven bandwidth selection procedures using the Lepski's method are employed  \citep{comte2013anisotropic, kappus2014adaptive} with deconvolution kernel estimators, obtaining adaptivity with respect to the smoothness of the underlying function or density.  
Other types of the adaptive deconvolution estimator have been proposed, for instance, the ridge deconvolution estimator \citep{hall2007ridge} and the thresholding wavelet deconvolution estimator \citep{fan2002wavelet}. 

On the other hand, Bayesian procedures are naturally suited for general nonparametric regression tasks because of their ability to adapt to the unknown smoothness and to allow quantifications of uncertainty. For classical density estimation problems with no measurement error, Bayesian nonparametric techniques including Dirichlet process Gaussian mixture model \citep{escobar1995bayesian,ferguson1973bayesian,lo1984class} have demonstrated success in various applications, where the unknown density is modeled as a mixture of normals with a Dirichlet process prior on the mixing distribution. For the errors-in-variables regression estimation problem, \cite{berry2002bayesian} were the first to develop a fully Bayesian procedure for the nonparametric regression problem using smoothing splines and P-splines.  
Variants of spline-based models are developed in Bayesian framework to approximate the density function and/or variance function in the heteroscedastic case \citep{sarkar2014bayesian,staudenmayer2008density}. More recently, \cite{cervone2015gaussian} developed a Bayesian analysis for Gaussian processes with location errors using hybrid Monte-Carlo techniques.

Bayesian approaches have been demonstrated to be very successful numerically, however, there is a clear dearth of theoretical results justifying these approaches.  Few existing results for deconvolution density estimation are available recently in the Bayesian literature such as \cite{gao2016posterior,donnet2018posterior,rousseau2021wasserstein}. To the best of our knowledge, a formal theoretical justification for the use of Bayesian procedures in the errors-in-variables regression problem is missing. As the main contribution of this paper, we propose a fully Bayesian framework for the errors-in-variables regression using a Gaussian process prior, and develop a new theoretical framework for studying its frequentist properties including consistency and the quantification of posterior convergence rates. As mentioned earlier, the optimal rate in the errors-in-variables problem with Gaussian error distribution has been proved to be extremely slow, rendering  inference infeasible in applications. However, allowing the error variance to decrease to zero with sample size at an appropriate rate plays a very important role in improving the rate of convergence \citep{carroll1999nonparametric, fan1992deconvolution}.  
In this paper, we reexamine this situation from a Bayesian perspective assuming that the measurement error standard deviation decays at the order of of $n^{-1/(2\beta+1)}$ where $\beta$ is the smoothness of the true covariate density. However, we intend to maintain adaptivity with respect to the smoothness level of the true function and the true covariate density. 


As the main contribution, we show that in the errors-in-variables regression problem, when the Gaussian error variance decreases to zero at a certain rate, under appropriate regularity conditions on the true marginal density and regression function, the posterior distribution obtained from a suitably chosen hierarchical Gaussian process model with a Dirichlet process Gaussian mixture prior on the marginal density of the covariates converges to the ground truth at their respective minimax optimal rates, adaptively over a range of H\"{o}lder classes.  By viewing density deconvolution as an inverse problem \citep{knapik2011bayesian,ray2013bayesian}, we follow the general recipe in Theorem 3.1 of \cite{ray2013bayesian} as sufficient conditions for posterior convergence in our setting.
However, the work of \cite{knapik2011bayesian} is restricted to conjugate priors, \cite{ray2013bayesian} considers only periodic function deconvolution using wavelets, and substantial technical hurdles remain.
To address these challenges, we exploit the concentration properties of deconvolution kernel estimators to
construct 
test functions with exponentially small type-I and type-II error bounds for the testing problem
\begin{eqnarray} \label{eq: testing problem}
H_{0}: p= p_0, ~~~\mbox{vs}~~~ H_{A}: p\in \{p: d(p,p_0)>\xi_{n}\}.
\end{eqnarray}
\cite{ray2013bayesian} used concentration properties of thresholded wavelet based estimators based on standard results on concentration of Gaussian priors.
However, analogous results for kernel density estimators suited to density deconvolution problems are lacking.
One of our key technical contributions is to develop sharp concentration inequalities of the deconvolution kernel  estimators to construct tests in (\ref{eq: testing problem}).  

On the computational side, although Bayesian spline models are quite successful in practice, the choice of knots as well as the number of basis functions are critical to obtain good empirical performance. This stimulates the development of other Bayesian approaches for modeling the unknown function of interest such as Gaussian process priors. Gaussian processes are routinely used for function estimation in a Bayesian context. However, their use in the context of measurement error in nonparametric regression models is limited, since the unobserved values of covariates are involved in the prior covariance matrix of Gaussian process and is no longer conditionally independent given the data. To alleviate this issue in errors-in-variables regression problem, we develop an approximation to the Gaussian process as a prior for the unknown regression function.  The Gaussian process surrogate is computationally efficient as it avoids repeated computation of the matrix inversion. In addition to the appealing property of preserving the covariance kernel, we also show that the resulting surrogate process converges weakly to the original Gaussian process. This hints at the fact that  the good properties of the original posterior distribution will be subsequently inherited by the surrogate posterior.  For implementation, in addition to standard hyperparameters of a Gaussian process that control the smoothness of the sample paths, the Gaussian process surrogate contains a truncation parameter. Our result on the accuracy of such an approximation suggests that inference on the regression function is robust to the choice of the truncation parameter as long as it is chosen to be appropriately large.  Hence the approximation retains all the potential advantages of a Gaussian process.   
\subsection{Review on Nonparametric Regression with Errors in Variables}
Consider the regression model with errors in variables defined in Equation \eqref{eq:model}, where $ \{(Y_i, W_i), i =1, \dots, n\}$ are independent and identical draws from the joint unknown distribution. Recall that $Y_i$'s denote the observed responses and $W_i$'s are contaminated covariates. 
It is well known that in absence of any replicated proxy per data-point, the optimal rate for a ``supersmooth" error distribution is only of the logarithmic order, rendering the estimators to be highly inefficient for practical purposes \citep{fan1993nonparametric}. In cases where the error distribution remains unknown, it can be estimated from the repeated observations or extra validation data 
\citep{hall2007semiparametric,johannes2009deconvolution,neumann2007deconvolution}.  For the regular deconvolution kernel estimator, the deconvolution kernel function is constructed based on a suitable kernel function $K(\cdot)$ and the empirical estimator of the Fourier transform of the marginal density $p$ of covariates.
One can derive the deconvolution kernel density estimator \citep{fan1993nonparametric} for both the marginal density $p$ and the regression function $f$ by
\begin{eqnarray}
\widehat{p}_n (x) &=& \frac{1}{nh} \sum_{i=1}^{n}K_n\{(x-W_i)/h\},\label{eq: pn}\\
\widehat{f}_n(x) &=&  \frac{1}{nh} \sum_{i=1}^{n}K_n\{(x-W_i)/h\} Y_i / \widehat{p}_n(x),\label{eq: fn}\\
K_n(x)&=&\frac{1}{2\pi}\int e^{-itx} \frac{\phi_K (t)}{\phi_{u} (t/h)} dt.\label{eq:K_n}
\end{eqnarray}
$K_n(\cdot)$ is the deconvolution kernel function, $\phi_K(\cdot)$ and $\phi_{u}(\cdot)$ are the Fourier transforms of the kernel function $K(\cdot)$ and the density of measurement error $g_{u}(\cdot)$, respectively. Usually $\phi_K(\cdot)$ is assumed to be compactly supported to ensure that the deconvolution kernel $K_n(\cdot)$ is well defined. Also, to achieve the rate optimality one requires that kernel function $K(\cdot)$ is a $k$th-order kernel function where $k$ represents the regularity level of the true density function. 
However, in practice such deconvolution kernels typically do not admit closed-form expressions, and the estimation could suffer from extra errors due to numerical integrations. 

\subsection{Bayesian Nonparametric Regression with Errors in Variables}
In this article, we focus on the normal distribution $\mbox{N}(0,\delta^2)$ with an unknown variance $\delta^2$ as the measurement error distribution. We consider the following generic Bayesian hierarchical model for the nonparametric regression with errors in variables:
\begin{equation}
\begin{aligned}\label{eq:Bmodel}
&Y_i = f(X_i) + \epsilon_i, \quad \epsilon_i \sim \mbox{N}(0, \sigma^2), \\
&W_i =X_i + u_i, \quad u_i \sim \mbox{N}(0,\delta^2), \quad X_i\sim p, \quad i=1, \ldots, n,\\
&f \sim \Pi_f,\quad p\sim \Pi_p,\quad \sigma^2 \sim \Pi_{\sigma^2}, \quad \delta^2 \sim \Pi_{\delta^2}.
\end{aligned}
\end{equation} 
We assume that $Y_i$ is conditionally independent of $W_i$ given $X_i$, for $i=1, \dots, n$.  In the Bayesian framework, we obtain the posterior distribution of unknown parameters $\theta= (f, p, \delta)$ given the observed values $D_n= \{(Y_i, W_i), i =1, \dots, n\}$ via Bayes' rule:
\begin{equation*}
P(\theta \mid D_n) = \frac{P(D_n\, |\,\theta) \,P(\theta)}{P(D_n)}.
\end{equation*}
This posterior distribution $P(\theta \mid D_n)$ can then be used to conduct statistical inference on marginal density $p$ and regression function $f$, such as constructing point estimators and their associated credible intervals or bands. 
Variants of the model defined in Equation \eqref{eq:Bmodel} are used in the context of Bayesian methods in errors-in-variables regression problem \citep{berry2002bayesian,sarkar2014bayesian}.  Although for practical purposes we assume a prior distribution on $\delta^2$, in the theoretical investigation, to obtain the minimax-optimal convergence rate results, we assume $\delta^2$ to be known and let $\delta^2$ decrease to $0$ at a certain rate depending on $n$. 
For practical purpose we assign an objective prior on $\sigma^2$, the details of which can be found in Appendix~\ref{sec:gibbcom}. Whereas, for theoretical investigation we assume $\sigma = 1$ to simplify the analysis. Extension to general $\sigma$ is straightforward.  

By assigning proper priors on $f$ and $p$, we show that the estimation of $f$ and $p$ can be made adaptive, which means the prior does not demand any knowledge on the smoothness of the true regression function, and yet a nearly optimal rate of posterior contraction can be achieved as if the smoothness is known. Different from the deconvolution kernel estimator, a Bayesian method does not require explicitly constructing a deconvolution kernel function $K_n(\cdot)$, but the existence of such kernel is used for constructing the test function aforementioned in the introduction. The details of choosing specific priors for $f$ and $p$ are discussed in the following section. We start describing the Gaussian process prior for $f$ which requires specifying a covariance kernel analogous to the kernel $K(\cdot)$. 

\subsection{Prior Specifications}
We consider a Gaussian process prior \citep{rasmussen2004gaussian} as the prior $\Pi_f$ for $f$, which is a distribution over a space of functions such that the joint distribution of any finite evaluations of the random function is multivariate Gaussian. A Gaussian process is completely defined by a mean function $m(x) = E \{f(x)\}$ and a covariance kernel function $c(x,x') = \Cov \{f(x), f(x')\}$ for any $x,x' \in \mb R$. Therefore, any finite collection of random observation points $\{y_1(x_1),\ldots,y_N(x_N)\}$ at locations $x_1,\ldots,x_N$ has a joint Gaussian distribution given by
\begin{equation*}
\{y_1(x_1),\ldots,y_N(x_N)\}\sim  \mbox{N} (m, \tau^2\Sigma),
\end{equation*}
where $m=\{m(x_1),\ldots,m(x_N)\}$ and $\Sigma$ is the $N\times N$ covariance matrix with the $(i,j)$th element $\Sigma_{ij}= c(x_i,x_j)$. The mean function reflects the expected center of the realization, and the covariance kernel function reflects its fluctuation and local dependence.  The hyperparameter $\tau$ attached to the covariance kernel function controls the fluctuation magnitude. We use the notation $f(\cdot) \sim \textsc{gp}(m(\cdot), c(\cdot,\cdot))$ to denote that function $f$ follows a Gaussian process with mean function $m$ and covariance kernel function $c$. For the regular Gaussian process regression in the noised case with noise level $\sigma$, the predictive formula \citep{rasmussen2004gaussian} is 
\begin{equation}
\begin{aligned}\label{gp_formula}
f(X^{*})\mid X, Y, X^{*} &\sim \mbox{N}({\bar{f}}^{*}, \Cov\{f(X^{*})\}),\\
{\bar{f}}^{*} &= c(X^{*},X)\{c(X,X)+\sigma^2 I\}^{-1}Y,\\
\Cov\{f(X^{*})\} &= c(X^{*}, X^{*}) - c(X^{*}, X) \{c(X,X)+\sigma^2 I\}^{-1} c(X, X^{*}), 
\end{aligned}
\end{equation}
where $X, Y$ are the given data, $X^{*}$ is a new data point, $f(X^{*})$ is the prediction at $X^{*}$ and $c(X^{*},X)$ denotes the covariance matrix between $X^{*}$ and $X$. The posterior is a multivariate normal involved with the original data and the new data point. Refer to \cite{rasmussen2004gaussian} for a detailed explanation of a Gaussian process. Choice of the covariance kernel $c$ is crucial to obtain a desirable functional estimation. A squared exponential covariance or more generally, a Mat\'{e}rn covariance kernel are commonly used in practice. Also, the covariance kernel is often associated with hyperparameters which control the smoothness of the sample paths \citep{adler1990}. We shall discuss specific choices  in Section~\ref{sec:assfp}. 

It might appear on the surface that one can assume a parametric distribution for the unknown $X$ if the interest is solely on recovering the unknown function $f$. However as we will show in the simulation studies and also observed in \cite{sarkar2014bayesian}, a parametric distribution on $X$ is not capable of recovering the unknown infinite dimensional parameters $f$. As a flexible prior distribution on the density $p$, we propose to use a Dirichlet process Gaussian mixture prior defined by
\begin{eqnarray} \label{eq:DPMM}
X \sim g(\cdot), \quad g(\cdot) = \int \phi_{\sqrt{\tau}}(\cdot - \mu)\, G(d\mu, d\tau), \quad G \sim \textsc{dp}(\alpha G_0).
\end{eqnarray}
Here $\phi_{\sqrt{\tau}}(\cdot - \mu)$ denotes the normal density function with mean $\mu$ and variance $\tau$. 
$\textsc{dp}(\alpha G_0)$ denotes a Dirichlet process prior \citep{ferguson1973bayesian} with $G_0$ as the base probability measure on $\mathbb{R} \times \mathbb{R}^+$ and $\alpha > 0$ is a precision parameter. Given a probability space $\m P$, for any $P \in \m P$ we define the measure space $(\m X, \Omega, P)$ with $\Omega$ denoting the Borel sets of $\m X$, a Dirichlet process satisfies that for any finite and measurable partition $B_1, \dots, B_k$ on $\m X$, $\{P(B_1),\dots,P(B_k)\} \sim \mbox{Dir}\{\alpha G_0(B_1),\dots, \alpha G_0(B_k)\}$, where $\mbox{Dir}\{a_1,\ldots,a_k\}$ denotes the Dirichlet distribution with parameters $a_1,\ldots,a_k$.  A Dirichlet process Gaussian mixture prior is known to be a highly flexible nonparametric prior on the space of densities having a common support as the base measure $G_0$ \citep{escobar1995bayesian}.  It has thus become a very popular Bayesian density estimation method which received considerable attention over the last two decades both from computational \citep{kalli2011slice,neal2000markov} and theoretical perspectives \citep{ghosal2007posterior,kruijer2010adaptive,shen2013adaptive}. Recently, Dirichlet process mixture models have also been commonly used for studying the posterior consistency and contraction rate for Bayesian deconvolution problem under various settings \citep{gao2016posterior,su2020nonparametric,rousseau2021wasserstein}.   
In the next section, we shall discuss in detail that applying a Gaussian process prior to recover the true regression combined with modeling the covariate density with a finite approximation of the Dirichlet process Gaussian mixture prior, we can correct for the bias due to the measurement error.

 \section{Posterior Contraction Properties} \label{theory}
In this section, we study the frequentist large sample properties of the proposed Bayesian errors-in-variables model. We begin with a description of notations used throughout the rest of the paper in Section~\ref{sec:notation}, then state assumptions on the true functions and priors in Section~\ref{sec:assfp}. Section~\ref{sec:main_result} contains our main result on the posterior contraction rate.

\subsection{Notation and Preliminaries}\label{sec:notation}

Let $\floorx$ denote the greatest integer that is strictly less than or equal to $x$ for all $x \in \mathbb{R}$. We define the $L_1$ norm as $\|f\|_1 = \int |f(x)| dx$, and define the supremum norm as $\|f\|_{\infty} = \sup_{x\in S}|f(x)|$, where $S$ is the domain of function $f$. We say a sequence of measures $P_n$  converges weakly to a measure $P$, denoted by $P_n \rightsquigarrow P$  if $\int \phi d P_n \to \int \phi d P$, for all bounded continuous function $\phi$. Denote by $\mathcal{C}[0,1]$ the space of continuous functions defined on $[0,1]$ and denote by $\mathcal{C}^{\beta}[0,1]$ the H$\ddot{\text{o}}$lder space of $\beta$-smooth functions $f : [0,1] \rightarrow \mathbb{R}$ satisfying
\begin{eqnarray*}
| f(x+y)^{\floorbeta}-f(x)^{\floorbeta}| \leq L|y|^{\beta-\floorbeta}, \quad (x,y)\in [0,1],
\end{eqnarray*}
for some constant $L>0$.  For any probability measure $F$ on $\mathbb{R}$, let $p_{F, \sigma}(x) = \int \phi_{\sigma}(x-z)dF(z)$ be the location mixture of normals induced by $F$. For any finite positive measure $\alpha$ write $\bar{\alpha}=\alpha/ \alpha(\mathbb{R})$, where $\alpha(\mathbb{R})$ denotes a measure on $\mathbb{R}$. Let $\textsc{dp}(\alpha)$ denote the Dirichlet process with the base measure $\alpha$. We denote the prior distribution by $\Pi(\cdot)$ and the posterior distribution by $\Pi_n (\cdot \mid D_n)$. 
For two positive sequences $a_n, b_n$, we write $a_n\asymp b_n$ if $a_n/b_n$ can be bounded from below and above by finite constants. In addition, we use ``$\lesssim$" (``$\gtrsim$") to indicate inequalities up to finite universal constants.

\subsection{Assumptions}\label{sec:assfp}
\begin{ass}\label{ass:f0}
The regression function $f_0 \in \mathcal{C}^{\beta}[0,1]$ with  $\beta > 1/2$. We also assume $\|f_0\|_{\infty} <A_0$ for some large enough constant $A_0$.
\end{ass}
We assume that $\beta$ is unknown while fitting the model and our optimal convergence rate results are adaptive for any choice of $\beta > 1/2$. This is achieved easily in a Bayesian paradigm through a suitable prior on the smoothness parameter of the Gaussian process. The finite upper bound assumption is common to achieve the adaptivity in the errors-in-variables problem, similar assumptions can be found in \cite{chesneau2010adaptive,chichignoud2017adaptive}. In practice, we can obtain a reasonable upper bound as a multiple of averaged responses from additional validation data sets \citep{yang2016bayesian}. 
The lower bound on the smoothness is also a common assumption in the random design regression problem, refer to \cite{baraud2002model,birge1979theoreme,brown2002asymptotic} for further discussion on this topic.  
%

\begin{ass}\label{ass:p0}
 The marginal density $p_0$ of the unobserved covariates $X$ is in $\mathcal{C}^{\beta'}[0,1]$ for some $\beta' \ge \beta$, where $\beta$ is defined in Assumption \ref{ass:f0}. Also, we assume there exists a finite constant $B>0$ such that $\inf_{x\in [0,1]} p_0(x) \geq B^{-1}$.
\end{ass}
Smoothness assumptions and the lower bound assumption on the marginal density ensure a better control of the numerator and the denominator of the deconvolution kernel estimator defined in Equation \eqref{eq: fn} separately. Analogous smoothness assumptions can be found in \cite{fan1993nonparametric}, that the regression function and marginal density are assumed to have the same smoothness level.  Refer also to \cite{delaigle2007nonparametric} where $f_0\,p_0$ and $p_0$ are assumed to have same regularity level.  

 The assumption $\beta' > \beta$ in Assumption \ref{ass:p0} requires discussion.  From model \eqref{eq:model}, the deconvolution density estimation problem for $p_0$ can be reduced to a random design regression function estimation problem for $f_0$ by conditioning  on a density $p$ in the parameter space. Hence the overall convergence rate will be determined by the slowest contraction rates for estimating $p_0$ and $f_0$. Although our theory is derived for compactly supported $p_0$, it can be extended to the unbounded support case under desirable tail conditions \citep{kruijer2010adaptive} on $p_0$.



In the Bayesian errors-in-variables model defined in Equation \eqref{eq:Bmodel}, we assign a centered and rescaled Gaussian process prior on $f$, denoted by $\textsc{gp}(0, c; A)$, associated with the squared exponential covariance kernel $c(x, x'; A)=\exp\{-A^2 \|x-x'\|^2\} $ with the rescaled random variable $A$ satisfying the following Assumption \ref{ass:A}. This choice is motivated by the fact that a properly scaled squared exponential covariance kernel is known to lead to the optimal rate of posterior convergence \citep{van2007bayesian,van2009adaptive}. In addition, we consider a Dirichlet process Gaussian mixture prior on the marginal density $p$ defined as $p_{F, \widetilde{\sigma}}$, with $F \sim \textsc{dp}(\alpha)$ and $\widetilde{\sigma} \sim G$, where $G$ satisfies Assumption \ref{ass:sig} below. 
\begin{ass}\label{ass:A}
We assume the rescaled parameter $A$ possesses a density $m$ satisfying for sufficiently large $a>0$,
\begin{eqnarray*}
C_1a^p\exp{(-D_1a \log^q a)} \le m(a) \le C_2a^p\exp{(-D_2a \log^q a)},
\end{eqnarray*}
for constants $C_1, C_2, D_1, D_2 >0$ and $p,q \geq 0$. We assume a conditional Gaussian process prior on the sets of all functions $\m{A}=\{f\in \m{C}[0,1]: \|f\|_{\infty} <A_0\}$, for the same constant $A_0$ in Assumption \ref{ass:f0}.
\end{ass}
Assumption \ref{ass:A} includes the gamma density as a special case when $q=0$. A similar assumption appears in \cite{van2009adaptive}. We restrict the Gaussian prior over the set $\m{A}$ based on Assumption \ref{ass:f0}. 
\begin{ass} \label{ass:sig}
The Dirichlet process Gaussian mixture prior on the marginal density $p(x)$ defined by $p_{F, \widetilde{\sigma}}$ with $F \sim \textsc{dp}(\alpha)$ and $\widetilde{\sigma} \sim G$, satisfy the following conditions:
\begin{eqnarray*}
1-\bar{\alpha}[-x,x] \le \exp(-b_1 x^{\tau_1}) \quad \text{for all sufficiently large} \ x>0,& \\
G(\widetilde{\sigma}^{-2} \geq x) \leq c_1\exp(-b_2 x^{\tau_2})\quad \text {for all sufficiently large} \ x>0,&\\
G(\widetilde{\sigma}^{-2} < x ) \leq c_2 x^{\tau_3} \quad \text{for all sufficiently small} \ x >0,& \\
G(s<\widetilde{\sigma}^{-2}<s(1+t)) \leq c_3 s^{c_4} t^{c_5}\exp(-b_3 s ^{1/2}) \quad \text {for}\ s>0 \ \text{and} \ t \in (0,1),&
\end{eqnarray*}
for positive constants $\tau_1, \tau_2, \tau_3, b_1, b_2, b_3, c_1, \ldots, c_5$.
\end{ass}
The inverse-gamma density on $\widetilde{\sigma}$ satisfies the above assumptions, whereas the inverse-gamma density on $\widetilde{\sigma}^2$ does not.  This is a fairly standard assumption in the Bayesian asymptotics literature on the Dirichlet process mixture of Gaussians, for similar assumptions refer to the posterior convergence analysis for density estimation in \cite{shen2013adaptive}.

\subsection{Main Theorem on Posterior Contraction}\label{sec:main_result}

For the model defined in Equation \eqref{eq:Bmodel}, we first define the marginal likelihood of random pairs $\{(Y_i, W_i), i=1,\ldots,n \}$ as $g_{f,p}(y,w) =(2\pi \delta)^{-1} \int \phi_1\{y-f(x)\}\, \phi_{\delta}(w-x)\, p(x)\, dx$ and denote its distribution measure by $G_{f,p}$. Recall that we assume the noise level of the random error $\sigma = 1$ to simplify the calculation.  Based on the Baye's rule, the posterior distribution given $n$ pairs of observations denoted by $\{Y_{1:n},W_{1:n}\}$ can be written as
\begin{eqnarray}\label{df_post}
\Pi_n\{(f,p)\in B \mid Y_{1:n}, W_{1:n}\} = \frac{\int_{B}\Pi_{j=1}^n g_{f,p}(Y_j, W_j) d\Pi(f)d\Pi(p)}{\int_\m P \Pi_{j=1}^n g_{f,p}(Y_j, W_j) d\Pi(f)d\Pi(p)},
\end{eqnarray}
where $B$ is any measurable subset of $\m P = \{(f, p): f \in \m C [0,1]; \ p: [0, 1] \to \mathbb{R}, \,\int p(x)dx=1 \}$.  

\begin{thm} \label{thm:mthm}
Suppose $f_0$ and $p_0$ satisfy Assumptions \ref{ass:f0} and \ref{ass:p0}, respectively, and the prior $\Pi$ on $(f,p)$ satisfies Assumptions \ref{ass:A} and \ref{ass:sig}.  Then for some sufficiently large constant $M>0$ 
and the standard deviation $\delta_n$ of normal measurement error, 
\begin{equation*}
\begin{aligned}
\Pi_n\{ (f,p): \|f-f_0\|_1 < M\,{\tcr \max(\epsilon_n,\delta_n^{\beta})}, \|p-p_0\|_1 < M\,{\tcr\max(\epsilon'_n, \delta_n^{\beta'})}  \mid Y_{1:n}, W_{1:n}\}\\
\to 1\ \text{almost surely in}\ G_{f_0,p_0}, \ \text{as} \ n \to \infty,
\end{aligned}
\end{equation*}
where {\tcr $\epsilon_n = n^{-\beta/(2\beta+1)}(\log n)^t$ and $\epsilon'_n = n^{-\beta'/(2\beta'+1)}(\log n)^{t}$ with $t=\max{\{(2\vee q)\beta/(2\beta+1),t'\}}$, where $t'>(\gamma+1/\beta')/(2+1/\beta')$ for some $\gamma>2$.} When $\delta_n \lesssim \epsilon_n^{1/\beta}$, {\tcr the convergence rate for recovering $f_0$ is a multiple of $\epsilon_n$.}
\end{thm}
It has been known that fixing $\delta_n \equiv 1$ leads to a logarithmic minimax error rate for errors-in-variables regression estimation. We remark that Theorem \ref{thm:mthm} does not yield the optimal rate in this case, as the current method to deliver posterior contraction rate for nonparametric models is sharp only up to logarithmic terms. However, when $\delta_n \lesssim n^{-1/(2\beta+1)}$, Theorem \ref{thm:mthm} shows {\tcr that optimal rates for regression and density estimation under the EIV setting is the same under the regular nonparametric setting without measurement errors, respectively. } 

The proof of Theorem \ref{thm:mthm} can be found in Appendix~\ref{app:mthm}. Existing contraction rate results in the frequentist deconvolution literature \citep{fan1993nonparametric} require the knowledge of the smoothness of both the true covariate density and the regression function to achieve the optimal convergence rate for the regression function.  Theorem \ref{thm:mthm}, on the other hand, achieves minimax optimal rate of posterior convergence adaptively over all smoothness levels $(\beta', \beta)$ for $\beta,\beta'$ defined in Assumptions \ref{ass:f0} and \ref{ass:p0}, given the knowledge of decaying rate of the error standard deviation $\delta_n$ (or the number of replications). To understand the implication of the posterior convergence rate of $f$ in Theorem \ref{thm:mthm} let us focus on the case where $\beta =1$. Since $\{f(X) - f_0(X)\} \asymp \{f(W) - f_0(W)\} + \{f'(W) + f_0'(W)\} (X-W)$, the convergence rate for estimating $f$ is limited by how fast the marginal density of $X$ can be recovered from observations $W$. This intuitively justifies the rate $\max(\epsilon_n,\epsilon'_n,\delta_n^{\beta})$ in EIV model.  {\tcr We remark that the rate results in Theorem \ref{thm:mthm} also hold for $f,p$ even if $\beta>\beta'$. In that case, the posterior of $p$ always attains the near-minimax rate in recovering the true density, whilst the best obtainable posterior rate for recovering $f_0$ is limited to $\epsilon'_n$, which is slower than $\epsilon_n$.  }

Analyzing the posterior distribution following the seminal work \citep{ghosal2000} requires upper-bounding the numerator of the posterior defined in Equation \eqref{df_post} over some set $B$ of interest and lower-bounding the marginal likelihood. In our proof, the numerator can be bounded above by constructing a sequence of test functions that is used to test the true model against models outside a small neighborhood of the truth under proper metric. As a key technical contribution, we obtain sharp bounds for Type I and Type II errors of the constructed tests by developing large deviation bounds for the deconvolution density estimator, which generalizes some results in \cite{JMLR:v16:pati15a} for random design regression to errors-in-variables problem.  To bound the marginal likelihood from below, it requires the priors assigned on the regression and covariate density assigning enough mass around the truth.  A component-wise Gaussian prior on the covariate cannot concentrate enough over a small neighborhood of the true locations, simply because the concentration of $n$-dimensional standard Gaussian vector cannot exploit the smoothness of the density and hence cannot assign enough mass within a small neighborhood around the true density.  On the other hand,  a mixture of normals prior allows borrowing of information, naturally exploits the smoothness and provides adequate concentration. A similar treatment to the covariate density can also be found in recent Bayesian deconvolution literature \citep{gao2016posterior,donnet2018posterior,rousseau2021wasserstein}.

\section{Posterior Computation}
In order to sample from the posterior distribution of $(f, p, \delta, \sigma)$, we employ a Gibbs
sampler and sample from each of the parameters given the others.  Posterior sampling methods for Bayesian density estimation using Dirichlet process Gaussian mixture prior is popular, refer to the P\'{o}lya urn sampler \citep{escobar1995bayesian,maceachern1998estimating} and blocked Gibbs sampler with stick-breaking representation \citep{ishwaran2001gibbs}.  In this article, we use the finite approximation of the Dirichlet process Gaussian mixture prior with the stick-breaking representation.  The major bottleneck of the computation stems from sampling the Gaussian process term $f$ which requires a) inversion of $n \times n$ matrices depending on the latent covariates and b) sampling from the conditional distribution of the true covariates, which is intractable.    Task a) makes 
the algorithm computationally inefficient and unstable specifically for the errors-in-variables regression problem, since it requires evaluating the inverse of the covariance matrix repeatedly along with the updates of covariates. To bypass $O(n^3)$ computation steps associated with inverting an unstructured $n \times n$ covariance matrix, numerous powerful techniques have been proposed in the last decade; fixed rank kriging \citep{banerjee2008gaussian,finley2009improving}, covariance tapering \citep{furrer2006covariance,kaufman2008covariance}, composite likelihood methods \citep{guan2006composite,heagerty1998composite}. In using these techniques,  often the original covariance kernel itself is not preserved, which means the covariance function of the approximate process differs from the covariance function of the original process.  More recently, \cite{stroud2017bayesian} and \cite{guinness2017circulant} derived a fast algorithm of sampling from stationary Gaussian processes on the large-scale lattice data, using the circulant embedding technique proposed in \cite{wood1994simulation}. Such techniques typically require the assumption of equally spaced covariates. In the absence of equally spaced design, the idea is to define a larger lattice and consider the prediction as missing data imputation \citep{guinness2017circulant,stroud2017bayesian}.  However, it is not straightforward to translate these ideas to the errors-in-variables regression problem as the true covariates are contaminated and the true marginal distribution remains unknown. Instead, we propose using a lower dimensional mapping to approximate the Gaussian process based on the random Fourier basis proposed by \cite{rahimi2008random}. And the random mapping to the Fourier domain preserves the covariance kernel associated with the original Gaussian process. This also avoids computing the inverse of covariance matrix by introducing moderate numbers of parameters associated with the Fourier basis. Moreover, this is suitable in applications where practitioners have a pre-conceived notion of using a particular covariance function and we require the approximated covariance to accurately reflect that prior opinion. The lower dimensional mapping is chosen to approximate the original Gaussian process arbitrarily well; refer to Theorem \ref{thm:appgp}.  We describe the approximate Gaussian process in the following Section~\ref{GPEV}. 

\subsection{An Approximation of the Gaussian Process} \label{GPEV}
 {\tcr The low-rank projection of a stationary covariance kernel on a random feature space is a popular approach to scale up kernel-based regression methods \citep{rahimi2008random}.  Theoretical properties of the random Fourier feature projection have been extensively studied in the last decade, mostly in terms of the approximation accuracy of the covariance kernel function \citep{sutherland2015error,sriperumbudur2015optimal}, properties of the induced RKHS \citep{rahimi2008uniform,rahimi2008weighted,bach2017equivalence}, and the expected risk bounds of an approximated kernel ridge regression estimator based on random Fourier features and their variants \citep{avron2017random,li2019towards,zhang2019low,yang2021exact}. For a detailed and categorized summary of existing results, one may refer to a recent work \citep{liu2021random}.} In this section, we develop a low-rank random Fourier basis projection as an approximate of a stationary zero-mean Gaussian process $\textsc{gp}(0, c)$, {\tcr which can be represented as a Bayesian linear model. Such a representation has been considered in \cite{wilson2020efficiently} where they used a random feature projection to approximate the original GP, and further approximated the obtained posterior distribution to speed up posterior computation. In our case we study the exact posterior distribution resulting from the approximated GP prior. }

Denote by the corresponding spectral density $\phi_c(\cdot)$ defined through $c(h) = \int e^{ihx} \phi_c(x) dx$.  For a suitably chosen large integer $N$, we define 
\begin{eqnarray}\label{eq:rfbp}
\widetilde{f}_N(x)= (2/N)^{1/2}\sum_{j=1}^{N} a_j \cos (w_j x + s_j),
\end{eqnarray}
where $a_j \sim \mbox{N}(0,1)$, $w_j \sim \phi_c$ and $s_j \sim \text{Unif}\,[0, 2\pi]$, i.i.d. for $j = 1,\ldots,N$.  
{The random process \tcr $\widetilde{f}_N(x)$ is an $N$-dimensional approximation to a GP such that its covariance function coincides with the kernel function of original GP. In addition, Theorem \ref{thm:appgp} shows that the approximate also converges to the original Gaussian process $\textsc{gp}(0, c)$ weakly. } 

\begin{thm}\label{thm:appgp}
Suppose $f$ is the original Gaussian process $\textsc{gp}(0, c)$ and $\widetilde{f}_N$ is defined in Equation \eqref{eq:rfbp}, we have 
\begin{eqnarray*}
\widetilde{f}_N \rightsquigarrow f, \quad \mbox{as} \quad N\to\infty. 
\end{eqnarray*}
Also, for any $x, y \in \mathbb{R}$, 
\begin{eqnarray*}
\mbox{E}\,\{\widetilde{f}_N(x)\}=0; \quad \Cov\{\widetilde{f}_N(x),\widetilde{f}_N(y)\}=c(x,y).
\end{eqnarray*}
\end{thm}


The proof of Theorem \ref{thm:appgp} is deferred to Appendix \ref{sec:appgp}. The construction $\widetilde{f}_N$ is related to the random feature mapping in the Fourier domain \citep{rahimi2008random}, used to project the kernel onto a lower-dimension space $\mathbb{R}^N$. It is straightforward to show that preservation of the covariance kernel associated with the original Gaussian process for the proposed process defined over the real area, due to the expression of Fourier features. However, the weak convergence result of $\widetilde{f}_N$ is non-trivial and the proof provides a framework to study the asymptotic property of random processes constructed based on Fourier projection.  

{\tcr 

Theorem \ref{thm:appgp} validates the usage of $\wt f_N$ to approximate a stationary GP in an asymptotic manner.  
Allowing $\wt f_N$ to be adaptive to the unknown smooth level of the true regression function, we now assume $\{\omega_j\}$ are independently and identically generated from the spectral measure of a rescaled squared exponential kernel function. The rescaling parameter is unknown and endowed with the prior satisfying Assumption \ref{ass:A}.  
Adopting the same notation of a rescaled mean-zero Gaussian process  $f \sim \textsc{gp}(0, c^{A})$ considered in Section \ref{theory} and the rescaling $A\sim g(\cdot)$ for some distribution $g$, analogously, we define a rescaled version of \eqref{eq:rfbp} as 
\begin{eqnarray*}
\wt{f}_{A,N}(x)= (2/N)^{1/2}\sum_{j=1}^{N} a_j \cos (w^{A}_j\ x + s_j),
\end{eqnarray*}  
where $\{a_j, s_j\}$ are same as in \eqref{eq:rfbp}, and for any $a>0$, $w^{A}_j | (A=a) \overset{i.i.d.}\sim \phi^{a}_c$ for $j=1,\ldots,N$, recall that $\phi^{a}_c(\lambda)=a^{-1}\phi_c(\lambda/a)$ denotes the spectral density of GP associated with a squared exponential kernel function indexed with the rescaling parameter $a$. It is straightforward to show that Theorem \ref{thm:appgp} holds with $A=a$ for any fixed $a>0$. In addition, Theorem \ref{thm:rf} below verifies that the posterior of $\wt{f}_{A,N}$ converges towards the true regression curve $f_0$ at a near minimax rate in the EIV regression problem as well, given an appropriate number of the random Fourier features. 

\begin{thm}\label{thm:rf}
Suppose $f_0$ and $p_0$ satisfy Assumptions \ref{ass:f0} and \ref{ass:p0}, respectively, and Assumptions \ref{ass:A} and \ref{ass:sig} hold.  Then for some fixed large constant $M>0$ and for the number of features $N$ satisfying $N\asymp n\epsilon_n^2 (\log n)^{\tilde{t}}$ for some constant $\tilde{t}>0$ which is free of $n,N$, 
and recall the standard deviation $\delta_n$ of normal measurement error, 
\begin{equation*}
\begin{aligned}
\Pi_n \left\{ (\wt{f}_{A,N},p): \|\wt f_{A,N}-f_0\|_1 < M\max(\epsilon_n, \delta_n^{\beta}), \|p-p_0\|_1 < M\max(\epsilon'_n, \delta_n^{\beta'})  \mid Y_{1:n}, W_{1:n}\right\}\\
\to 1\ \text{almost surely in}\ G_{f_0,p_0}, \ \text{as} \ n \to \infty,
\end{aligned}
\end{equation*}
where $\epsilon_n, \epsilon'_n$ are same as in Theorem \ref{thm:mthm}. Again, when $\delta_n \lesssim \epsilon_n^{1/\beta}$, the posterior contraction rate of $\wt{f}_{A,N}$ is a multiple of $\epsilon_n$.
\end{thm} 

Theorem \ref{thm:rf} provides an asymptotic result on the posterior distribution of approximated GP, provided the rank of the random feature projection increases at a certain rate with the sample size. 
To the best of our knowledge, this is the first theoretical result on low-rank random feature projection of GPs in $L_1$ norm under a Bayesian framework. This result can be easily adapted to other regression/learning problems beyond the EIV context, such as nonparametric regression with random designs.  The proof of Theorem \ref{thm:rf} is deferred to Appendix, which follows a similar line of arguments as in the proof of Theorem \ref{thm:mthm}. Theorem \ref{thm:rf} delineates a specific increasing rate of the number of random features in order to attain the best rate. A minimum requirement on the number of random features has been determined in literature on KRR with random features, which conveys the idea that larger is the number of the random features, the better is the approximation of the RFF to the original KRR estimator. However, when all $\{a_j, \omega_j, s_j\}$ in \eqref{eq:rfbp} treated as random parameters, the number of random features $N$ cannot increase too fast in order to retain a minimum prior concentration over a small KL-neighborhood of the truth, due to the concentration of measure phenomenon of high-dimensional Gaussian random vectors. 
}

 To implement $\widetilde{f}_{A,N}$, it suffices to treat $\{a_j, w_j,s_j\}$ as unknown parameters endowed with suitable independent priors. More details of the posterior computation is deferred to Appendix \ref{sec:gibbcom}. 
{\tcr  We remark that the order of the the number of the features in Theorem \ref{thm:rf}  is primarily of theoretical interest as the smoothness of the function is unknown. Although we do not have adaptive results, in the empirical study, we find the Gaussian process surrogate performs almost as well as the original Gaussian process when $N$ is chosen within the range $(n/8,  n/2)$ for data sets of moderate sizes. }  

\section{Numerical Results} \label{sec:sims}
In this section, we empirically illustrate applications of the proposed Gaussian process surrogate and its variants to Bayesian errors-in-variables model in the following synthetic examples. 
We consider a uniform marginal distribution $X\sim \text{Unif}\,[-3,3]$ and the regression function:  $f(x) = \sin(\pi x/2)/[1+2x^2\{\text{sign}(x)+1\}]$. We consider three choices of sample size $ n\in\{ 100, 250, 500\}$, and consider additive normal regression errors independently and identically drawn from $\mbox{N}(0, \sigma^2)$ with a fixed noise level $\sigma = 0.2$. We confine ourself to the centered normal distribution $\mbox{N}(0, \delta^2)$ for the measurement error with a sequence of gradually increasing variances $\{\delta^2\}$, for the purpose of checking empirical performance of proposed methods in the presence of measurement errors of varying degrees. Specifically, for $n = 100,250$, we consider $\delta^2 \in \{ 0.01,0.2,0.4,0.6,0.8,1\}$; for $n=500$, we consider $\delta^2 \in \{0.001,0.005,0.01,0.1,0.5,1\}$. Under each setting, we compare the following methods:
\begin{enumerate}
\item \textsc{gpev}$_a$: Approximated Gaussian process model described in Section\ \ref{GPEV} with a Dirichlet process Gaussian mixture prior on the marginal density.
\item \textsc{gpev}$_f$: Full scale Gaussian process model using the predictive formula in Equation \eqref{gp_formula}, with a Dirichlet process Gaussian mixture prior on the marginal density.
\item \textsc{gpev}$_n$: Approximated Gaussian process model described in Section\ \ref{GPEV} with a univariate normal prior on the covariate component-wise.
\item \textsc{gp}: Full scale Gaussian process model that ignores the measurement error.
\item \text{decon}: Deconvolution kernel method from \url{https://github.com/TimothyHyndman/deconvolve}.
\end{enumerate}
To implement \textsc{gpev}$_a$ and \textsc{gpev}$_n$, we consider the following combinations of the sample size $n$ and the number of Fourier basis functions $N$: $(n,N) \in \{(100,40),(250, 60),(500,80)\}$. We remark that the values of $N$ are chosen based on preliminary numerical experiments. We only present the numerical results for $n=100,500$ in this section, the result for $n=250$ is similar and thus deferred to Appendix~\ref{app:addres}. For Bayesian approaches, we ran the Gibbs sampler with 2,000 iterations and discarded the first 1,500 iterations as a burn-in. The derivation of a full conditional and detail on hyperparameter choices can be found in Appendix~\ref{sec:gibbcom}. The investigation on the mixing behavior of the Gibbs sampler for estimated marginal and regression functions as well as other diagnostic checks are deferred to Appendix~\ref{app:addres}.
For the Bayesian methods, the posterior mean denoted by $\widehat{f}$, is our estimator of the unknown regression function $f$ and its pointwise $95\%$ credible interval is obtained by constructing $U(x)$ and $L(x)$ such that
\begin{eqnarray*}
\Pi_n\{ f(x) \in [L(x), U(x)] \mid  D_n\} = 0.95.
\end{eqnarray*}
We also consider simultaneous credible bands centered at the posterior mean $\widehat{f}$ with level $\gamma\in(0,1)$,
\begin{align*}
\mbox{CB}_n(\gamma) = \Big\{f:\, \big\|f - \widehat{f}\big\|_\infty \leq r\Big\},
\end{align*}
where the half length $r$ is chosen so that posterior probability of $f$ falling into the credible band is $\gamma$, 
\begin{align*}
\Pi_n\big\{ f \in \mbox{CB}_n(\gamma)\, \big|\,  D_n \big\}= \gamma.
\end{align*}
Computation of $\mbox{CB}_n(\gamma)$ can be found in Appendix \ref{sec:gibbcom}.

\begin{table}[htbp]
\footnotesize
\begin{center}
    \begin{tabular}{clccccccc} \toprule \footnotesize
     &  &  \multicolumn{6}{c}{$\delta^2$} \\
         \cmidrule(lr){3-8}
   $n$ & Method  & 0$\cdot$01 & 0$\cdot$2 & 0$\cdot$4 & 0$\cdot$6 & 0$\cdot$8 & 1 \\
     \cmidrule(lr){1-8}
  \multirow{5}{*}{$100$}
   &\textsc{gpev}$_a$    &  0$\cdot$58 (0$\cdot$43)  &   1$\cdot$82 (1$\cdot$23)   &  3$\cdot$89 (4$\cdot$00)  &   4$\cdot$64 (3$\cdot$81)   &  5$\cdot$88 (5$\cdot$55)  &  6$\cdot$31 (5$\cdot$02)  \\
    & \textsc{gpev}$_f$     &  0$\cdot$55 (0$\cdot$41)   &  1$\cdot$85 (1$\cdot$22)    &  3$\cdot$20 (2$\cdot$83)  &  4$\cdot$24 (3$\cdot$19)  &  5$\cdot$54 (5$\cdot$21)  &  5$\cdot$82 (4$\cdot$96)  \\
    & \textsc{gpev}$_n$       &  0$\cdot$60 (0$\cdot$45) &  4$\cdot$82 (2$\cdot$35)  &  10$\cdot$98 (4$\cdot$91)  &  15$\cdot$29 (6$\cdot$07)  &  19$\cdot$26 (7$\cdot$78)  &  20$\cdot$98 (9$\cdot$36)  \\
   &  \textsc{gp}       &  3$\cdot$29 (0$\cdot$31)  &  6$\cdot$11 (1$\cdot$39)  &  9$\cdot$35 (2$\cdot$34)  &  12$\cdot$00 (2$\cdot$99) &  14$\cdot$80 (3$\cdot$54)  &  16$\cdot$81 (3$\cdot$94)  \\
 & \text{decon}   &  1$\cdot$18 (1$\cdot$00)   &  5$\cdot$07 (2$\cdot$46)  & 10$\cdot$46 (3$\cdot$79)  &  14$\cdot$72 (4$\cdot$08)  &  18$\cdot$25 (3$\cdot$96)  &  20$\cdot$59 (3$\cdot$52)  \\\toprule 
     & &  \multicolumn{6}{c}{$\delta^2$} \\
         \cmidrule(lr){3-8}
   $n$ &   Method  & 0$\cdot$001 & 0$\cdot$005 & 0$\cdot$01 & 0$\cdot$1 & 0$\cdot$5 & 1 \\
     \cmidrule(lr){1-8}
  \multirow{5}{*}{$500$}
     &\textsc{gpev}$_a$   &0$\cdot$11 (0$\cdot$04)   & 0$\cdot$12 (0$\cdot$04)   &  0$\cdot$13 (0$\cdot$04)  &  0$\cdot$37 (0$\cdot$21)  &  1$\cdot$69 (1$\cdot$20)   &  3$\cdot$35 (3$\cdot$27) \\
  &\textsc{gpev}$_f$   & 0$\cdot$10 (0$\cdot$04)   &  0$\cdot$11 (0$\cdot$04)  &  0$\cdot$12 (0$\cdot$04)  &   0$\cdot$35 (0$\cdot$19)   & 1$\cdot$59 (1$\cdot$04)   &  3$\cdot$94 (5$\cdot$76)  \\
   & \textsc{gpev}$_n$  &0$\cdot$11 (0$\cdot$04)  & 0$\cdot$12 (0$\cdot$05)   &  0$\cdot$14 (0$\cdot$05) &  1$\cdot$51 (0$\cdot$44)   & 12$\cdot$09 (2$\cdot$02)   &  20$\cdot$38 (4$\cdot$12) \\
   & \textsc{gp}       & 1$\cdot$78 (0$\cdot$08)  &   1$\cdot$80 (0$\cdot$09)   &  1$\cdot$80 (0$\cdot$09)   &  2$\cdot$57 (0$\cdot$26)   &  8$\cdot$45 (1$\cdot$08)   &  14$\cdot$37 (1$\cdot$60) \\
   & \text{decon}  & 0$\cdot$35 (0$\cdot$21) & 0$\cdot$38 (0$\cdot$26)   & 0$\cdot$38 (0$\cdot$26)  & 1$\cdot$14 (0$\cdot$46)   & 9$\cdot$48 (1$\cdot$55)   &  18$\cdot$03 (1$\cdot$53) \\ \toprule
\end{tabular}
\end{center}
  \caption{Averaged mean squared errors (\textsc{amse}) defined as $\bbE\, [\, K^{-1}\sum_{k=1}^K \{\widehat{f}(t_k) - f(t_k) \}^2 \,]$ ($\widehat{f}(\cdot)$ denotes the estimator of $f$, $\bbE(\cdot)$ denotes taking average over replicates) on an evenly spaced test grid $\{t_1, \ldots, t_K\}$ of size $K=100$ over $[-3, 3]$ with standard errors ($\times 10^{2}$) over 50 replicated data sets of sizes $n =100, 500$.} 
\label{tab:f1_mse}
\end{table}

\normalsize
\begin{figure}[h!]
\centering 



\begin{multicols}{3}
    \includegraphics[scale = 0.045]{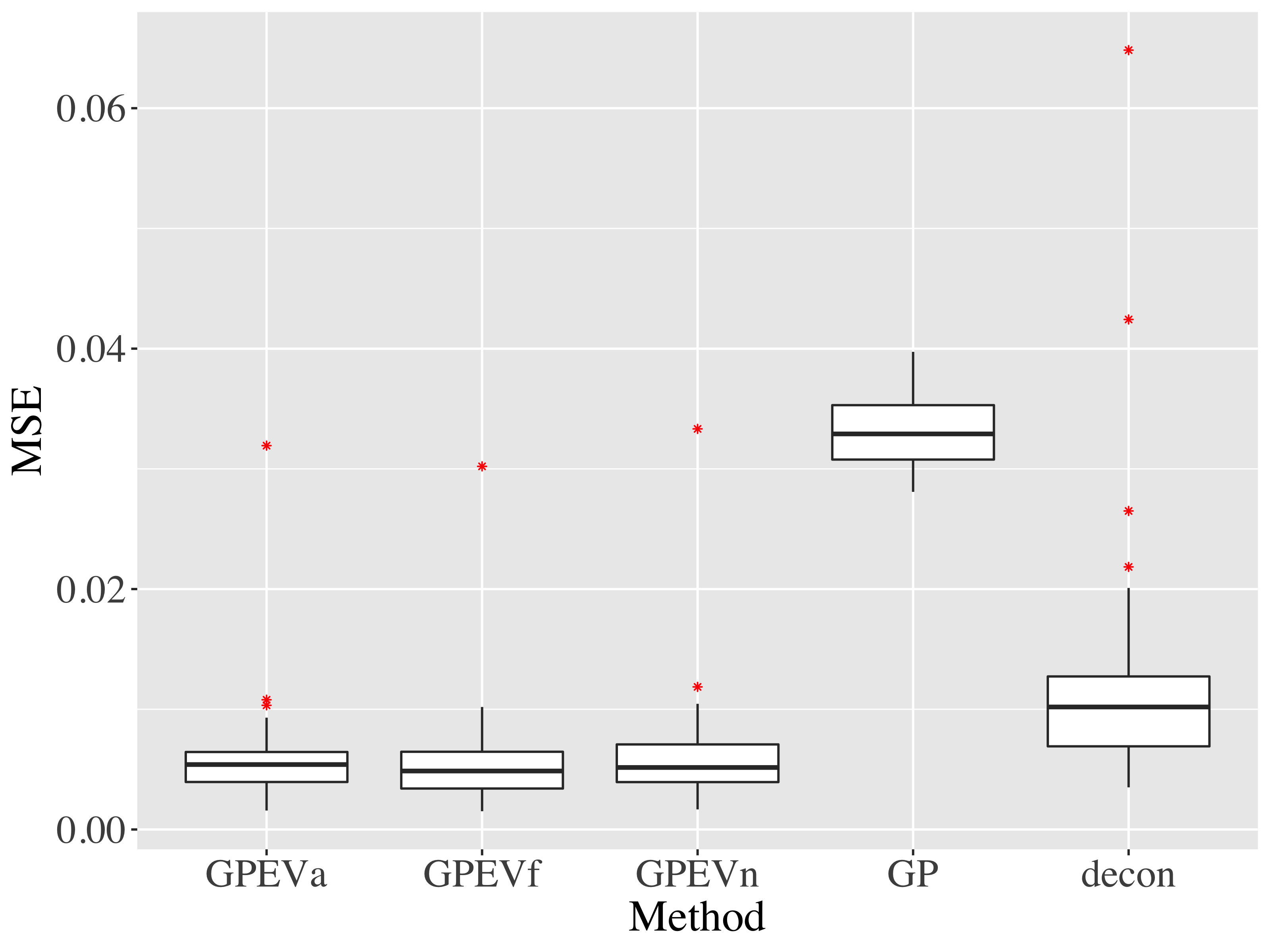}
    \includegraphics[scale = 0.045]{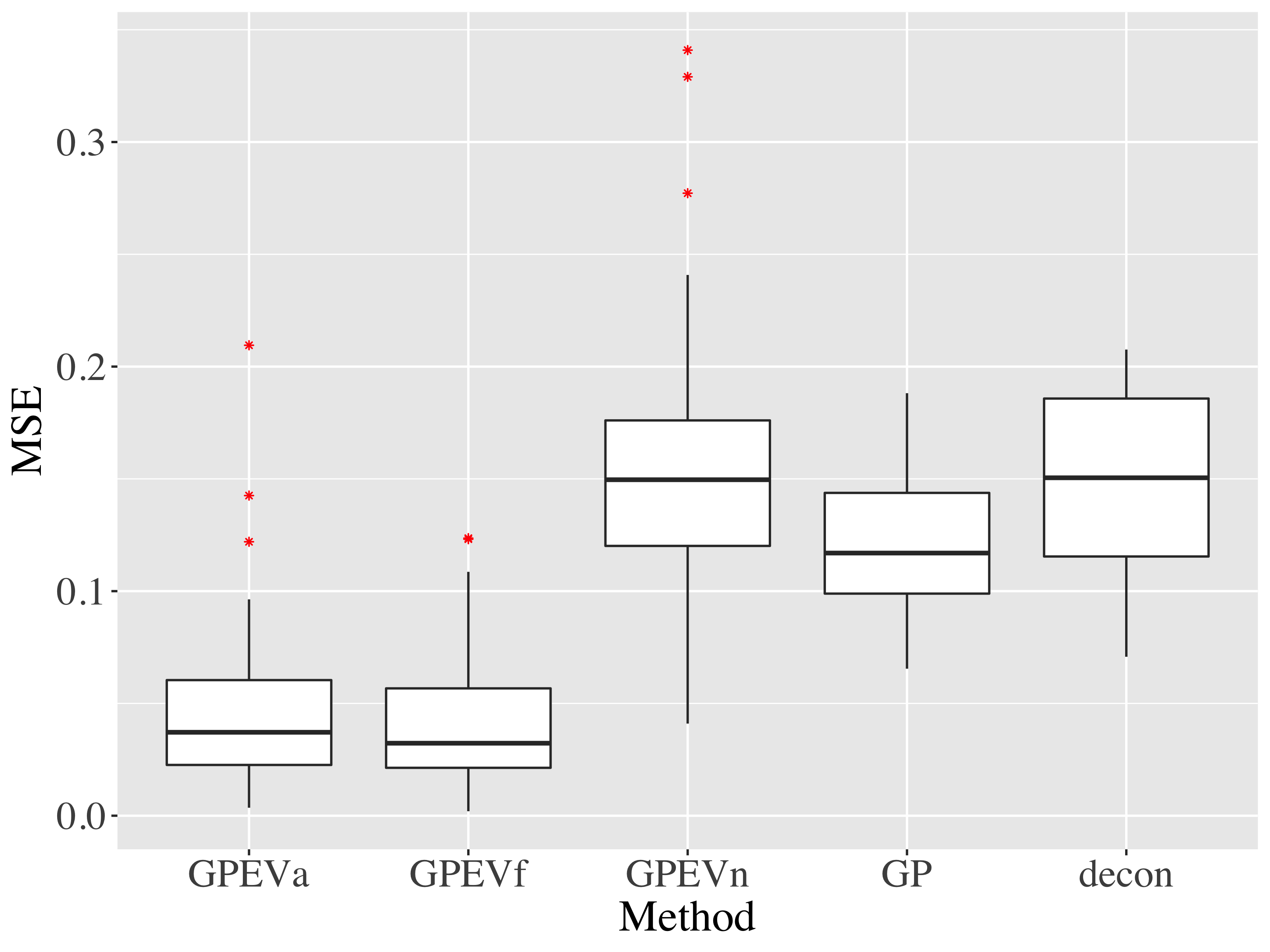}
    \includegraphics[scale = 0.045]{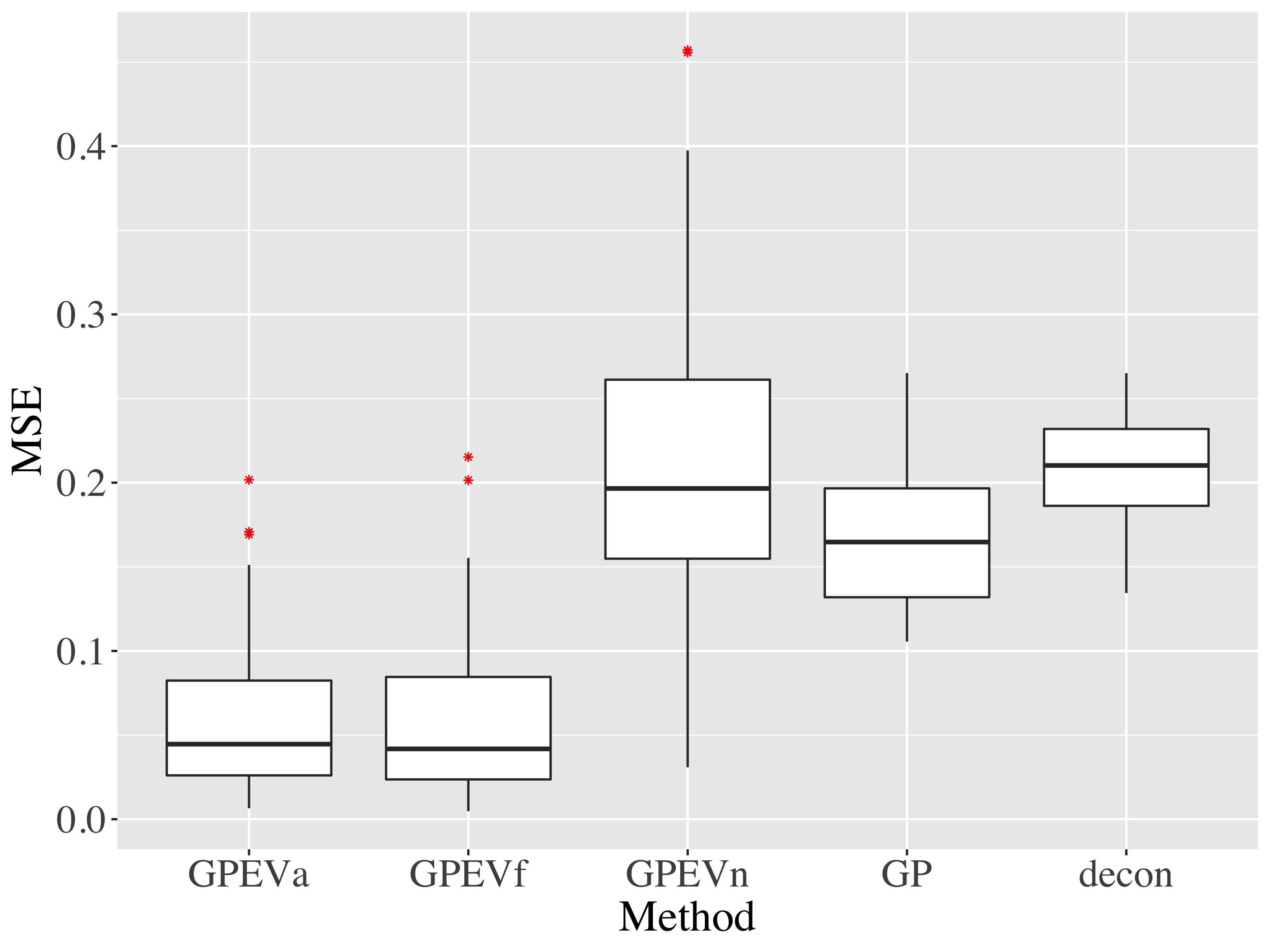}
\end{multicols}
\begin{multicols}{3}
   \includegraphics[scale = 0.045]{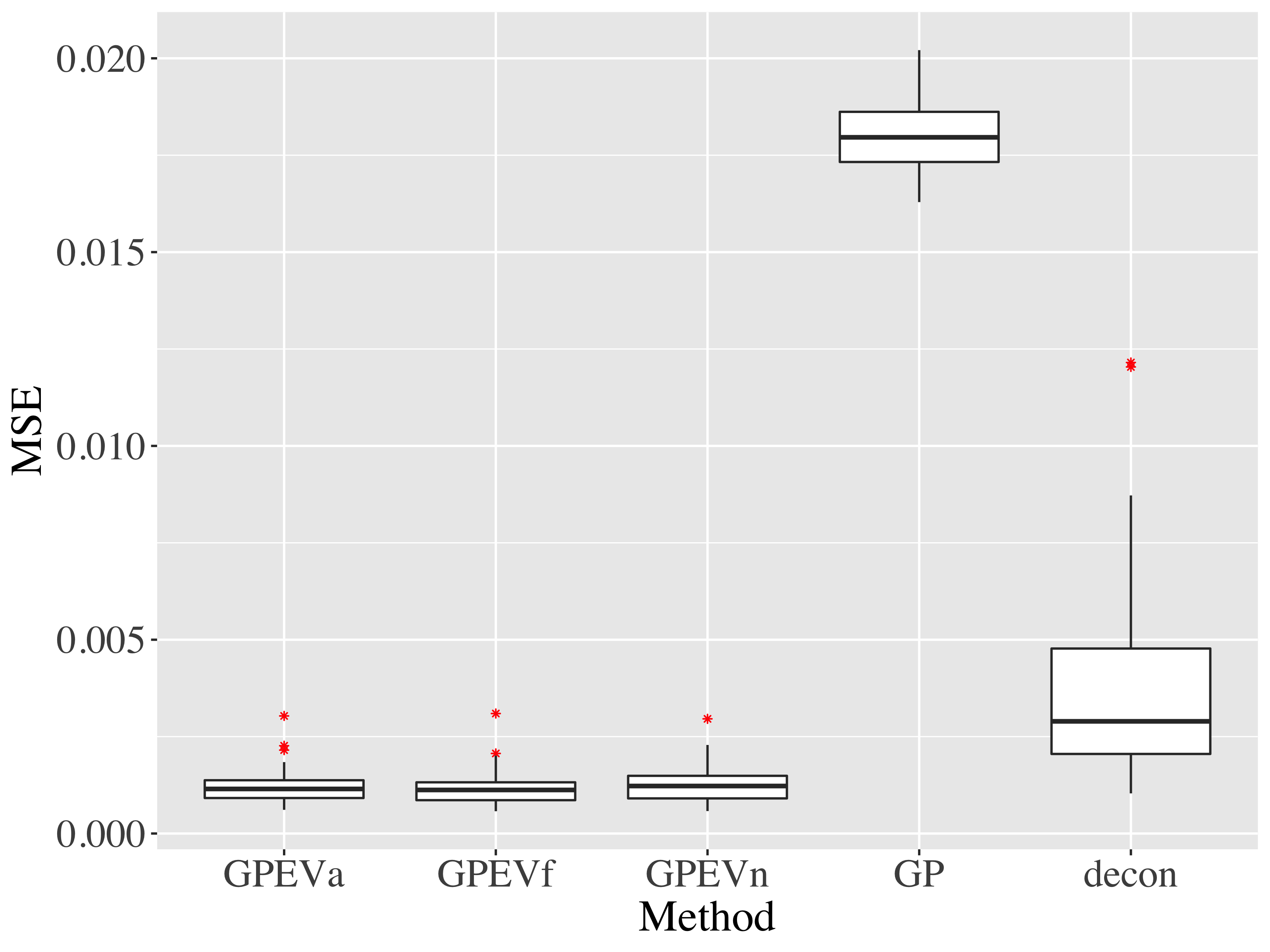}
    \includegraphics[scale = 0.045]{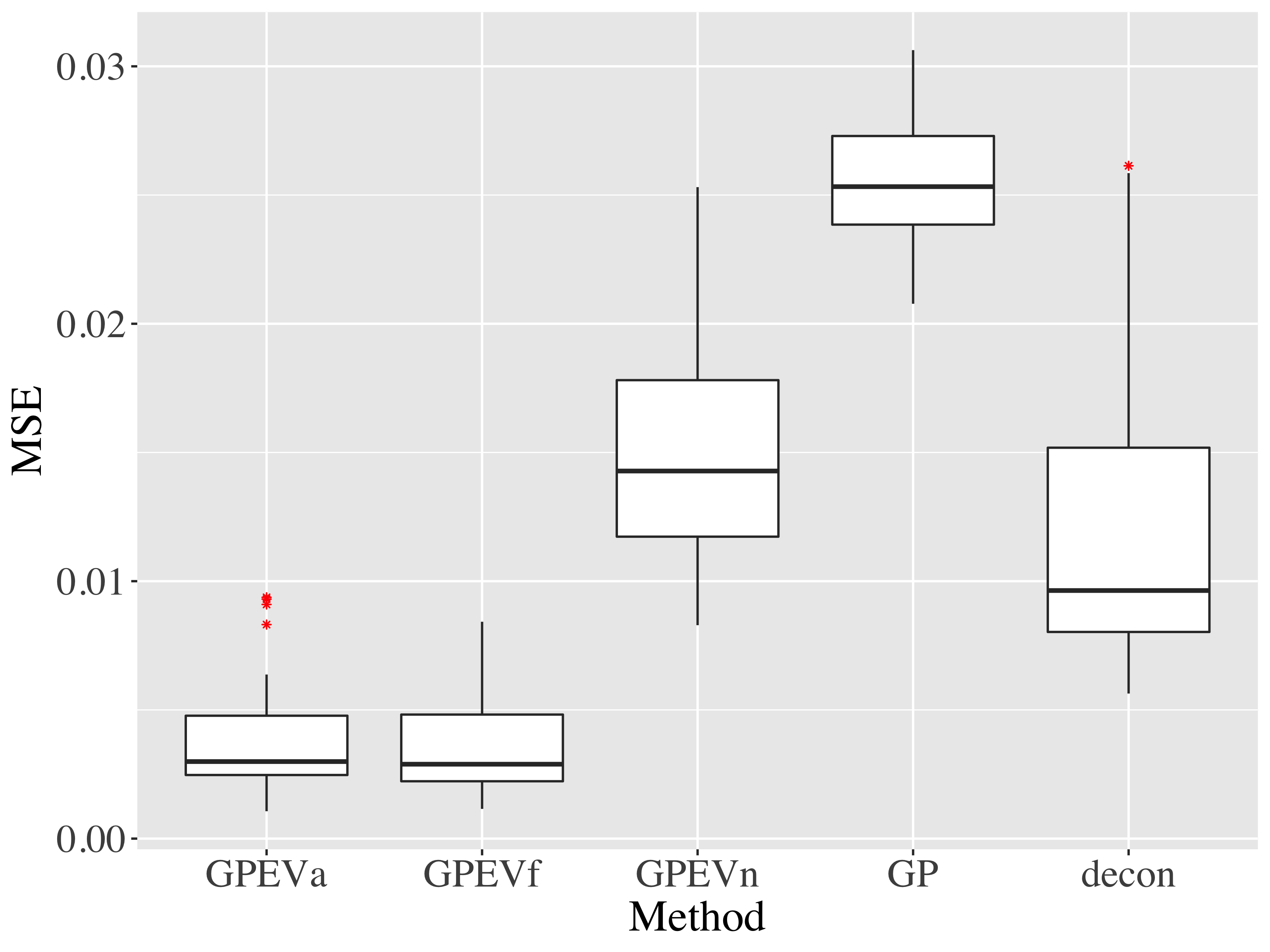}
   \includegraphics[scale = 0.045]{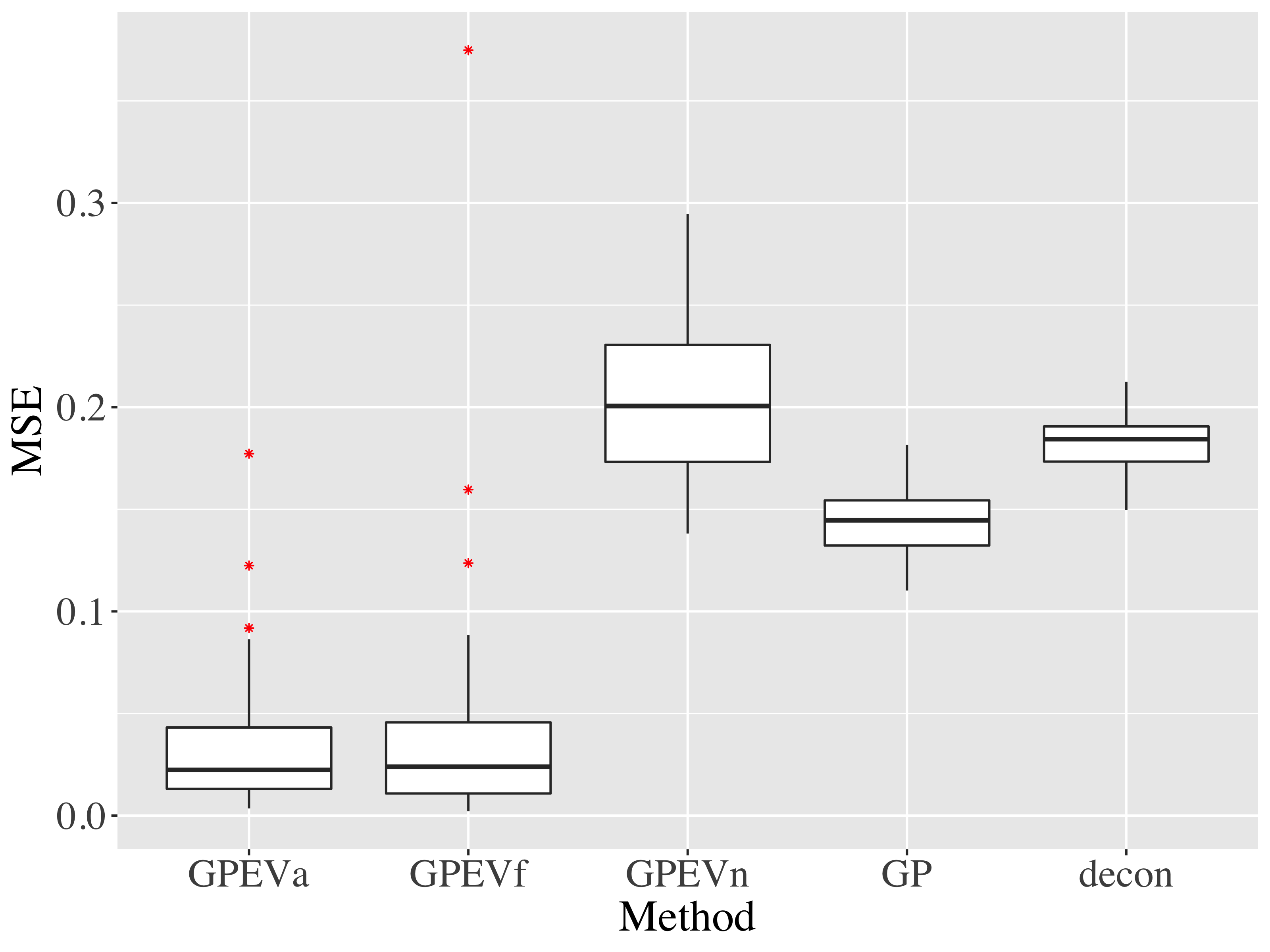}
\end{multicols}
\caption{Boxplots of mean squared values for $f(x)$ when $n =100$ (the top row) and $n=500$ (the bottom row) over 50 replicated data sets. For $n=100$, set $\delta^2 = 0.01$ (left panel), $\delta^2 = 0.6$ (middle panel) and  $\delta^2 = 1$ (right panel); for $n=500$, set $\delta^2 = 0.005$ (left panel), $\delta^2 = 0.1$ (middle panel) and  $\delta^2 = 1$ (right panel). In each panel the methods displayed from left to right are \textsc{gpev}$_a$, \textsc{gpev}$_f$, \textsc{gpev}$_n$, \textsc{gp} and \text{decon}.} \small
\label{fig:boxplot_n500}
\end{figure}

\begin{figure}[h!]
\centering 

\begin{multicols}{3}
    \includegraphics[scale = 0.045]{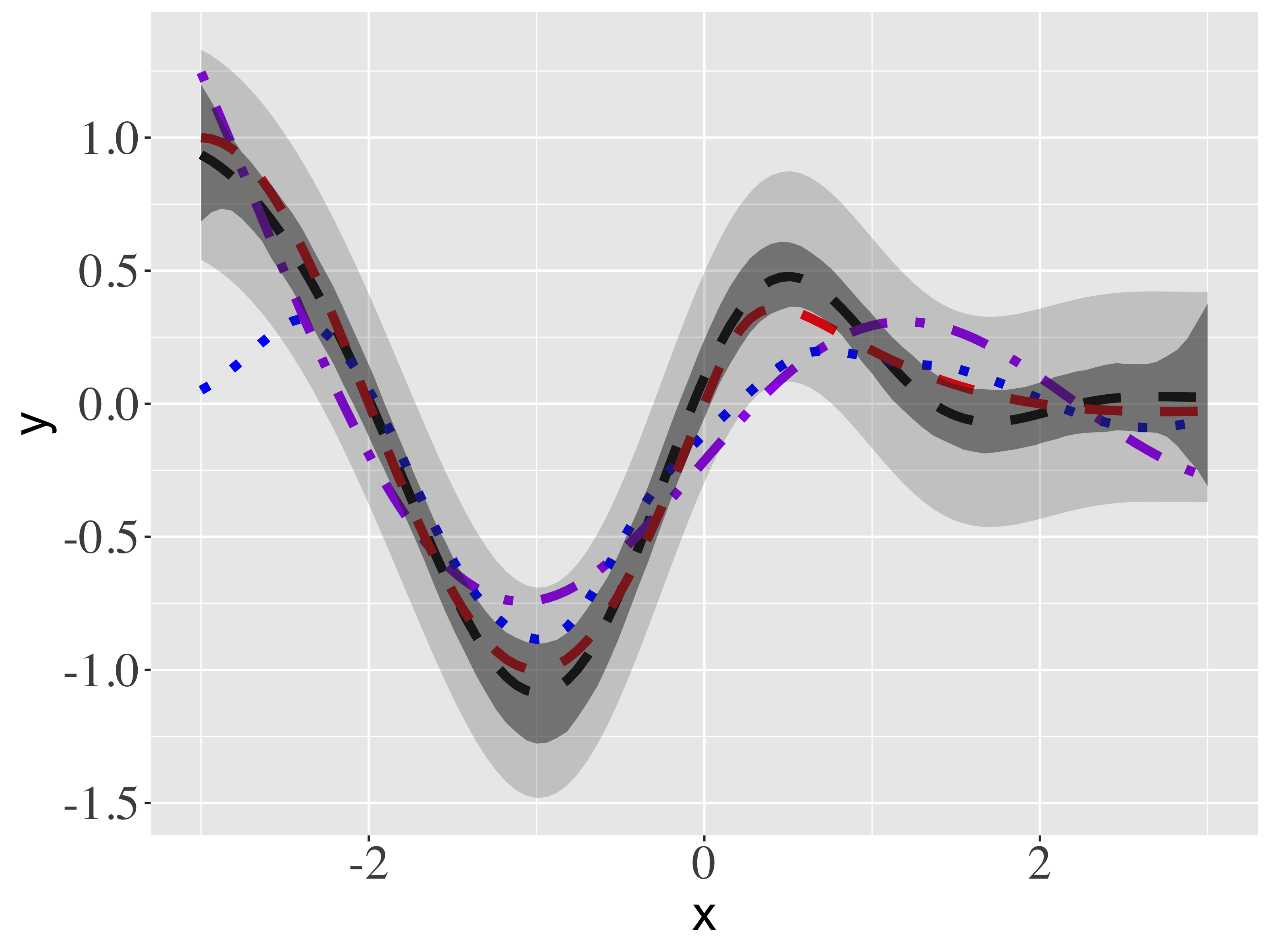}
    \includegraphics[scale = 0.045]{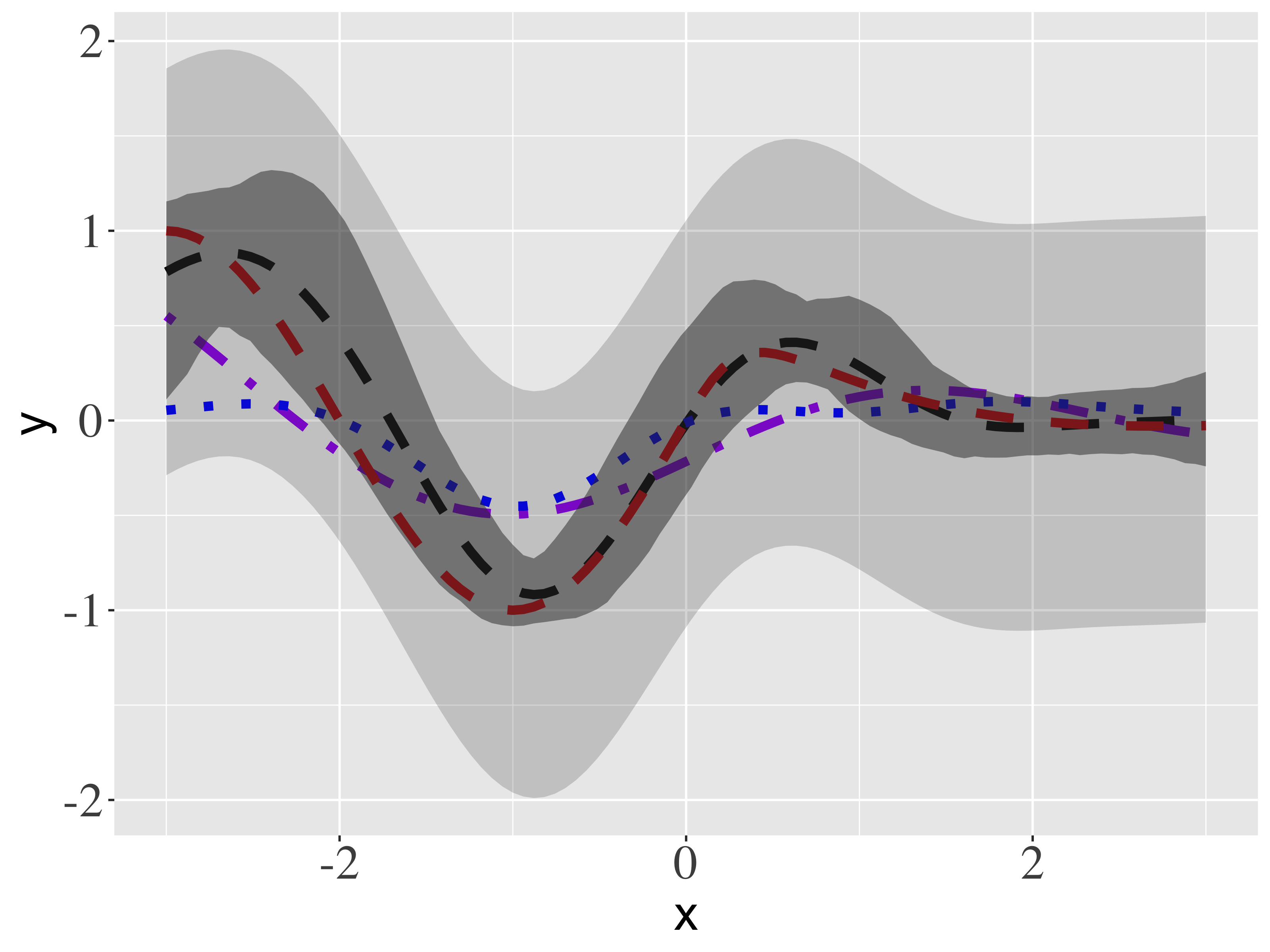}
    \includegraphics[scale = 0.045]{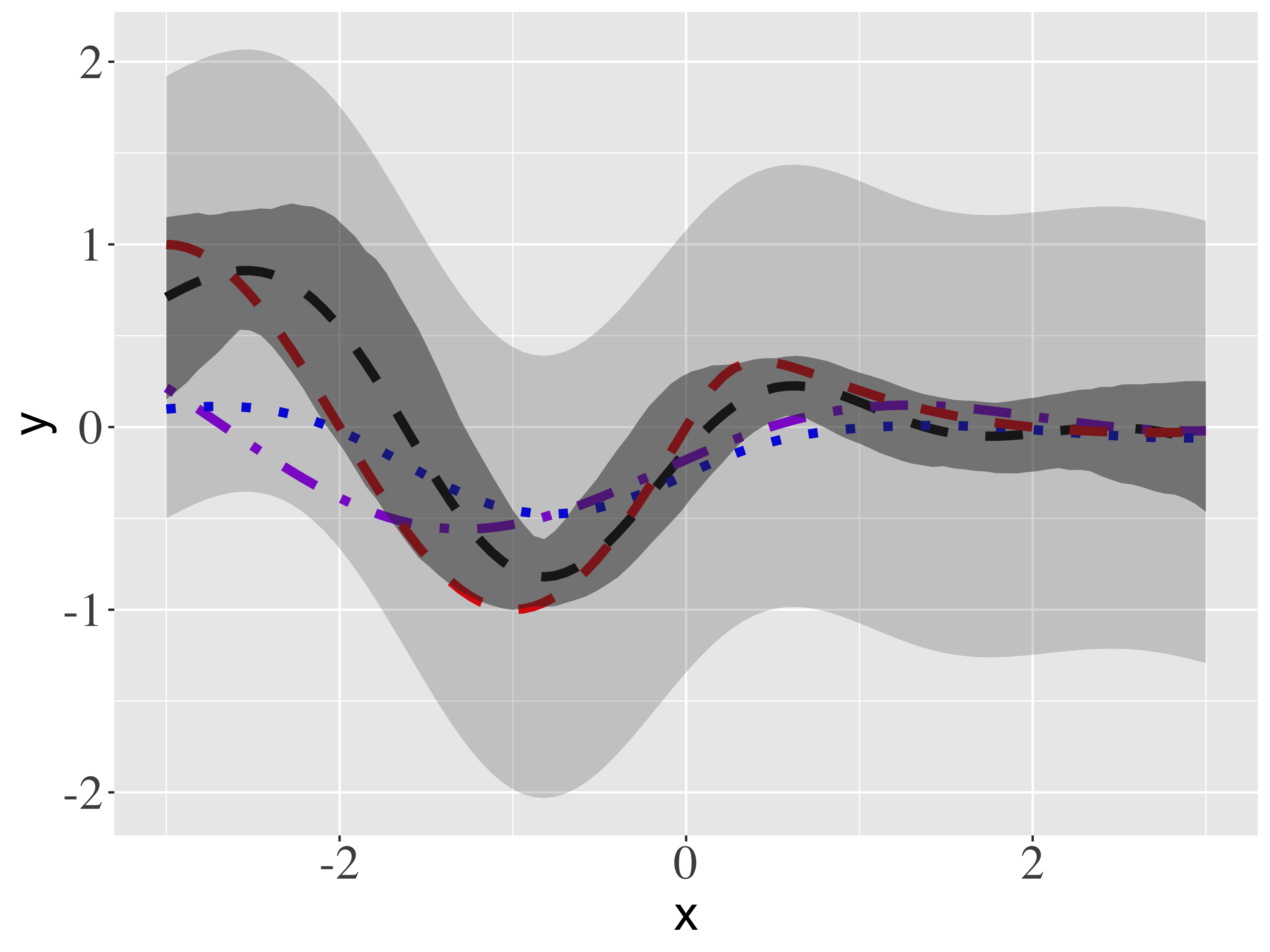}
\end{multicols}
\begin{multicols}{3}

  \includegraphics[scale = 0.045]{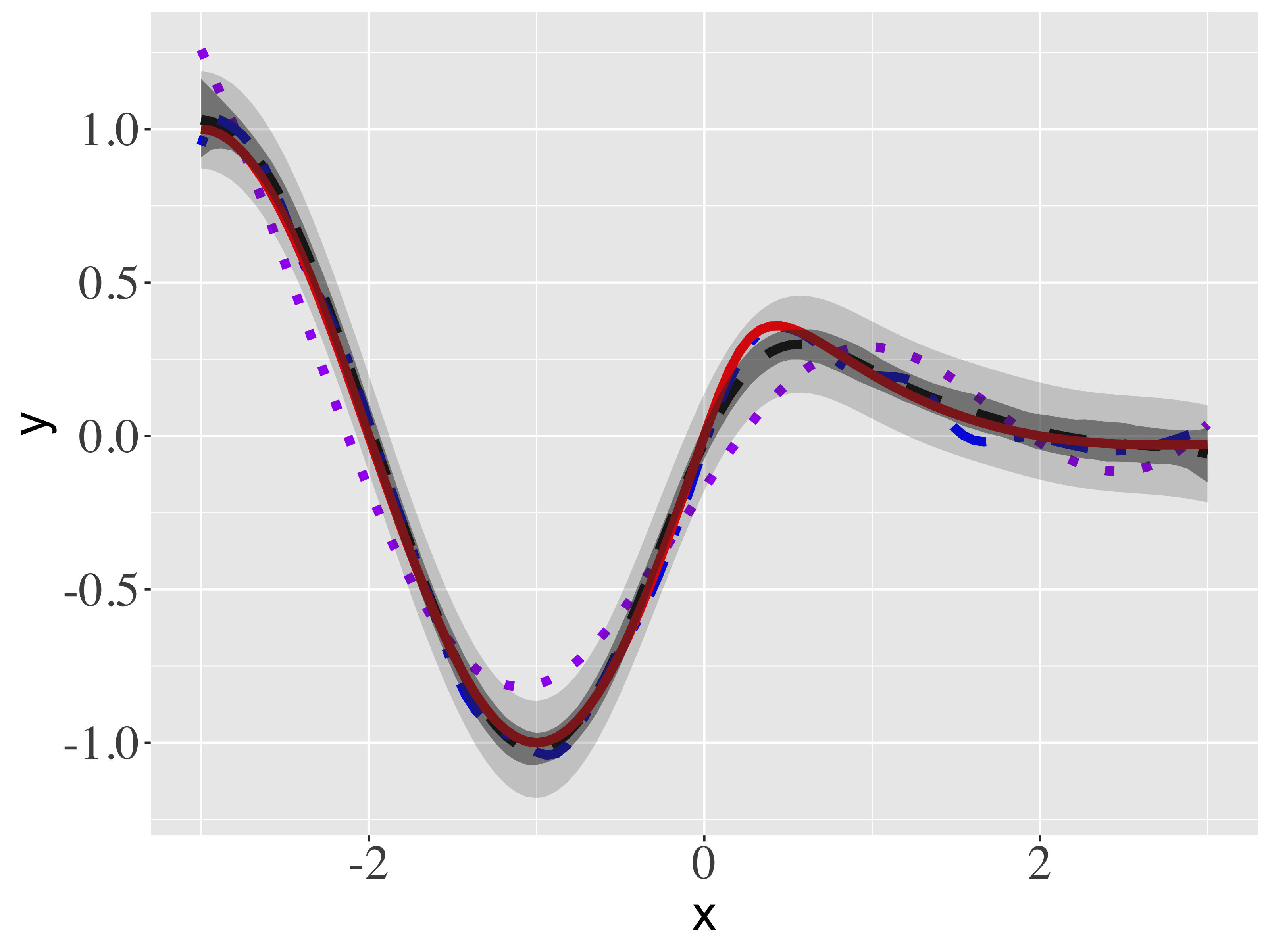}
    \includegraphics[scale = 0.045]{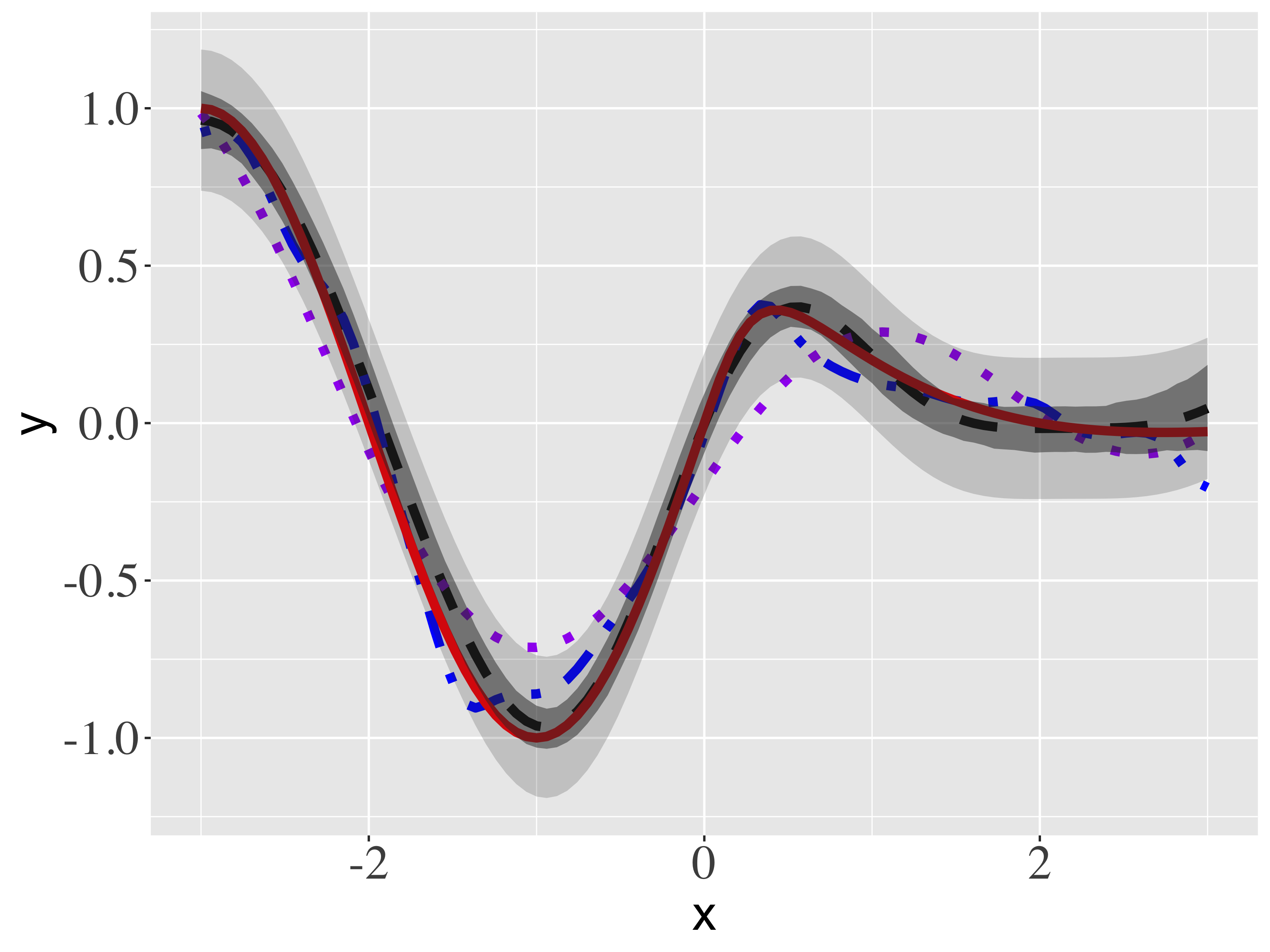}
    \includegraphics[scale = 0.045]{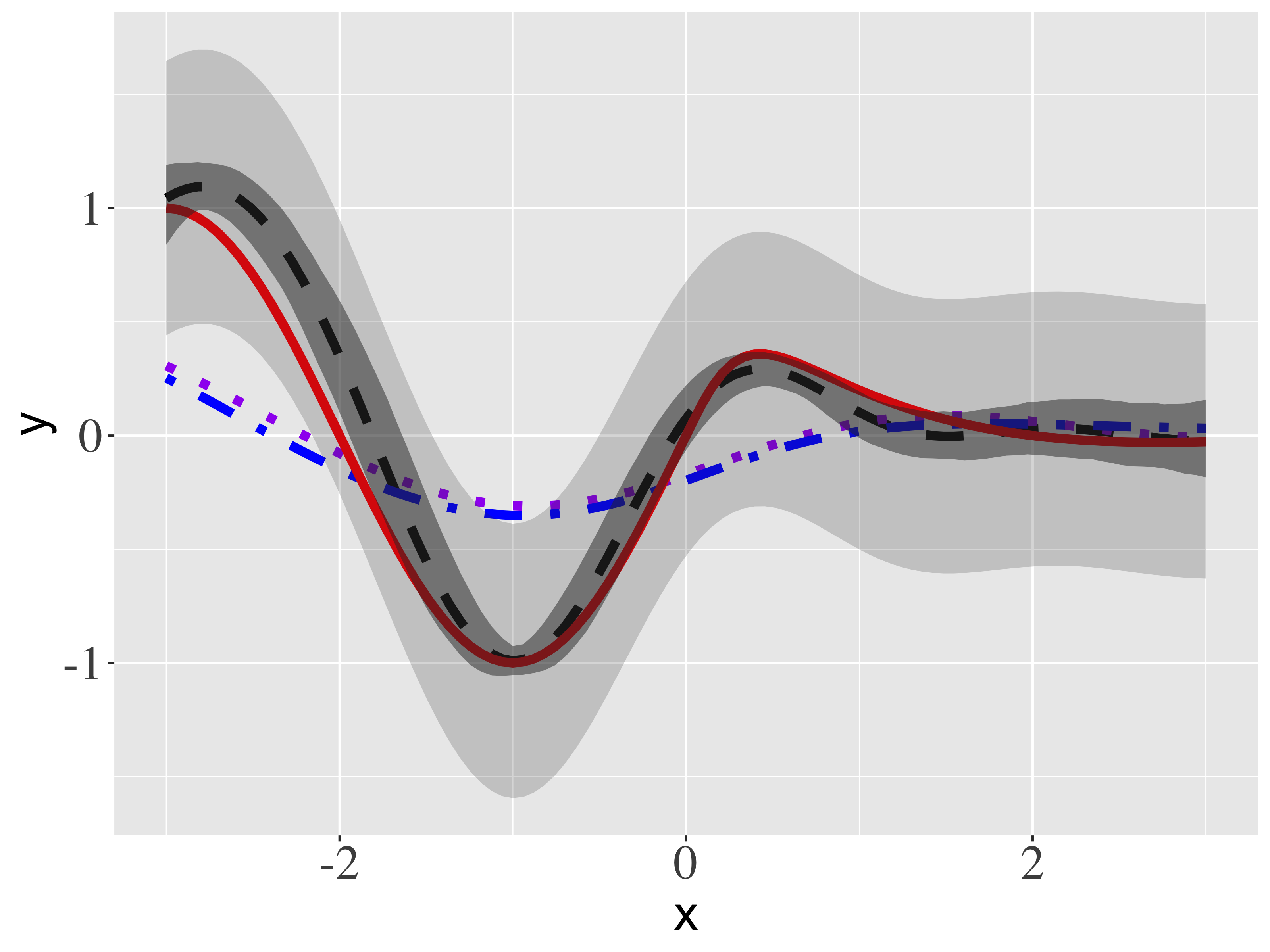}

\end{multicols}

\caption{Out-of-sample predictions of $f(x)$ with $n =100$ (first row) and $n=500$ (second row). For $n=100$, set $\delta^2 = 0.01$ (left panel), $\delta^2 = 0.6$ (middle panel) and  $\delta^2 = 1$ (right panel). For $n=500$, set $\delta^2 = 0.005$ (left panel), $\delta^2 = 0.1$ (middle panel) and  $\delta^2 = 1$ (right panel). In each plot, the red solid line is the true function, the black dashed line stands for the prediction based on \textsc{gpev}$_a$, the blue dot-dashed line based on \text{decon}, the purple dotted line based on \textsc{gp}. The darker and the lighter shades are the pointwise and simultaneous $95\%$ credible intervals obtained with \textsc{gpev}$_a$.} \small

\label{fig:fit_n500}
\end{figure}

Table \ref{tab:f1_mse} summarizes out-of-sample prediction results for all methods in terms of the averaged mean squared errors (\textsc{amse}) given different values of $\delta^2$. The results show \textsc{gpev}$_f$ performs the best among compared methods.  However, we observe that the performance of \textsc{gpev}$_a$ is very close to that of \textsc{gpev}$_f$ for all combinations of $n$ and $\delta^2$. This observation suggests that the approximation error of the proposed \textsc{gp} surrogate to the original \textsc{gp} is almost negligible in out-of-sample prediction despite the magnitude of the measurement error. We now investigate the performance of considered models in detail against the noise level of measurement errors.  When $\delta^2$ is small, all considered methods are robust to the measurement error except the \textsc{gp} model, implying that ignoring measurement errors could compromise the estimation significantly even though covariates are mildly contaminated.  As $\delta^2$ increases,  the \textsc{amse}s for \textsc{gpev}$_a$ and \textsc{gpev}$_f$ increase only by a relatively small margin, whereas other methods have suffered a drastic increase in \textsc{amse} values. For instance, the \textsc{amse} values obtained by \textsc{gp} and \text{decon} models are three times greater than those by \textsc{gpev}$_a$ and \textsc{gpev}$_f$ when $\delta^2\ge 0.6$. The robustness of \textsc{gpev}-based models to large measurement errors empirically justifies our theoretical claim that a \textsc{dpmm} prior is necessary for recovering the covariate density and thus allows the regression recovery to be robust to measurement errors. 

Similar results can be also observed from the boxplots of mean squared error (\textsc{mse}) values in Figure \ref{fig:boxplot_n500}. The increasing amount of \textsc{mse}s for all methods along with $\delta^2$ is due to that the true covariate density $p_0$ turns harder to recover when the errors in covariates become more disturbing. On the other hand, this implies that increasing the number of replicates can improve the performance significantly of the Bayesian estimator in practice. Beyond the investigation on \textsc{mse}s, the model fitting result in Figure \ref{fig:fit_n500} graphically displays the prediction performance of compared methods over various values of $\delta^2$. In particular, one can observe that when $\delta^2$ increases the performance of \text{decon} and \textsc{gp} deteriorates fast and both fail to recover the curvature of the true function. On the contrary, even when $\delta^2 =1$, the posterior mean of \textsc{gpev}$_a$ aligns with the true curve closely and its $95\%$ pointwise credible interval contains the whole true function. A wider credible interval for larger values of $\delta^2$ is expected due to an increasing amount of uncertainty in retrieving the covariate density. Overall, the \textsc{gpev}-based models stand out among other competitors in terms of prediction.

A careful inspection of our theory implies that placing a component-wise normal prior on the covariate results in a slower posterior contraction rate in recovering both the true covariate density and the true regression function. This is supported by the empirical observation that much larger \textsc{amse} values are obtained by \textsc{gpev}$_n$ when $\delta^2$ becomes large. Additional investigation regarding the recovery of covariate can be found in Figure \ref{fig:p0} of Appendix~\ref{app:addres}. By comparing the posterior density function of covariates based on \textsc{gpev}$_a$ and \textsc{gpev}$_n$, one can see that a component-wise normal prior is unable to identify the true covariate from the contaminated observations. 
{\tcr In Figure \ref{fig:p0} in Appendix \ref{app:addres}, we display a few examples of the posterior marginal density function obtained by \textsc{gpev}$_a$ and \textsc{gpev}$_n$, when $n=500$ and $\delta^2 =0.001, 0.1, 0.5$, respectively. Recall that the true covariate distribution is $\text{Unif}[-3,3]$ in our simulation setting. When $\delta^2=0.001$, both \textsc{gpev}$_a$ and \textsc{gpev}$_n$ recover the true underlying density quite well,  indicating that a DPMM prior on the marginal density has a similar performance with independent normal priors on the locations. When $\delta^2=0.1,0.5$, the performance of both methods deteriorate dramatically in estimating the marginal density function in Figure \ref{fig:p0}, which is expected since the best obtainable rate of convergence becomes slower with large $\delta^2$. However, as $\delta^2$ increases, one can still notice an improvement in estimating the covariate density using \textsc{gpev}$_a$. The posterior density function of the covariate obtained from \textsc{gpev}$_n$  is more erratic, suggesting that assigning independent normal priors on locations results in a poor recovery when the measurement errors are more significant. In addition, we compared these two methods in terms of the averaged mean squared error in recovering true locations in all cases of sample sizes in Table \ref{x_mse} in Appendix~\ref{app:addres}, which tells a similar story regarding the performance of \textsc{gpev}$_a$ and \textsc{gpev}$_n$. In each iteration of the Gibbs sampler, we update the covariate values and update the rest of parameters upon those, a better performance of recovering the true locations leads to a better result in updating other parameters, which partially explains the outperformance of \textsc{gpev}$_a$ in estimating the regression curve.   
} 

In addition to a comparable performance in prediction, \textsc{gpev}$_a$  is more computationally efficient than \textsc{gpev}$_f$. \textsc{gpev}$_a$ avoids repeated computation of the inverse of covariance matrix associated with a full \textsc{gp}, at a price of updating hyperparameters of a relatively moderate size (a fraction of sample size) related to Fourier basis functions. This is particularly beneficial for the errors-in-variables problem as covariates are treated as unknown parameters and both covariates and the covariance matrices are  updated in each iteration. Also, for \textsc{gp} models, the choice of covariance kernel and treatment to the associated hyperparameters play an important role in the mixing of the Markov chains \citep{murray2010slice}. To implement \textsc{gpev}$_a$, we  consider a squared exponential covariance kernel associated witha bandwidth parameter, denoted by $\lambda$, which is treated as an unknown parameter. The conjugate form of its spectral density induces a closed-form conditional of the bandwidth parameter  $\lambda$ based on the Fourier basis representation. Figure \ref{fig:mixing_lambda} in Appendix~\ref{app:addres} shows the trace plots of posterior samples of the bandwidth parameter $\lambda$, where one can observe that the mixing of the chain based on \textsc{gpev}$_a$ is much better than that based on \textsc{gpev}$_f$. We also remark that auto-correlation of the Markov chains obtained from \textsc{gpev}$_a$ and \textsc{gpev}$_f$ are similar, which is from the boxplots of the effective sample sizes (\textsc{ess}) of estimated function values based on \textsc{gpev}$_a$ and \textsc{gpev}$_f$ over training data points in Figure \ref{fig:ess-boxplot} of Appendix~\ref{app:addres}. To gauge the computational efficiency of \textsc{gpev}-based methods, we report that the computation time of \textsc{gpev}$_a$, \textsc{gpev}$_n$, \textsc{gpev}$_f$ for a single Markov chain iteration when $n = 500$ are $0.025,0.022,0.197$ second separately, on an 8-Core Intel Core i9 computer with 32 GB RAM. It is evident that implementing the proposed \textsc{gp} surrogate improves the computation speed substantially and the improvement  becomes more pronounced as the sample size increases. In conclusion, \textsc{gpev}$_a$ stands out as a more robust and computationally efficient method for tackling the errors-in-variables regression problem. 

\section{A Case Study}\label{sec:real}
We re-analyzed the real data set studied in \cite{berry2002bayesian} using the proposed \textsc{gpev} method. As described in \cite{berry2002bayesian}, the data set was collected from a randomized study where the actual content is not allowed to be disclosed. Basically, the data contains a treatment group and a control group. In each group we have the surrogate measurement $W$ evaluated at baseline, and the observed response $Y$ evaluated at the end of study. We know smaller values of $W$ and $Y$ indicate a worse case in the study. As discussed in \cite{berry2002bayesian}, the quantity of interest is the change from the baseline $\Delta (X) = f(X)-X$.  We assume a normal zero-mean measurement error with two choices of variance, 1) a fixed variance $\delta^2 = 0.35$, adopting the estimated value from the study; and 2) an unknown variance $\delta^2$ which will be treated as an unknown parameter in the model.  To implement the \textsc{gpev}$_a$ model, we choose $N=60$ based on the simulation results, and consider an exponential prior $\exp{(\lambda_0)}$ with hyperparameter $\lambda_0=1.5$ on the bandwidth parameter $\lambda$ associated with the squared exponential kernel. Besides, we treat the response error variance $\sigma^2$ as unknown and we consider an objective prior for $\sigma^2$, namely, $\Pi(\sigma^2) \propto 1/\sigma^2$, allowing the data to update the parameter. To update $\sigma^2$, we refer to Step 7 of the Gibbs sampler in Appendix~\ref{sec:gibbcom}. For both cases of $\delta^2$, we ran the Gibbs sampler with $1,500$ iterations with the first $1,000$ being discarded as a burn-in. We consider the posterior mean as our Bayesian estimator and report the $95\%$ pointwise credible interval. 

\begin{figure}[h!]
\begin{center}
    \includegraphics[scale = 0.1]{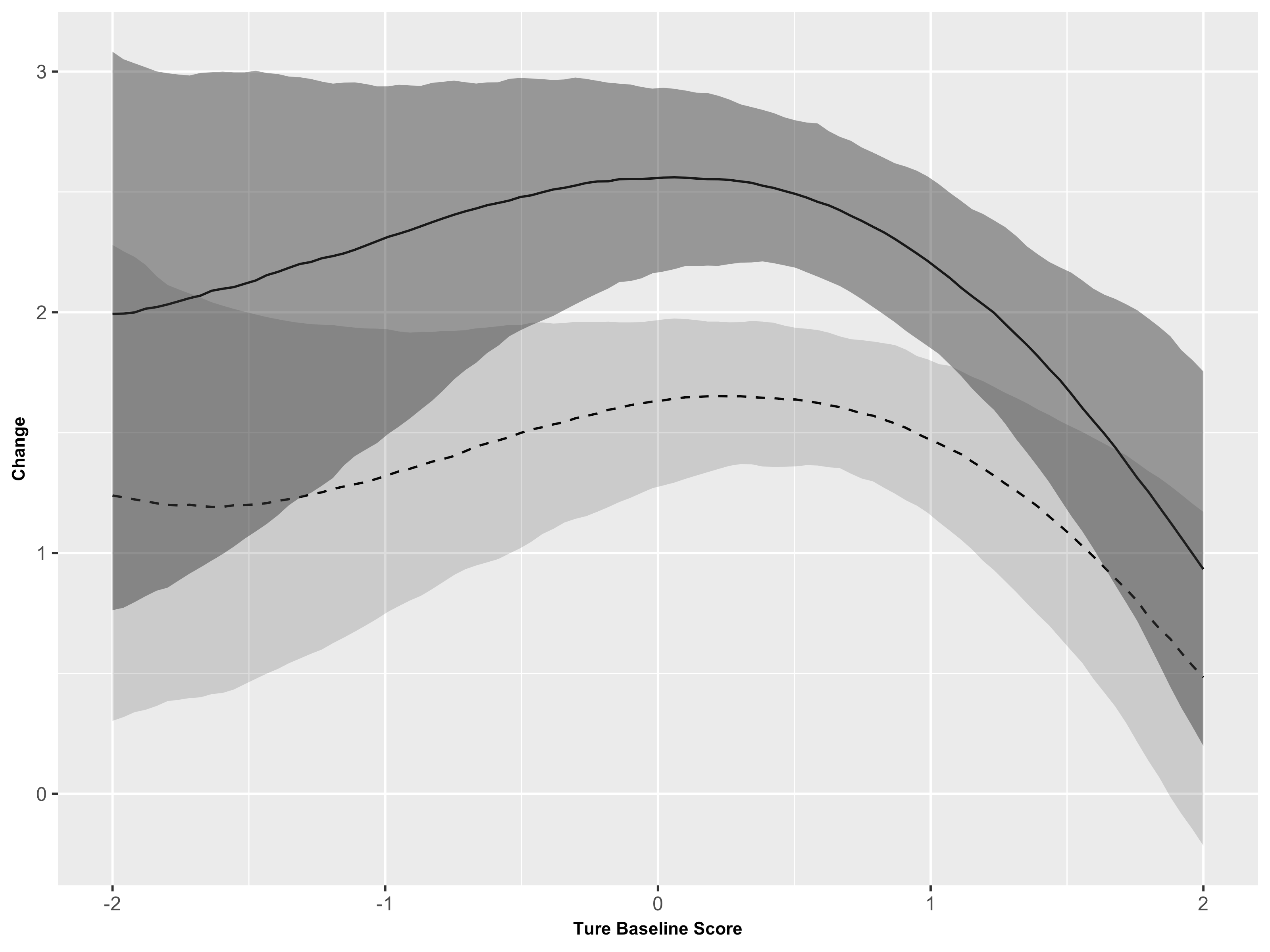}
 \end{center}
\caption{Estimate of $\Delta(X)$ on an evenly spaced test grid over $[-2,2]$ with $\delta^2 = 0.35$. The solid line indicates the treatment group with the darker shade as its $95\%$ pointwise credible interval and the dashed line indicates the control group with the lighter shade as its $95\%$ pointwise credible interval.}
\label{fig:plot}
\end{figure}

\begin{figure}[h!]
\begin{center}
    \includegraphics[scale = 0.068]{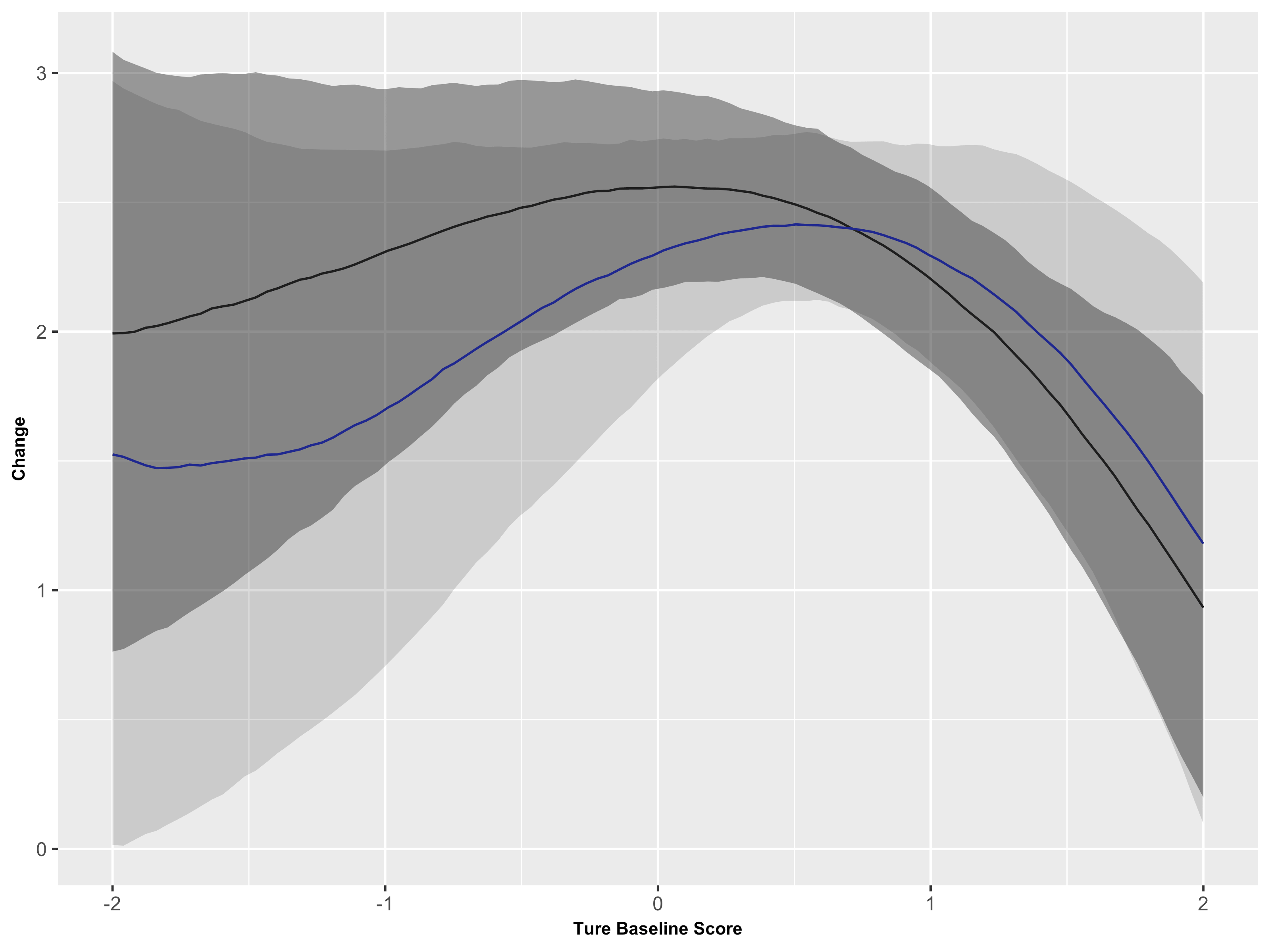}
    \includegraphics[scale = 0.068]{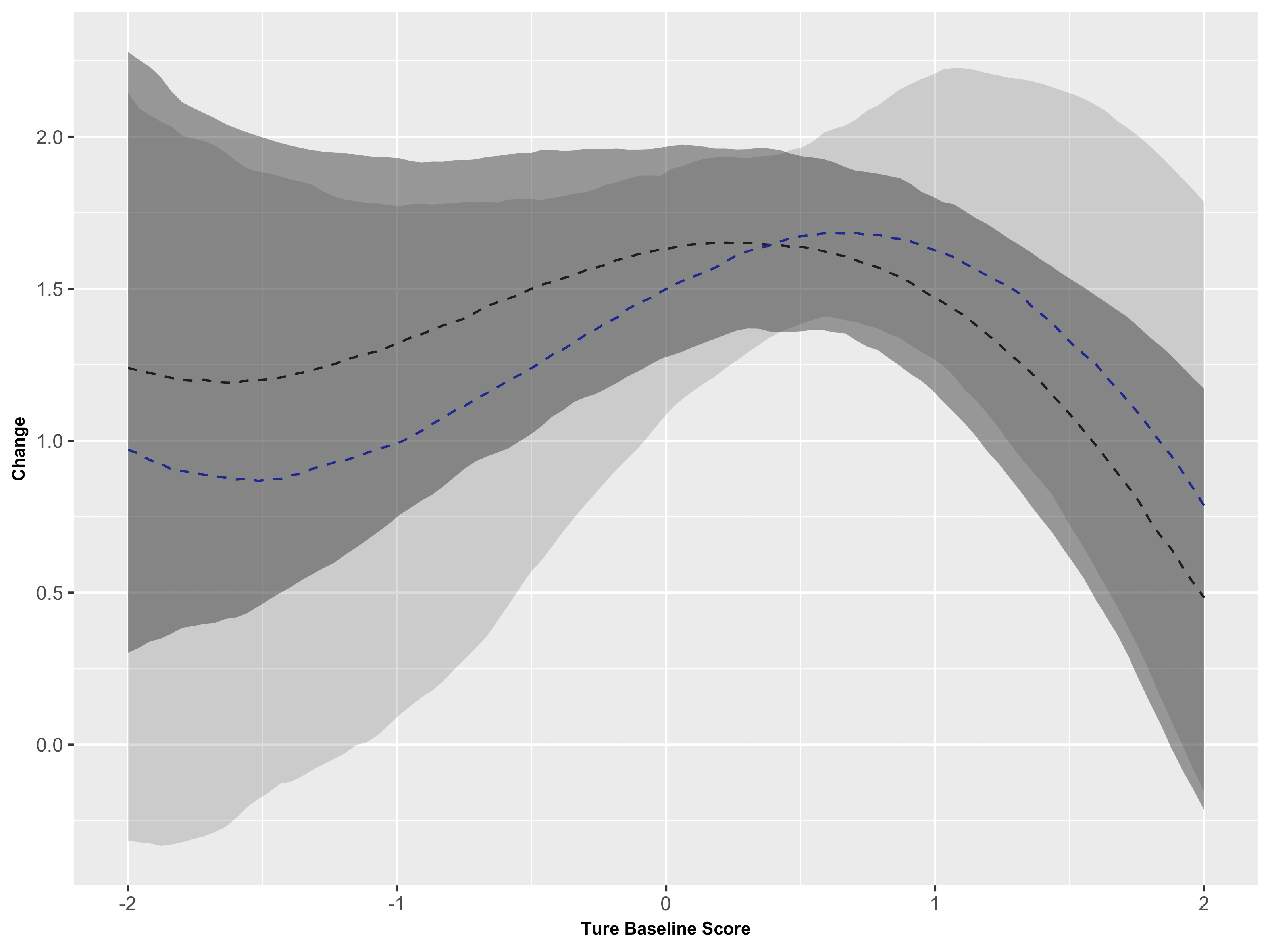}
    \end{center}
\caption{Estimate of $\Delta(X)$ on an evenly spaced test grid over $[-2,2]$ with different choices of $\delta^2$. The left panel indicates the treatment group and the right panel indicates the control group. In both panels, the black lines stand for the posterior mean with $\delta^2 = 0.35$ with the dark shades as $95\%$ pointwise credible intervals, the blue lines stand for the posterior mean with unknown $\delta^2$ with the light shades as $95\%$ pointwise credible intervals.}
\label{fig:error}
\end{figure}

Figure \ref{fig:plot} shows the prediction results of the changes by \textsc{gpev}$_a$ with $\delta^2 = 0.35$. We observe that for both the treatment and control groups, the change from the baseline increases first and then decreases as the true baseline score increases, which coincides with the  results presented in \cite{berry2002bayesian}. In Figure \ref{fig:error}, we compare the estimated changes by \textsc{gpev}$_a$ with fixed $\delta^2$ and unknown $\delta^2$ for both groups. We observe that for both treatment and control groups, an objective prior on $\delta^2$ results in a similar prediction of $\Delta(X)$ as that with fixed $\delta^2$. {\tcr Since the true regression and true covariates are unknown, we compute the mean squared error as $\textsc{mse} = \frac{1}{n}\sum_{i=1}^n\{(\hat{f}(\hat{x}_i) - y_i)^2\}$, where $\{y_i\}$ are observed responses and $\{\hat{x}_i\}$ denote the posterior mean of covariates obtained by \textsc{gpev}$_a$. The defined \textsc{mse} value accounts for randomness in the responses and errors in estimating the regression function and the covariates. Although this MSE value does not directly reflect the accuracy of predicting the true function, it provides some insights when comparing the performance of various methods. For the treatment group, the \textsc{mse} values for \textsc{gpev}$_a$ with $\delta^2=0.35$, \textsc{gpev}$_a$ with unknown $\delta^2$ and the decon method are $4.35$, $4.54$, and $8.23$ separately; and for the control group, the \textsc{mse} for the three competitors are $3.87$, $4.35$ and $5.18$, respectively. Theoretical results have shown that with relatively large $\delta^2$, all methods may obtain an extremely slow rate of convergence \citep{fan1991optimal,fan1993nonparametric}, which explains the large \textsc{mse} values.  However, MSE values for \textsc{gpev}$_a$ are smaller than those of decon for both control and treatment data, despite of knowledge of the measurement error variance, showing a superior performance to decon in the real data example.}  The diagnostic results show the mixing of Markov chains for $\{w_j, s_j, x_j\}$ are good in both scenarios for the \textsc{gpev}$_a$ model, refer to Figure \ref{fig:trace_treat} in Appendix \ref{app:addres} for more details on trace plots and density plots of posterior samples of selected subsets of $\{a_i,\omega_j,s_j\}$. 

\section{Discussion}
The article revisits error-in-variables regression problem from a Bayesian framework and addresses two fundamental challenges.  Theoretical guarantees on the convergence of the posterior are established for the first time in a Bayesian framework.  More specifically, optimal rates of posterior convergence are obtained simultaneously for the regression function as well as the covariate density. From a computational perspective, we provide a new Gaussian process approximation which facilitates posterior sampling and avoids costly matrix operations associated with a standard Gaussian process framework. 


In addition to showing weak convergence of the approximate Gaussian process to the original ones, we have also shown that when it is employed as a prior process, the resultant posterior maintains same contraction results as those of original GPs in recovering both the regression curve and covariate density function in EIV problem. As the procedure can be easily generalized to other nonparametric setting, this result implies some statistical guarantee of the predictive performance of projection technique with the random Fourier features under a Bayesian framework, which is a new addition to the theoretical investigations of the random Fourier features.



\acks{Pati's research was partially supported by NSF DMS (1613156, 1854731, 1916371) and Yang's research was partially supported by NSF DMS 1810831. The research of Wang and Carroll was supported by a grant from the National Cancer Institute (U01-CA057030).}


\newpage

\appendix
\section{Technical Results}\label{app}

Section~\ref{app:notns} introduces notations used throughout the rest of the paper and some background knowledge on the Gaussian process prior and its associated reproducing kernel Hilbert space. Section~\ref{sec:aux} collects all auxiliary results used to prove Theorem \ref{thm:mthm}. 

\subsection{Notations and Backgrounds} \label{app:notns}

We first introduce some notations used in the proofs. Denote by $E_{X}$ the marginal expectation with respect to random variable $X$; denote $\mbox{P}^{f,p}_{X,Y}$ as the probability measure of random pair $(X,Y)$ which has a joint density denoted by $(f,p)$. Let $*$ denote the convolution, say, for two functions $f$ and $g$ we define $f*g (\cdot) = \int f(\cdot-t)g(t)\,dt$. Denote the Kullback--Leibler divergence between functions $f$ and $g$ with respect to the Lebesgue measure $\mu$ by $KL(f,g) = \int f\log(f/g)\,d\mu$ and denote the second moment of the Kullback--Leibler divergence by $V(f,g) =\int f\{\log(f/g)\}^2\,d\mu$. Define the $\epsilon$-Kullback--Leibler neighborhood around $f_0$ as $B_{f_0}(\epsilon)=\{f: KL(f_0,f) \le \epsilon^2,V(f_0,f)  \le \epsilon^2\}$. We also define the Hellinger distance between two densities $f$ and $g$ as $H(f,g) = \{\int (\sqrt{f}-\sqrt{g})^2d\mu\}^{1/2}$. And define the $L_2$-norm as $\|f\|_2 = \sqrt{\int f(x)^2dx}$. Let $\mathds{1}_C(\cdot)$ denote the indicator function on any set $C \subset \mb R$. For two sets $A,B$, we denote their Cartesian product by $A\otimes B$, the set contains all pairs $(x,y)$, where $x\in A$ and $y \in B$. For two positive sequences $a_n, b_n$, we write $a_n = O(b_n)$ if $\lim_{n\to \infty}(a_n/b_n)=c$ for some constant $c>0$, and $a_n=o(b_n)$ if $\lim_{n\to \infty}(a_n/b_n)=0$. At last, we define a $k$th order kernel function $K(\cdot)$ that satisfies,
\begin{eqnarray}\label{eq:kernel1}
\int K(u) \,du =1, \quad \int K^2(u)\,du < \infty, \quad \int u^{\floorbeta} K(u)\, du \neq 0,\nonumber\\
\int u^{i-1} K(u)\, du = 0, \quad \text{for} \ i=1, \dots, \floorbeta-1, \, \beta \ge 2.
\end{eqnarray}

Now we briefly recall the definition of the reproducing kernel
Hilbert space of a Gaussian process prior; a detailed review can be found in \cite{van2008reproducing}. A Borel measurable random element $W$ with values in a separable Banach space denoted by $(\mathbb{B}, \norm{\cdot})$, for instance, the space of continuous functions $C[0,1]$, is called Gaussian if the random variable $b^*W$ is normally distributed for any element $b^* \in \mathbb{B}^*$, the dual space of $\mathbb{B}$. The reproducing kernel Hilbert space $\mathbb{H}$ attached to a zero-mean Gaussian process $W$ is
defined as the completion of the linear space of functions $t \mapsto EW(t)H$
relative to the inner product
\begin{eqnarray*}
\langle \mbox{E} (W(\cdot)H_1); \mbox{E} (W(\cdot )H_2)\rangle_{\mathbb{H}} = \mbox{E} (H_1H_2),
\end{eqnarray*}
where $H, H_1$ and $H_2$ are finite linear combinations of the form $\sum_{i}a_{i}W(s_i)$
with $a_i \in \mathbb{R}$ and $s_i$ in the index set of $W$.

Let $W = (W_t: t \in \mathbb{R})$ be a Gaussian process associated with a squared exponential covariance kernel, which is 
\begin{eqnarray*}
C(t,t') = e^{-(t-t')^2}.
\end{eqnarray*}
 The spectral measure $m_{w}$ of $W$ is absolutely continuous with respect to the Lebesgue measure $\lambda$ on $\mathbb{R}$ with the Radon-Nikodym derivative given by
\begin{eqnarray*}
\frac{dm_{w}}{d\lambda}(x) = \frac{1}{(2\pi)^{1/2}}e^{-x^2/4}.
\end{eqnarray*}

Define a scaled Gaussian process $W^a=(W_{at}: t \in [0,1])$, viewed as a map in $C[0,1]$. Let $\mathbb{H}^a$ denote the reproducing kernel Hilbert space of $W^a$, with the corresponding norm $\norm{\cdot}_{\mathbb{H}^a}$. The unit balls in reproducing kernel Hilbert space and in the Banach space are denoted by $\mathbb{H}^a_1$ and $\mathbb{B}_1$, respectively.

Next we describe the construction of the sieve $\mathcal{P}_n$ on the parameter space of $(f,p)$, the parameter space of $p$. For fixed constants $m, \underline{\sigma}, \bar{\sigma} >0$ and integer $H\ge 1$. Let
\begin{align*}
\mathcal{F}=\bigg\{ p_{F,\widetilde{\sigma}} = \phi_{\widetilde{\sigma}} * F: F=\sum_{h=1}^{\infty} \pi_h \delta_{z_h}, z_h \in [-m,m], h \le H, \sum_{h>H} \pi_h < \epsilon_n, \underline{\sigma} \leq \widetilde{\sigma} < \bar{\sigma} \bigg\}.
\end{align*}
Set  $\mathcal{P}_n = \widetilde{B}_n \otimes \mathcal{F}$, where $\widetilde{B}_n=B_n \cap \m{A}$ with $B_n= M_n\mathbb{H}_1^{a_n} + \epsilon_n\mathbb{B}_1$ and $\m{A}$ as in Assumption \ref{ass:A}.


\subsection{Auxiliary Results} \label{sec:aux} 
In this section, we collect auxiliary results that are needed for the proofs of main theorems. The proofs of Lemmata \ref{lem:sieve}-\ref{lem:gkl} are deferred to Appendix \ref{sec:auxiliary}. 

\begin{lem}\label{lem:sieve}
Suppose Assumptions \ref{ass:f0}, \ref{ass:p0}, \ref{ass:A} and \ref{ass:sig} hold, by taking $M_n=a_n=\epsilon_n^{-1/\beta}$, $H\lesssim n\epsilon^2_n, m^{\tau_1}\lesssim n, \underline{\sigma}\lesssim n^{-1/2\tau_2}$ and $\bar{\sigma}^{2\tau_3}\lesssim e^n$, we have $\Pi({\mathcal{P}_n^c}) \le e^{-n\epsilon_n^2}$ with $\epsilon_n = n^{-\beta/(2\beta+1)}(\log n)^t$, where $t=\max{\{(2\vee q)\beta/(2\beta+1),1\}}$.
\end{lem}


 \begin{lem}\label{lem:test}
 For model \eqref{eq:model}, $\widehat{p}_n$ and $\widehat{f}_n$ defined in Equations \eqref{eq: pn} and \eqref{eq: fn} and $f, p \in \mathcal{P}_n$ for any small constant $\epsilon_0 >0$,
 \begin{eqnarray}
   P_{W,X}^p(|| \widehat{p}_n - p ||_{\infty} > \epsilon_0) \le e^{-C_1n\epsilon_0 h_n^2},\label{test1}  \\
 P_{W,X}^{p}(||\widehat{p}_n-p||_1> \epsilon_n) \le e^{-n\epsilon_n^2}, \label{test2}\\
  P_{Y,W,X}^{f,p}(||\widehat{f}_n\, \widehat{p}_n-f\, p||_1> \epsilon_n) \le e^{-n\epsilon_n^2},\label{test3}
\end{eqnarray}
where $h_n \asymp \epsilon_n^{1/\beta}$, $\epsilon_n = n^{-\beta/(2\beta+1)}(\log n)^t$ with $t=\max{\{(2\vee q)\beta/(2\beta+1),1\}}$ and some constant $C_1 >0$.
\end{lem}

\begin{lem}\label{lem:priorthick}
Suppose Assumptions \ref{ass:p0}, \ref{ass:A} and \ref{ass:sig} hold, then $\Pi\{KL(p_0,\epsilon_n)\} \ge e^{-n\epsilon_n^2}$, where $\epsilon_n = n^{-\beta/(2\beta+1)}(\log n)^t$ with $t=\max{\{(2\vee q)\beta/(2\beta+1),1\}}$.
\end{lem}

\begin{lem}\label{lem:gkl}
Under the conditions in Theorem \ref{thm:mthm} and suppose Lemma \ref{lem:priorthick} hold, for sufficiently large $n$,
\begin{align*}
\Pi \big\{B_{(f_0, p_0)}(\epsilon_n) \big\} \ge e^{-n\epsilon_n^2},
\end{align*}
where $B_{(f_0, p_0)}(\epsilon_n)$ is defined in \eqref{KL_fp} and 
$\epsilon_n = n^{-\beta/(2\beta+1)}(\log n)^t$ with $t=\max{\{(2\vee q)\beta/(2\beta+1),1\}}$.
\end{lem}

\begin{lem} \label{lem:Talagrand}
(Theorem 7.3 in \citealt{bousquetconcentration}) Suppose $\m{G}$ is a countable set of functions $g: \m{X} \to \mathbb{R}$ and assume all functions $g\in \m{G}$ are measurable, squared-integrable and satisfy $E\{g(X_k)\}=0$. Assume $\sup_{g\in\m{G}} \operatorname*{ess\,sup} g$ is bounded and define $Z=\sup_{g\in\m{G}}\sum_{k=1}^n g(X_k)$. Let $\sigma_{\m{G}}$ be a positive real number such that $n\sigma^2_{\m{G}} \ge \sum_{k=1}^n \sup_{g\in\m{G}} E \{g^2(X_k)\}$, then for all $t >0$ with $\nu=n\sigma^2_{\m{G}}+2 E(Z)$, we have
\begin{align*}
P\bigg\{Z \ge \mbox{E}(Z)+(2t\nu)^{1/2}+\frac{t}{3}\bigg\} \le e^{-t}.
\end{align*}
\end{lem}

\begin{lem}\label{lem:Borell}
(Borell's inequality in \citealt{adler1990}) Let $\{f(x):x\in [0,1]\} $ be a centered Gaussian process and denote $\|f\|_{\infty}=\sup_{x\in[0,1]}f(x)$ and $\sigma_{f}^2=\sup_{x\in[0,1]}E\{f^2(x)\}$. Then $E(\|f\|_{\infty}) < \infty$ and for any $t>0$,
\begin{eqnarray*}
P(| \|f\|_{\infty}- E\|f\|_{\infty} | >t) \le 2e^{-\frac{1}{2}t^2/\sigma^2_{f}}.
\end{eqnarray*}
\end{lem}

\section{Proof of Theorem \ref{thm:mthm}} \label{app:mthm}
In this section, we provide the proof of Theorem \ref{thm:mthm}. Given $\epsilon_n,\epsilon_n'$ in Theorem \ref{thm:mthm}, define $U_n=\{f,p: ||f-f_0||_1 < M\epsilon_n, ||p-p_0||_1 < M\epsilon'_n \}$, our goal is to show $ \Pi_n (U_n^c \mid Y_{1:n}, W_{1:n}) \to 0$ almost surely in $G_{f_0,p_0}$ as $n\to \infty$. To that end, note that 
\begin{align}\label{post_prob}
 &\Pi_n (U_n^c \mid Y_{1:n}, W_{1:n})\nonumber\\
 &\le  \Pi_n (f,p: ||f-f_0||_1 > M\epsilon_n \mid Y_{1:n}, W_{1:n}) +\Pi_n (p:  ||p-p_0||_1 > M\epsilon'_n \mid W_{1:n}) := S_1+S_2.
 \end{align} 
It suffices to estimate $S_1$ and $S_2$ in the preceding separately. We shall analyze term $S_1$ in detail and only provide a brief discussion about bounding term $S_2$ as it can be considered as an immediate application of existing results. \\ 
 
\noindent{\em Bounding term $S_1$ in Equation \eqref{post_prob}.} Define the $\epsilon_n$-Kullback--Leibler neighborhood around $(f_0, p_0)$ as
\begin{eqnarray}\label{KL_fp}
B_{f_0,p_0} (\epsilon_n) = \bigg\{  \int g_{f_0,p_0} \log\frac{g_{f_0,p_0}}{g_{f,p}}  \leq \epsilon_n^2, \quad  \int g_{f_0,p_0}\bigg(\log\frac{g_{f_0,p_0}}{g_{f,p}}\bigg)^2 \leq \epsilon_n^2 \bigg\}.
\end{eqnarray}
The following Theorem provides sufficient conditions showing $S_1\to 0$ almost surely as $n\to \infty$.  A sketch of the proof is provided in the following.
\begin{thm} {(Contraction Theorem)} \label{thm:thm2}
Consider model \eqref{eq:model} and under the conditions in Theorem \ref{thm:mthm}, let $\mathcal{U}_n = \{||f-f_0||_1 > M\epsilon_n\}$. If there exist a sequence of $\epsilon_n \to 0$ and $n\epsilon_n^2\to \infty$ and a sequence of sieve $\mathcal{P}_n \subset \mathcal{P}$, and a sequence of test functions $\phi_n = \mathds{1}_{\{||\widehat{f}_n -f_0||_1 >(M-1)\epsilon_n\}}$ satisfying the following conditions,
\begin{align}
G_{f_0,p_0}\, \phi_n \le e^{-(C+4)n\epsilon_n^2}, \quad \sup_{(f,p) \in \mathcal{P}_n \cap\ \mathcal{U}_n} G_{f, p} \,(1- \phi_n) \le e^{-(C+4)n\epsilon_n^2},\label{eq:test}\\
\Pi \big \{B_{f_0,p_0} (\epsilon_n) \big\}  \ge e^{-n\epsilon_n^2}, \label{eq:priorconc}\\
\Pi({\mathcal{P}_n^c}) \le e^{-(C+4)n\epsilon_n^2} \label{eq:sieve},
\end{align}
 for some constant $C>0$, then $\Pi_n( \mathcal{U}^c_n \mid Y_{1:n}, W_{1:n}) \to 0$ almost surely in $G_{f_0,p_0}$, for
the constant $M$ same as in Theorem \ref{thm:mthm}.
\end{thm}
\begin{proof} {(Sketch)}
Define the set
\begin{eqnarray*}
C_n=\bigg\{\displaystyle \int {\frac{\Pi_{j=1}^n g_{f,p}\,(Y_j, W_j)}{\Pi_{j=1}^n g_{f_0,p_0}\,(Y_j, W_j)} }d\Pi(f)d\Pi(p) \ge e^{-(C+3)n\epsilon^2_n} \,\Pi\big\{B_{f_0,p_0} (\epsilon_n)\big\}\bigg\}.
\end{eqnarray*}
Under the conditions in Theorem \ref{thm:mthm}, from Lemma 8.1 in  \cite{ghosal2000}, it follows  $G_{f_0,p_0}(C_n) \ge 1-1/(C'n\epsilon_n^2)$, for some constant $C'>0$.    Hence for any sequence of test functions $\{\phi_n\}$,
\begin{align*}
\Pi_n( \mathcal{U}_n^c \mid Y_{1:n}, W_{1:n}) &\le G_{f_0,p_0}\phi_n+G_{f_0,p_0}\{ (C_n^c)
+G_{f_0,p_0}\Pi(\mathcal{P}_n^c \mid Y_{1:n}, W_{1:n})\mathds{1}_{C_n}\}\\
&+G_{f_0,p_0}\{\Pi(\mathcal{U}_n\cap\mathcal{P}_n\mid Y_{1:n}, W_{1:n})\,(1-\phi_n)\mathds{1}_{C_n}\}.
\end{align*}
We suppress the term ``almost surely" in the following argument. According to Conditions \eqref{eq:priorconc} and \eqref{eq:sieve}, the third term in the above display goes to 0.  From Conditions \eqref{eq:test} and \eqref{eq:priorconc}, the first and the fourth terms in the preceding go to 0. Then we have shown $S_1\to 0$ as $n\to \infty$. 
\end{proof}

We now verify three conditions in Theorem \ref{thm:thm2} under the conditions in Theorem \ref{thm:mthm}, based on the auxiliary results summarized in Appendix~\ref{sec:aux}. The main steps are
\begin{itemize}
\item Condition \eqref{eq:priorconc} of Theorem \ref{thm:thm2}:  Follows from Lemma \ref{lem:gkl} under the conditions of Theorem \ref{thm:mthm}.
\item Condition \eqref{eq:sieve} of  Theorem \ref{thm:thm2}:  Follows from Lemma \ref{lem:sieve} under the conditions of Theorem \ref{thm:mthm}.
\item Condition \eqref{eq:test} of Theorem \ref{thm:thm2}:  For model \eqref{eq:model}, recall $\widehat{p}_n$ and $\widehat{f}_n$ defined in Equations \eqref{eq: pn} and \eqref{eq: fn}, and $f, p \in \mathcal{P}_n$,
it suffices to estimate $P_{Y,W,X}^{f_0,p_0} (||\widehat{f}_n - f_0||_1 > \epsilon_n)$ and $P_{Y,W,X}^{f,p} (||\widehat{f}_n - f||_1 > \epsilon_n)$. Following a similar line of argument in \cite{Meister2009}, for any marginal density $p_0$ satisfying Assumption \ref{ass:p0}, for any $p\in \m{P}_n\cup p_0$ define $\Delta p=(\widehat{p}_n-p)/p$, then for any $f \in \m{P}_n\cup f_0$ we have
\begin{align*}
|\widehat{f}_n - f| \le \frac{|\widehat{f}_n\widehat{p}_n - fp|}{|p|}\bigg(\frac{|\Delta p|}{|\Delta p+1|}+1\bigg) + |f| \frac{|\Delta p|}{|\Delta p+1|}.
\end{align*}
By Assumption \ref{ass:p0}, $p_0$ is lower-bounded by some constant $B^{-1}>0$. Then applying the Equantion \eqref{test1} in Lemma \ref{lem:test}, for any constant $\epsilon_0>0$ we have $||\widehat{p}_n -p||_{\infty} < \epsilon_0$ with probability at least $1-e^{-n\epsilon_0h_n^2}$. Define the set $\m A_\epsilon = \{p: ||\widehat{p}_n -p||_{\infty} < \epsilon_0\}$. Thus for $p\in \m{P}_n\cap \m A_\epsilon$, $||p-p_0||_{\infty} \le ||\widehat{p}_n-p_0||_{\infty} + ||\widehat{p}_n-p||_{\infty} \le 2\epsilon_0$. Then $||p||_{\infty} \ge ||p_0||_{\infty} -||p-p_0||_{\infty}\ge B_1$, for some constant $B_1>0$ by choosing $\epsilon_0<B^{-1}/2$. Thus for $f\in\m{P}_n\cup f_0$ and $p\in \m{P}_n\cap \m A_\epsilon$, we have
\begin{align}\label{eq:tbnd}
||\widehat{f}_n - f||_1 \le \frac{1}{B_1} ||\widehat{f}_n\widehat{p}_n - fp ||_1 \bigg( \bigg| \bigg|\frac{\Delta p}{\Delta p+1} \bigg| \bigg|_{\infty}+1 \bigg) + ||f||_{\infty}  \bigg| \bigg|\frac{\Delta p}{\Delta p+1} \bigg| \bigg|_1.
\end{align}
Since $\||\Delta p||_{\infty}\le\epsilon_0/B_1$, choosing $\epsilon_0$ such that $\epsilon_0/B_1\le 1/2$, then we have $||\Delta p/(\Delta p +1)||_{\infty} \le 1$ and $1/2 \le ||\Delta p+1||_{\infty}\le 3/2$ and therefore $1/||\Delta p+1||_{\infty} \le2$. Thus we have,
\begin{align*}
\bigg|\bigg|\frac{\Delta p}{\Delta p+1}\bigg|\bigg|_1 \le \frac{1}{||\Delta p+1||_{\infty} ||p||_{\infty}} \int_0^1 |\widehat{p}_n(x)-p(x)|dx \le \frac{2}{B_1}||\widehat{p}_n -p||_1.
\end{align*}
Similarly for $p=p_0 \in \m A_\epsilon$, we bound $||\Delta p/(\Delta p+1)||_1\le 2B ||\widehat{p}_n-p_0||_1$.
Combining the above results and the result in Equation (\ref{eq:tbnd}), we obtain,
\begin{align*}
P(||\widehat{f}_n - f||_1 > \epsilon_n) &\le P(||\widehat{f}_n\cdot \widehat{p}_n-f\cdot p||_1>  B_1\epsilon_n/4)\\
&~~ + P\{||\widehat{p}_n-p||_1>B_1\epsilon_n/(4||f||_{\infty})\}\\
&~~+ P(||\widehat{p}_n-p||_{\infty}>\epsilon_0).
\end{align*}
Since we assume $f_0$ and $f\in \m{P}_n$ are bounded, applying Lemma \ref{lem:test} verifies Condition \eqref{eq:test}.
\end{itemize}

 \noindent{\em Bounding term $S_2$ in Equation \eqref{post_prob}.} To estimate $S_2$, we apply an inversion inequality built upon a special kernel function (the $sinc$ kernel) considered in \cite{donnet2018posterior}, then apply the existing posterior contraction result for the direct density problem. 
Recall the Fourier transform of the error density $\hat{\phi}_{\delta_n}(t) \asymp \delta_n e^{-\pi^2\delta_n^2t^2}$. Then with a careful inspection of the proof of Proposition 1 in \cite{donnet2018posterior}, one can obtain the inversion inequality 
\begin{align}\label{inv_inq}
\|p-p_0\|^2_2 &\lesssim \delta^{2\beta'} + \|\phi_{\delta_n} * p - \phi_{\delta_n} * p_0\|_1^2\times \int_{|t|\le 1/\delta}|\hat{\phi}_{\delta_n}|^{-2}dt \nonumber\\
&\lesssim \delta_n^{2\beta'} +\|\phi_{\delta_n} * p - \phi_{\delta_n} * p_0\|_1^2,
\end{align}
where $\beta'$ denotes the regularity level of $p_0$. The last inequality in Equation \eqref{inv_inq} holds by choosing $\delta \asymp \delta_n$ and the fact that $\int_{|t|\le 1/\delta}|\hat{\phi}_{\delta_n}|^{-2}dt \lesssim (\delta_n/\delta) e^{2\pi^2(\delta_n/\delta)^2} =O(1)$. Denote the observed density and the true density of $W$ by $f_W=\phi_{\delta_n} * p$ and $f_{0W}=\phi_{\delta_n} * p_0$ separately. By Cauchy-Schwarz inequality, Equation \eqref{inv_inq} implies $\|p-p_0\|_1 \lesssim  \max\{\delta_n^{\beta'}, \|f_W-f_{0W}\|_1\}$. Then under the Assumptions \ref{ass:p0} and \ref{ass:sig}, for $\delta_n,\epsilon'_n$ defined in Theorem \ref{thm:mthm}, one can easily show $S_2=o(1)$, by applying posterior contraction results for direct density estimation problem from the seminal work \citep{ghosal2001entropies,shen2013adaptive}, which leads to the error rate $\epsilon'_n = n^{-\beta'/(2\beta'+1)}(\log n)^{t'}$ with $t'>(\gamma+1/\beta')/(2+1/\beta')$ for some $\gamma>2$ under Assumption \ref{ass:p0}. 

Combining above results for terms $S_1$ and $S_2$ completes the proof of Theorem \ref{thm:mthm}, and letting $\epsilon_n = n^{-\beta/(2\beta+1)}(\log n)^t$ and $\epsilon'_n = n^{-\beta'/(2\beta'+1)}(\log n)^{t}$ with $t=\max{\{(2\vee q)\beta/(2\beta+1),t'\}}$ yields the desired rates in Theorem \ref{thm:mthm}. 



\section{\bf{Proof of Theorem \ref{thm:appgp}}}  \label{sec:appgp}
In this section, we provide a proof of Theorem \ref{thm:appgp}. In Part I, we shall first show the weak convergence of $\widetilde{f}_N$ to the original Gaussian process; in Part II, we derive expressions of  expectation and covariance of $\widetilde{f}_N$. 

\noindent{\em Part I.} We now show $\widetilde{f}_N$ weakly converges to the Gaussian process $f$. Based on Theorem 1.5.7 in \cite{van1996weak}, it suffices to show the marginal weak convergence and asymptotical tightness of $\widetilde{f}_N$. 

First, we show the marginal weak convergence. For any finite sequence $ \{x_1,\dots,x_k\}$ in $[0,1]$ of size $k$ where $k$ is arbitrary positive integer, applying multivariate central limit theorem with the expectation and covariance of $\wt{f}_N$ derived in Part II, one can easily show that as $N \to \infty$,
$$
\{\widetilde{f}_N(x_1),\dots, \widetilde{f}_N(x_k)\} \to \m \mbox{N}(0, c_{k,k}),
$$
in distribution, where $c_{k,k} = (c_{ij})$ is a $k\times k$ covariance matrix with the $(i,j)$th element $c_{ij} = c(x_i,x_j)$. 

Next, we show the asymptotic tightness of $\widetilde{f}_N$. By definition, it suffices to verify the following three conditions. First, it is straightforward that $[0,1]$ is totally bounded. Second, for any fixed $x_0\in[0,1]$, we shall show the tightness of $\widetilde{f}_N(x_0)$. It is equivalent to show, by definition, for any $\epsilon>0$, there exists a compact set $K$ such that,
\begin{eqnarray}\label{eq:tight}
P\{\widetilde{f} (x_0) \in K\} > 1-\epsilon. 
\end{eqnarray}
For any $x_0\in[0,1]$, we bound $\widetilde{f}_N(x_0)$ from above as 
 $$| \widetilde{f}_N(x_0) | \le (2/N)^{1/2}\,\sum_{i=1}^N |a_j|,$$
where $a_j\overset{i.i.d.}\sim \mbox{N}(0,1),\, j=1,\ldots,N$. It is well-known that $|a_j|$ is a sub-gaussian random variable for $j=1,\ldots,N$.  For any $t>0$, we have 
 \begin{align*}
P\{|\widetilde{f}_N(x_0)| \ge t \} \le P\bigg\{\,(2/N)^{1/2}\, \sum_{i=1}^N |a_j|\ge t \bigg\}  \le 2 \exp(-ct^2),
 \end{align*} 
for some constant $c>0$. 
For any $\epsilon >0$, we choose $t = \{2\log(1/\epsilon)\}^{1/2}$ and the set $K = \{ |\widetilde{f}(x_0)| \le t \}$, then Equation \eqref{eq:tight} holds, thus we show the tightness of $\widetilde{f}_N(x_0)$ for any $x_0\in [0,1]$. 

Third, we shall show $\widetilde{f}_N$ is asymptotically uniformly eqicontinuous with respect to the Euclidean norm, which is defined as $d(x,y) = |x-y|$, for $x,y \in \mathbb{R}$. It suffices to show that for any $\epsilon, \eta >0$, there exists some $\delta >0$ such that
\begin{eqnarray} \label{eq:eqi_cont}
\limsup_{N\to\infty}P\bigg\{\sup_{d(x,y) < \delta} |\widetilde{f}_N(x) -\widetilde{f}_N(y)|> \epsilon\bigg\} < \eta.
\end{eqnarray} 
 Without loss of the generality, we assume $0\le x\le y\le 1$. Then 
\begin{align*}
\sup_{|x-y|\le \delta}|\tilde{f}_N(x) - \tilde{f}_N(y)| &= \sup_{|x-y|\le \delta} \bigg|\frac{1}{\sqrt{N}}\sum_{j=1}^N a_j \{\cos(w_j x+u_j)-\cos(w_j y+u_j)\}\bigg| \\
&\le \sup_{\theta\in [0,1]^N}\bigg|\frac{1}{\sqrt{N}}\sum_{j=1}^N a_j\,w_j\sin(w_j\theta_j+u_j)\bigg|\,\delta.
\end{align*}  
The inequality in the preceding holds by applying the mean-value theorem, namely, there exists a sequence $\{\theta_1,\ldots,\theta_N\}$ such that we have $\theta_j\in (x,y)$ satisfying $\cos(w_j x +u_j) - \cos(w_j y+u_j) = -w_j\,\sin(w_j\theta_j+u_j)(x-y)$ for $j=1,\ldots,N$.
To show Equation \eqref{eq:eqi_cont}, it remains to show
\begin{align*}
\limsup_{N\to\infty}P\bigg\{\sup_{\theta\in [0,1]^N}\bigg|\frac{1}{N}\sum_{j=1}^N a_j\,w_j\sin(w_j\theta_j+u_j)\bigg|> \epsilon/\delta\bigg\} < \eta.
\end{align*}
For any fixed $\lambda>0$, recall that $w_j \overset{i.i.d.}\sim \mbox{N}(0, 2/\lambda)$, for $j=1,\ldots,N$. Then $(\lambda/2)\sum_{j=1}^N w_j^2$ is a chi-square random variable with the degree of freedom $N$. Let $K=c\sqrt{2/\lambda}$ with some constant $c \in (1,2)$, then by the sub-exponential tail bound of a chi-square random variable, we have 
\begin{align}\label{eq_P2}
P\bigg(\frac{1}{N}\sum_{j=1}^N w_j^2 > K\bigg) \le \exp(-NK^2/8).
\end{align}
Now define the set $\m A=\{w \in \mb R^N: (1/N)\sum_{j=1}^N w_j^2 \le K\}$ and define the truncated variable $\tilde{w} = w\ind_{\m A}(w)$ over the set $\m A$, the density function of $\tilde{w}$ follows as $\Pi_{\tilde{w}}(\cdot) = \mbox{N}(\cdot;0,2/\lambda)\ind_{\m A}(\cdot)/P(w\in \m A)$. For any fixed $N$, 
\begin{align}\label{eq_cont3}
&P\bigg(\sup_{\theta\in [0,1]^N}\bigg|\frac{1}{\sqrt{N}}\sum_{j=1}^N a_j\,w_j\sin(w_j\theta_j+u_j)\bigg|> \epsilon/\delta\bigg) \nonumber\\
&\le P\bigg(\bigg\{\sup_{\theta\in [0,1]^N}\bigg|\frac{1}{\sqrt{N}}\sum_{j=1}^N a_j\,w_j\sin(w_j\theta_j+u_j)\bigg|> \epsilon/\delta\bigg\} \cap \m A\bigg) +P(\m A^c).
\end{align}
By Equation \eqref{eq_P2}, we see that $\lim_{N\to \infty}P(\m A^c) =0$.  Now we estimate the first term on the right hand side of Equation \eqref{eq_cont3}. First, we consider 
\begin{align*}
 P\bigg(\bigg\{\sup_{\theta\in [0,1]^N}&\bigg|\frac{1}{\sqrt{N}}\sum_{j=1}^N a_j\,w_j\sin(w_j\theta_j+u_j)\bigg|> \epsilon/\delta \bigg\}\cap \m A\bigg) \\
&=  E_{\tilde{w},u}\bigg\{P_{a\mid \tilde{w},u}\bigg(\sup_{\theta\in [0,1]^N}\bigg|\frac{1}{\sqrt{N}}\sum_{j=1}^N a_j\,\tilde{w}_j\sin(\tilde{w}_j\theta_j+u_j)\bigg|> \epsilon/\delta \mid \tilde{w},u \bigg)\bigg\}.
\end{align*}
With fixed $\{\tilde{w_j},u_j\}$, define the set of indexes $J_N = \{1\le j\le N: a_j\tilde{w}_j \ge 0\}$, then 
\begin{align}\label{eq_bound}
\sup_{\theta\in [0,1]^N}\bigg|\frac{1}{\sqrt{N}}\sum_{j=1}^N a_j\,\tilde{w}_j\sin(\tilde{w}_j\theta_j+u_j)\bigg| \le \frac{1}{\sqrt{N}} \bigg| \sum_{j'=1}^{J_N} a_{j'}\,\tilde{w}_{j'}\bigg|.
\end{align}
Then 
\begin{align*}
&P_{a\mid \tilde{w},u}\bigg(\sup_{\theta\in [0,1]^N}\bigg|\frac{1}{\sqrt{N}}\sum_{j=1}^N a_j\,\tilde{w}_j\sin(\tilde{w}_j\theta_j+u_j)\bigg|> \epsilon/\delta \mid \tilde{w},u \bigg) \\\nonumber
&\le P_{a\mid \tilde{w},u}\bigg(\bigg|\sum_{j'=1}^{J_N} a_{j'}\,\tilde{w}_{j'}\bigg|> \epsilon\sqrt{N}/\delta \mid \tilde{w},u \bigg) \\\nonumber
&\le 2 \exp\bigg(-\frac{c\,N \epsilon^2}{2\delta^2\sum_{j'=1}^{J_N} \tilde{w}^2_{j'}}\bigg) \le 2 \exp\bigg(-\frac{c\,\epsilon^2}{2K\delta^2}\bigg),
\end{align*}
where $c>0$ is some constant. 
The first inequality in the preceding applies the bound in Equation \eqref{eq_bound}; the second inequality holds by applying the general Hoeffding's inequality for independent Gaussian random variables; the third inequality is due to the fact that $J_N\le N$ for any fixed $N$. Therefore we have  
$$\limsup_{N\to\infty}P\bigg(\bigg\{\sup_{\theta\in [0,1]^N}\bigg|\frac{1}{\sqrt{N}}\sum_{j=1}^N a_j\,w_j\sin(w_j\theta_j+u_j)\bigg|> \epsilon/\delta\bigg\} \cap \m A\bigg) \le 2 \exp\bigg(-\frac{c\,\epsilon^2}{2K\delta^2}\bigg).$$
Combining the above result with the bound for $P(\m A^c)$, we show that for any $\epsilon,\eta>0$, Equation \eqref{eq:eqi_cont} holds by choosing $\delta =\sqrt{c\epsilon^2/\{K\log(1/\eta)\}}$. Therefore, we have verified that $\widetilde{f}_N$ is asymptotically uniformly eqicontinuous with respect to Euclidean norm. Then we complete the proof of weak convergence of $\widetilde{f}_N$ to the original Gaussian process.  \\

\noindent{\em Part II.} Now we compute the expectation and covariance of $\widetilde{f}_N$. For any $x\in \mathbb{R}$,
\begin{align*}
&E\{\widetilde{f}_N(x)\} \\
&=(2/N)^{-1/2}\sum_{j=1}^N \int\int\,\frac{1}{2\pi}\,\cos(w_j x+s_j)\phi_c(w_j)\,dw_j\,ds_j\\
&= (2/N)^{-1/2} \sum_{j=1}^N \int \frac{1}{2\pi}\,\bigg\{ \cos(w_j x) \int_{-\pi}^{\pi}\cos s_j ds_j -  \sin(w_j x)\int_{-\pi}^{\pi}\sin s_j ds_j \bigg\}\phi_c(w_j)\,d w_j = 0.
\end{align*}
For any $x,y\in \mathbb{R}$,  
\begin{align*}
&\Cov\{\widetilde{f}_N(x),\widetilde{f}_N(y)\}\\
&= (2/N)\sum_{j=1}^N \Cov\{\cos(w_j x+s_j),\cos(w_j x+s_j)\}=2E_{w,s}\cos(xw+s)^2\\
&=\frac{1}{2\pi}\int_w\int_{-\pi}^{\pi} [\cos\{(x+y)w+2s\}+\cos\{(x-y)w\}] \phi_c(w) ds dw\\
&=\frac{1}{2\pi}\int_w\bigg( \int_{-\pi}^{\pi}[\cos\{(x+y)w\}\sin(2s)+\sin\{(x+y)w\}\cos(2s)]\,ds + \cos\{(x-y)w\}\bigg)\phi_c(w)\,dw\\
&=\frac{1}{2\pi}\int_w\cos\{(x-y)w\}\phi_c(w)\,dw =c(x,y).
\end{align*}
   
We now have completed the proof of Theorem \ref{thm:appgp}.

{\tcr \section{Proof of Theorem \ref{thm:rf}}}\label{sec:rf_proof} 
To prove Theorem \ref{thm:rf}, it suffices to verify Conditions \eqref{eq:test}, \eqref{eq:priorconc}, and \eqref{eq:sieve} in Theorem \ref{thm:thm2}. Below we only highlight different steps from the proof of Theorem \ref{thm:mthm}. We use bold letters $\bsm a, \bsm \omega, \bsm s$ to denote the vector form of parameters $\{a_j\}, \{\omega_j\},\{s_j\}$, respectively. 

First we define the sequence of sieves for parameters $\{a_j,\omega_j, s_j\}_{j=1}^N$ in $\wt f_N$ as 
\begin{align}\label{sieve_fn}
\m D_N= \left\{(a_i,\omega_j,s_j): a_j\in [-\sqrt{n}\epsilon_n,\sqrt{n}\epsilon_n],\ \omega_j\in [-n\epsilon_n^2,n\epsilon_n^2],\ s_j\in [0,2\pi],\ j=1,\ldots N. \right\},
\end{align}
for any fixed positive integer $N>0$ and for $\epsilon_n$ defined in Theorem \ref{thm:mthm}. Denote $\wt {\m D}_N = \m D_N\cap \wt {\m A}_N$, where $\wt {\m A}_N = \{(a_i,\omega_j,s_j): \|\wt f_N\|_\infty \le A_0, \ j=1,\ldots N.\}$ for the same constant $A_0$ defined in Assumption 3.

We now start from verifying Condition \eqref{eq:priorconc}, which suffices to show Lemma \ref{lem:gkl} for $\wt f_N$.  Define the set $\wt B_{(f_0, p_0)}(\epsilon_n)$ as the KL--neighborhood of $\wt f_N$ centered at $(f_0,p_0)$, by replacing $f$ with $\wt f_N$ in the definition of $B_{(f_0, p_0)}(\epsilon_n)$ in \eqref{KL_fp}. It suffices to lower bound 
\begin{align}\label{KL_lower}
\Pi \big\{\wt B_{(f_0, p_0)}(\epsilon_n) \big\} \ge \int_{r_0}^{r_1} \Pi \big\{\wt B_{(f_0, p_0)}(\epsilon_n) \mid A=a \big\} g(a) da,
\end{align}
for arbitrary fixed constants $r_0,r_1>0$.
 Then, following a same argument in the proof of Lemma \ref{lem:gkl} leads to
\begin{align*}
\Pi \big\{\wt B_{(f_0, p_0)}(\epsilon_n) \mid A=a \big\} \ge \Pi(\|\wt f_{A,N}-f_0\|_\infty \le \epsilon_n \mid A=a)\,\Pi(KL(p_0,p)\le \epsilon^2_n).
\end{align*}
It suffices to lower bound the first probability term on the right hand side of the preceding. For simplicity, we use the shorthand $\wt f_{a,N}$ for $\wt f_{A=a,N}$.  Fixing $A=a$, recall an original Gaussian process $f\sim \textsc{gp}(0,c^a(\cdot,\cdot))$, and applying the triangular inequality one obtains
\begin{align}\label{fn_prior}
\Pi(\|\wt f_{a,N}-f_0\|_\infty \le \epsilon_n) &\ge \Pi( \|\wt f_{a,N}\|_\infty +\|f\|_\infty +\|f-f_0\|_\infty \le \epsilon_n)\notag\\
&\ge\Pi(\|\wt f_{a,N}\|_\infty \le \epsilon_n/3)\, \Pi(\|f\|_\infty \le \epsilon_n/3)\, \Pi(\|f-f_0\|_\infty \le \epsilon_n/3).
\end{align}
In \cite{van2009adaptive}, it has verified that $\Pi(\|f\|_\infty \le \epsilon_n/3)\gtrsim e^{-c n\epsilon_n^2}$ and $\Pi(\|f-f_0\|_\infty \le \epsilon_n/3)\gtrsim e^{-c' n\epsilon_n^2}$ for some constants $c,c'>0$. Then to bound \eqref{fn_prior}, it suffices to bound 
 \begin{align*}
\Pi(\|\wt f_{a,N}\|_\infty \le \epsilon_n/3) &=\Pi\bigg(\bigg\|(2/N)^{1/2}\sum_{j=1}^{N} a_j \cos (w_j x + s_j)\bigg\|_\infty \le \epsilon_n/3 \bigg) \\
&\ge \Pi \bigg((2/N)^{1/2}\sum_{j=1}^{N} |a_j| \le \epsilon_n/3 \bigg)\ge \Pi \big(\|\bsm a\|^2 \le \epsilon^2_n/18\big),
 \end{align*}
where $\|\bsm a\|^2 =\sum_{j=1}^N a_j^2\sim \chi^2_N$, which is a chi-square random variable with the degree of freedom $N$. Further, we have
\begin{align*}
\Pi \big(\|\bsm a\|^2 \le \epsilon^2_n/18\big) &=\int_{0}^{\epsilon^2_n/18} 2^{-N/2}\{\Gamma(N/2)\}^{-1}x^{N/2-1}\exp(-x/2) dx \\
&\ge (\epsilon_n^2/36)^{N/2} \sqrt{2\pi N} (N/2)^{-N/2}\exp(-\epsilon^2_n/36)\\
& \ge (\epsilon_n/(C\sqrt{N}))^{N} \asymp \exp(-N\log(C\sqrt{N}/\epsilon_n)) \gtrsim\exp(-cn\epsilon_n^2),
\end{align*}
by choosing $N$ such that $N\log(\sqrt{N}/\epsilon_n) \lesssim n\epsilon_n^2$. Notice that the above bound holds uniformly for all $a>0$, then, invoke the above result in \eqref{fn_prior} and \eqref{KL_lower}, we have $\Pi \{\wt B_{(f_0, p_0)}(\epsilon_n)\} \gtrsim e^{-c_1 n\epsilon_n^2}P(r_0<A<r_1) \gtrsim e^{-c_2 n\epsilon_n^2}$. The last inequality holds by choosing $R_n$, based on the result in \cite{van2009adaptive}. Therefore, we have verified Condition \eqref{eq:priorconc}.

Next, we verify Condition \eqref{eq:sieve}, it suffices to show the desired bound for $\Pi(\wt {\m D}_N^c)$. Note that $\Pi(\wt {\m D}_N^c) = \Pi(\m D^c_N\mid \wt {\m A}_N) \le\Pi(\m D_N^c)/\Pi(\wt {\m A}_N)$. Similar to the proof of Theorem \ref{thm:mthm}, it is easy to show that 
\begin{align}\label{A_cpl}
\Pi(\wt {\m A}_N) &\ge \int_{r_0}^{r_1} \Pi(\|\wt f_{a,N}\|_\infty \le A_0) g(a) da\notag\\
 &\ge \Pi\bigg((2/N)^{1/2}\sum_{j=1}^{N} |a_j| \le A_0 \bigg) P(r_0<A<r_1)  \ge \exp\{-N\log(C'\sqrt{N})\},
\end{align}
for some constants $r_0,r_1,c',C'>0$. Next we have 
\begin{align}\label{cpl}
\Pi( {\m D}_N^c) \le 2 \Big\{\Pi_{\bsm a,N}\big(\big\{\cap_{j=1}^N [-\sqrt{n}\epsilon_n,\sqrt{n}\epsilon_n]\big\}^c\big) +  \Pi_{\bsm\omega,N}\big(\big\{\cap_{j=1}^N [-n\epsilon^2_n, n\epsilon^2_n]\big\}^c\big) \Big\}.
\end{align}
Recall $a_j \overset{iid}\sim N(0,1)$ for $j=1,\ldots,N$. First, we can show that 
\begin{align}\label{a_upp}
\Pi_{\bsm a,N}\big(\big\{\cap_{j=1}^N [-\sqrt{n}\epsilon_n,\sqrt{n}\epsilon_n]\big\}^c\big) =  \Pi_{\bsm a,N} \bigg(\max_{1\le j\le N}\,|a_j|\ge\sqrt{n}\epsilon_n\bigg)\notag\\
 \le  \Pi_{\bsm a,N} \bigg(\max_{1\le j\le N}\,|a_j|- \bbE \max_{1\le j\le N}\,|a_j| \ge\sqrt{n}\epsilon_n/2\bigg) \le \exp(-n\epsilon_n^2/8).
\end{align} 
The inequality holds by the known result that $\bbE \max_{1\le j\le N}\,|a_j| \le c \sigma^2_{\max}\sqrt{2\log N}$ for $\{a_j\}$ are independent centered Gaussian random variables with $\sigma^2_{\max} = \max_j\{\Var(a_j)\}=1$. Given the chosen $N$, it is obvious that $ \bbE \max_{1\le j\le N}\,|a_j| < n\epsilon_n^2/2$.  The last inequality uses the tail bound for the maximum of independent Gaussian random variables.  

For any fixed $a>0$, denote $\omega_j|(A=a)$ by $\omega^a_j \overset{iid}\sim N(0, a^2)$ for $j=1,\ldots,N$. Similarly, we can show 
\begin{align}\label{omega_bdd}
&\Pi_{\bsm \omega,N}\big(\big\{\cap_{j=1}^N [-n\epsilon^2_n, n\epsilon^2_n]\big\}^c \mid A=a \big) = \Pi \bigg(\max_{1\le j\le N}\,|\omega^a_j|\ge n\epsilon^2_n\bigg)\notag\\
& \le \Pi \bigg(\max_{1\le j\le N}\,|\omega^a_j|- \bbE \max_{1\le j\le N}\,|\omega^a_j| \ge n\epsilon^2_n/2\bigg) \le \exp\big\{-(n\epsilon_n^2)^2/8a^2\big\}.
\end{align}
The last inequality holds due to facts that $\bbE \max_{1\le j\le N}\,|\omega_j| \le c' a\sqrt{2\log (2N)}$ for some constant $c'>0$ and choosing $R_n \asymp \sqrt{n}\epsilon_n$, and an application of the concentration bound of the maximum of independent Gaussian random variables. 

Then we have for some $R_n>0$ that depends on $n$,
\begin{align*}
\Pi_{\bsm\omega,N}\big(\big\{\cap_{j=1}^N [-n\epsilon^2_n, n\epsilon^2_n]\big\}^c\big) \le \int_0^{R_n} \Pi_{\omega^a,N}\big(\big\{\cap_{j=1}^N [-n\epsilon^2_n, n\epsilon^2_n]\big\}^c\big) g(a)da +P(A>R_n). 
\end{align*}
Based on the final bound in \eqref{omega_bdd}, the second term on the right hand of the preceding can be upper bounded by
\begin{align*}
\int_0^{R_n} \exp\big\{-(n\epsilon_n^2)^2/(8a^2)\big\} g(a)da \le \exp\big\{-(n\epsilon_n^2)^2/(8R_n^2)\big\} \asymp \exp\{-c'n\epsilon^2_n\},
\end{align*}
for some constant $c'>0$. Choosing $R_n \asymp \sqrt{n}\epsilon_n$ leads to the final bound in the preceding. And we have $P(A>R_n) \lesssim e^{-\tilde{c}n\epsilon_n^2}$ for some constant $\tilde{c}>0$ based on the results in \cite{van2009adaptive}. Then invoking these results in \eqref{cpl} combined with result in \eqref{a_upp} leads to the desired result that $\Pi(\wt {\m D}_N^c) \lesssim e^{-c_3n\epsilon_n^2}$ for some constant $c_3>0$, by choosing $N$ such that $N\log N \lesssim n\epsilon_n^2$.  

At last, we verify Condition \eqref{eq:test}. First, we estimate the entropy of the sieves. For arbitrary two parameter vectors $\theta = (a,\omega, s)$ and $\theta'=(a',\omega',s')\in \wt {\m D}_N$, denote $\wt f_N, \wt f'_N$ associated with $\theta,\theta'$, respectively. Here we denote by $a,a'$ two coefficient parameters associated with the Fourier feature functions, rather than the rescaling parameter.   Then 
\begin{align*}
\|\wt f_N - \wt f'_N\|_\infty &=\sqrt{\frac{2}{N}} \Big\|\sum_{j=1}^N\, a_j\cos(\omega_j x+s_j) - \sum_{j=1}^N\, a'_j\cos(\omega'_j x+s'_j)\Big\|_\infty\\
&\le \sqrt{\frac{2}{N}} \bigg[ \Big\|\sum_{j=1}^N \, (a_j-a'_j)\cos(\omega_j x+s_j)\Big\|_\infty \\
&~~~~~~~~~~~~~~~~~ + \Big\|\sum_{j=1}^N \, a'_j \{\cos(\omega_j x+s_j)-\cos(\omega_j x+s_j)\}\Big\|_\infty \bigg]\\
&\le \sqrt{\frac{2}{N}}\big\{ \|a-a'\|_1 + \sqrt{n}\epsilon_n (\|\omega-\omega'\|_1+\|s -s'\|_1) \big\}.
\end{align*}

We now consider the partition $S_n$ of length $\epsilon_n/\sqrt{N}$ on the interval $[-\sqrt{n}\epsilon_n,\sqrt{n}\epsilon_n]$ for each of $\{a_j\}$ for all $j$ and the partition $O_n$ of length $1/\sqrt{nN}$ on the interval $[-n\epsilon_n^2, n\epsilon_n^2]$ for $\{\omega_j\}$ and the partition $M_n$ of length $1/\sqrt{nN}$ on the interval $[0,2\pi]$ for $\{s_j\}$. For any $\wt f_N\in \m D_N$, we can always find $\{a'_j,\omega'_j,s'_j\}$ with $a'_j\in S_n$, $\omega'_j \in O_n$ and $s'_j \in M_n$ for $j=1,\ldots,N$, such that $\wt f'_N (x)= \sqrt{2/N}\sum_{j=1}^N a'_j \cos(\omega_j x +s_j)$ satisfies $\|\wt f_N -\wt f'_N\|_\infty \le \epsilon_n$. By definition, it is obvious that $N(\epsilon_n, \wt{\m D}_N, \|\cdot\|_\infty)\le N(\epsilon_n, \m D_N, \|\cdot\|_\infty)$. Then, it suffices to bound 
\begin{align*}
N(\epsilon_n, \m D_N, \|\cdot\|_\infty) &\le [2\sqrt{n}\epsilon_N/(\epsilon_n/\sqrt{N})+1]^N \times [2n\epsilon_n^2/(\sqrt{nN})^{-1}+1]^N\times [4\pi\sqrt{nN}]^N\\
& \lesssim (n^{3/2}N^{1/2}\epsilon_n^2)^N.
\end{align*}
It is easy to see that $\log N(\epsilon_n, \m D_N, \|\cdot\|_\infty) \lesssim N\log n\lesssim n\epsilon_n^2$ by choosing $N$ such that $N\log n \lesssim n\epsilon_n^2$. Therefore we have verified the entropy condition. 

Based on a same argument of verifying Condition \eqref{eq:test} in the proof of Theorem \ref{thm:mthm}, to complete the proof, it suffices to verify that Lemma \ref{lem:test} holds for all $\wt f_N \in \wt {\m D}_N$. For any $\wt f_N \in \wt {\m D}_N$, it is easy to show that Proposition \ref{bias} holds since $\|\wt f_N\|_\infty\le A_0$ and $\wt f_N$ is infinitely differentiable, which are key points to verify equations \eqref{KL_lower} and  \eqref{fn_prior}). And under the assumption that $\|\wt f_N\|_\infty$, we can verify equations \eqref{test1}, \eqref{test2}, and \eqref{test3}, which completes the proof of Lemma \ref{lem:test}. Putting all pieces together, we have shown Theorem \ref{thm:thm2} for $\wt f_N$, leading to the desired result in Theorem \ref{thm:rf}.  \\

\section{Proof of Auxiliary Results}\label{sec:auxiliary} 
In this section, we provide the proofs of Lemmata \ref{lem:sieve}, \ref{lem:test}, \ref{lem:priorthick} and \ref{lem:gkl} consecutively. 
\subsection{Proof of Lemma \ref{lem:sieve}} 
Based on the definition of sieves $\mathcal{P}_n$, one has $\mathcal{P}_n^c=(B_n^c \otimes \mathcal{F})\cup(B_n \otimes \mathcal{F}^c)\cup(B_n^c \otimes \mathcal{F}^c)$ and $\Pi(\mathcal{P}_n^c)\le 2\{\Pi(B_n^c) +\Pi(\mathcal{F}^c)\}$.  We first bound $\Pi(\mathcal{F}^c)$. Under Assumptions \ref{ass:A} and \ref{ass:sig},
\begin{align*}
\Pi(\mathcal{F}^c) &\le H\bar{\alpha}([-m,m]^c) + P(\widetilde{\sigma} \not\in [\underline{\sigma}, \bar{\sigma}]) + P\bigg(\sum\limits_{h>H} \pi_h >\epsilon\bigg) \\
&\le He^{-b_1 m^{\tau_1}}+ c_2 \bar{\sigma}^{-2\tau_3} +c_1 e^{-b_2\underline{\sigma}^{-2\tau_2}}+\bigg(\frac{e|\alpha|}{H}\log\frac{1}{\epsilon}\bigg)^H.
\end{align*}
Choosing $m^{\tau_1}\lesssim n, \underline{\sigma}\lesssim n^{-1/2\tau_2}$ and $\bar{\sigma}^{2\tau_3}\lesssim e^n$ with $\epsilon = \epsilon_n$ for $\epsilon_n$ defined in Theorem \ref{thm:mthm}, the first three terms on the right hand side of second line in the preceding can be bounded by a multiple of $e^{-n}$, and by taking $H\lesssim n\epsilon^2_n$ the last term in the same line can be bounded from above by,
\begin{align*}
\bigg(\frac{e|\alpha|}{H}\log\frac{1}{\epsilon}\bigg)^H\lesssim e^{-H\log(H\log n)}\lesssim e^{-\frac{1}{2\alpha+1} n^{1/(2\alpha +1)}(\log n)^{2t+1}}\lesssim e^{-c_4 n^{1/(2\alpha +1)}(\log n)^{2t}}.
\end{align*}
Thus $\Pi(\mathcal{F}^c)\lesssim e^{-c_4 n\epsilon_n^2}$ for every $c_4>0$. 

Now we bound $\Pi(\widetilde{B}_n^c)$. By definition, $\Pi(\widetilde{B}_n^c)=\Pi(B_n^c \mid \m{A})\le\Pi(B_n^c)/P(\m{A})$, with $\m{A}$ defined in Assumption \ref{ass:A}. Based on the facts that $E(||f||_{\infty}) < \infty$ and $\sigma^2_f = \sup_{x\in[0,1]} \mbox{E}\{f(x)\}^2 < \infty$, applying Borell's inequality in Lemma \ref{lem:Borell}, we have $P(\m{A})=P(||f||_{\infty} < A_0) \ge 1- e^{-A^2_0/2\sigma^2_f}\ge a_0$, for some $A_0>0$ and $a_0\in(0,1)$. Thus $\Pi(\widetilde{B}_n^c)\lesssim \Pi(B_n^c)\lesssim e^{-n\epsilon_n^2}$ with $M_n^2\lesssim n\epsilon_n^2$ and $a_n^2\lesssim n\epsilon_n^2$. More details can be found in the proof of Theorem 3.1 in \cite{van2009adaptive}.

\subsection{Proof of Lemma \ref{lem:test}}  

To prove Lemma \ref{lem:test}, we will prove the inequality in Equation \eqref{test3} in detail and only mention the key elements in the proof of results in Equations \eqref{test1} and \eqref{test2} since they all follow the similar line of argument. The key elements of the proof are applications of Talagrand's inequality stated in Lemma \ref{lem:Talagrand}, bounded $L_1$-norm of the deconvolution kernel $K_n$ and tight bounds on the bias terms of deconvolution estimators based on the construction in Equations \eqref{eq: pn} and \eqref{eq: fn}. The last two results are stated in the following Proposition \ref{L1norm} and Proposition \ref{bias}.
\begin{prop}\label{L1norm}  
For any kernel function $K$ satisfying conditions in Equation \eqref{eq:kernel1} and $K_n$ defined in Equation \eqref{eq:K_n}, we have $||K_n||_1<C_1$, for some constant $C_1>0$.
\end{prop}
\begin{proof}
There exists a symmetric and integrable kernel function $K$ such that Equation \eqref{eq:kernel1} hold and the Fourier transform $\phi_K(t) = \mathds{1}_{[-1,1]}/(2\pi)$, which is symmetric, real-valued, bounded infinitely smooth function with a compact support. We remark that one example of kernels that satisfy the above conditions is the $sinc$ kernel. For any fixed positive constant $a$, $\int |K_n(s)|\,ds = \int_{|s|\le a} |K_n(s)|\,ds+\int_{|s|>a}|K_n(s)|\,ds$. We have
\begin{align*}
|K_n(s)| \le \int |e^{-its}| \frac{|\phi_K(t)|}{|\phi_{\delta}(t/h_n)|}dt\le \int_{-1}^1 \frac{|\phi_K(t)|}{|\phi_{\delta}(t/h_n)|}dt\lesssim \exp(\delta_n^2 /2h_n^2),
\end{align*}
thus $\int_{|s|\le a} |K_n(s)|\,ds\lesssim \exp(\delta_n^2 /2h_n^2) = O(1)$. For $|s|>a$, by Cauchy-Schwarz inequality,
\begin{align*}
\int_{|s|>a}|K_n(s)|ds \le \bigg(\int_{|s|>a} \frac{1}{s^4}ds\bigg)^{1/2}\bigg\{\int_{|s|>a} s^4K_n(s)^2 ds\bigg\}^{1/2}.
\end{align*}
By Parseval's theorem, $\int \{s^2 K_n(s)\}^2\,ds=\int{\{g''(t)\}^2}\,dt$ with
\begin{align*}
g(t) = \phi_K(t)/\phi_{\delta}(t/h_n) =  \frac{1}{2\pi} e^{-t^2\delta^2/(2h_n^2)}~\mathds{1}_{[-1,1]}.
\end{align*}
Since $g''(t)$ is the Fourier transform of $(is)^2K_n(s)$, also $g(t), g'(t), g''(t)$ are continuous and therefore bounded on $[-1,1]$. Thus $\int \{s^2 K_n(s)\}^2\,ds$ is bounded and so is $\int_{|s|>a} 1/s^4\,ds$, which yields the result that $\int |K_n(s)|\,ds$ is bounded. 
\end{proof}

The following Proposition provides tight bounds on the bias terms of $\widehat{p}_n$ and $\widehat{f}_n\widehat{p}_n$ separately.
\begin{prop}\label{bias} 
For $\widehat{p}_n$ and $\widehat{f}_n$ defined in Equation \eqref{eq: pn} and Equation \eqref{eq: fn} and for any $f,p \in \m{P}_n$ we have
\begin{eqnarray*}
|| E_{W,X} (\widehat{p}_n) - p||_1\lesssim \epsilon_n, \quad \text{and} \quad || E_{Y,W,X} (\widehat{f}_n\widehat{p}_n) - fp||_1\lesssim \epsilon_n,
\end{eqnarray*}
 with $\epsilon_n$ defined in Theorem \ref{thm:mthm}.
\end{prop}

\begin{proof}
By Fourier inversion theorem, it is easy to show that $E_{W,X} (\widehat{p}_n) = K_{h_n}*p(x)$ and $E_{Y,W,X} (\widehat{f}_n\widehat{p}_n) = K_{h_n}* (fp)$ with $K_{h_n}=K(\cdot/h_n)/h_n$. First for any $p = \phi_{\widetilde{\sigma}}*F$, by Cauchy-Schwarz inequality we have  $|| K_{h_n}*p - p||_1 \le || K_{h_n}*p - p||_2$.
Recall the Fourier transform of the kernel function $K$ is denoted by $\phi_K(t)$, applying Parseval's theorem again,
\begin{align*}
|| K_{h_n}*p - p||_2^2 &= \int |2\pi\phi_K (h_n t) -1|^2 |\widehat{p}(t)|^2 dt=\int_{|t|>1/h_n} |\widehat{F}(t)|^2|\widehat{\phi}_{\widetilde{\sigma}}(t)|^2 dt\\
 &\le \int_{|t|>1/h_n} |\widehat{\phi}_{\sigma}(t)|^2 dt \le (h_n/\underline{\sigma}^2) e^{-(\underline{\sigma}/h_n)^2/2} \\
 &\lesssim h_n^{-1} (\log n)^{-t_3} e^{-K^2(\log n)^{2t_3}/2}\lesssim \epsilon_n^2,
\end{align*}
for all $\widetilde{\sigma} \ge \underline{\sigma}$. Let $h_n \asymp \epsilon_n^{1/\beta}$ with $\epsilon_n$ defined in Theorem \ref{thm:mthm} and by Lemma \ref{lem:sieve} we have $\underline{\sigma}\lesssim n^{-1/(2\tau_2)}$, where $\tau_2$ is chosen such that $\underline{\sigma} = K h_n (\log n)^{t_3}$ for some constants $K,t_3$ satisfying $K^2/2>1$ and $t_3>1/2$.

Now we bound the bias term of $\widehat{f}_n\widehat{p}_n$. By triangle inequality,
\begin{align}\label{eq:bias}
|| K_{h_n}*(fp) - fp||_1 \le || K_{h_n}*(fp) - p  K_{h_n}*f ||_1 + ||p  K_{h_n}*f  - fp||_1.
\end{align}
By Cauchy-Schwarz inequality, the first term of the right hand side of Equation \eqref{eq:bias} can be bounded as 
\begin{align} \label{bdd_L1}
||  K_{h_n}*(fp) - p  K_{h_n}*f ||_1 &= \int\int | K_{h_n}(x-y) \{p(y) - p(x)\}f(y)\,dy|\,dx \notag \\
&\le || K_{h_n}*p - p||_2 \,||f||_2\lesssim || K_{h_n}*p - p||_2,
\end{align}
since $||f||_2 \le ||f||_{\infty} \le A_0$ under Assumption \ref{ass:A}.
The second term on the right hand side of Equation \eqref{eq:bias} can be bounded
\begin{align}\label{bdd_sup}
||p  K_{h_n}*f  - fp||_1 \le ||p||_1 ||K_{h_n}*f  - f||_{\infty} = ||K_{h_n}*f  - f||_{\infty}\lesssim \epsilon_n.
\end{align}
The last inequality in the preceding holds based on the properties of higher order kernel as in Lemma 4.3 of \cite{van2009adaptive}.  
\end{proof}

\noindent{\em Proof of Equation \eqref{test3}.} Now we are ready to prove the inequality in Equation \eqref{test3}. By triangle inequality,
\begin{align}\label{I_terms}
 ||\widehat{f}_n\widehat{p}_n - fp||_1 &\le ||\widehat{f}_n\widehat{p}_n-E_{Y,W|X}(\widehat{f}_n\widehat{p}_n)||_1
+||E_{Y,W|X}(\widehat{f}_n\widehat{p}_n)-E_{Y,W,X}(\widehat{f}_n\widehat{p}_n)||_1 \nonumber\\
&+ ||E_{Y,W,X}(\widehat{f}_n\widehat{p}_n)-f\cdot p||_1:= I_{1,n}+I_{2,n}+I_{3,n}.
\end{align}
First we estimate $P(I_{1,n}>\epsilon_n/2)$ for $I_{1,n}$ in Equation \eqref{I_terms}. By definition,
\begin{align}\label{T_terms}
&\widehat{f}_n\widehat{p}_n-E_{Y,W|X}(\widehat{f}_n\widehat{p}_n)\nonumber\\
 &=\frac{1}{2\pi nh_n}\sum_{j=1}^{n} \int e^{-\frac{itx}{h_n}}\bigg\{e^{itW_j/h_n}Y_j - E_{W|X}\Big(e^{itW_j/h_n}\Big) E_{Y|X}(Y_j)\bigg\} \frac{\phi_K(t)}{\phi_{u}(t/h_n)}dt\nonumber\\
&= \frac{1}{2\pi nh_n}\sum_{j=1}^{n} \int  e^{-\frac{it(x-W_j)}{h_n}} \frac{\phi_K(t)}{\phi_{u}(t/h_n)}\,dt\ \{Y_j-E_{Y|X}(Y_j)\}\nonumber\\
&+\frac{1}{2\pi nh_n}\sum_{j=1}^{n} \int e^{-\frac{itx}{h_n}}\bigg\{e^{itW_j/h_n}-E_{W|X}\Big(e^{itW_j/h_n}\Big)\bigg\}\frac{\phi_K(t)}{\phi_{u}(t/h_n)}\,dt \ E_{Y|X}(Y_j)\nonumber\\
&:=T_{1,n}+T_{2,n}.
\end{align}
First, we estimate $P(||T_{2,n}||_1>\epsilon_n/2)$ with $T_{2,n}$ defined in Equation \eqref{T_terms}. By Hahn-Banach Theorem, there exists a bounded linear functional $T$ such that $T(h) = \int T_{2,n}(x) h(x)dx$ for all $h\in L_{\infty}[0,1]$, namely, for all $h(x)$ such that $\sup_{x\in[0,1]}|h(x)|<\infty$. And $||T_{2,n}||_1 = ||T||_{\m F_1}$ where $||T||_{\m F_1} = \sup_{h\in \m F_1}|T(h)|$ and $\m{F}_1$ is a countable and dense subset of $L_{\infty}[0,1]$. Thus we have
\begin{align*}
\m{K}=\bigg\{k(u,v): (u,v)\mapsto \frac{1}{h_n}\int_0^1 \bigg[K_n\bigg(\frac{x-u}{h_n}\bigg)-E_{W|X} &\bigg\{K_n\bigg(\frac{x-W}{h_n}\bigg)\bigg\}\bigg]f(v)h(x)dx, \\
&\text{for all}\ h\in\m{F}_1\bigg\},
\end{align*}
and $||nT_{2,n}||_1 = \sup_{k\in\m{K}} |\sum_{j=1}^n k(W_j, X_j)|$. To apply Lemma \ref{lem:Talagrand}, we need to estimate the following quantities, $\sup_{k\in\m{K}}||k(u,v)||_{\infty}$, $\sigma^2_{\m{K}}=E_{W|X}\{\sup k^2(W,X)\}$ and $E\{\sup_{k\in\m{K}}k(W,X)\}$. Based on the Assumptions \ref{ass:A} and \ref{ass:sig} we have $||f||_{\infty}\le C_0$ and $||h||_{\infty}\le 1$, then for any $k\in\m{K}$,
\begin{align*}
|k(u,v)|\le \frac{C_2}{h_n} \bigg[\int_0^1 \bigg|K_n\bigg(\frac{x-u}{h_n}\bigg)\bigg|dx+\int_0^1\bigg|E_{W|X} \bigg\{K_n\bigg(\frac{x-W}{h_n}\bigg)\bigg\}\bigg|dx \bigg],
\end{align*}
for some constant $C_2>0$. For any $u$, by change of variables $s=(x-u)/h_n$, for any fixed positive constant $a$, one has $\int_0^1|K_n\{(x-u)/h_n\}/h_n | ~dx \le \int |K_n(s)|ds \le C'$ for some constant $C'$. The second inequality holds by Proposition \ref{L1norm}. Given $W\mid X\sim \mbox{N}(X,\delta^2_n)$,
\begin{align*}
E_{W\mid X}\bigg\{K_n\bigg(\frac{x-W}{h_n}\bigg)\bigg\}&=\frac{1}{2\pi}\int E_{W\mid X}\bigg(e^{-it[\{x-X-(W-X)\}/h_n]}\bigg)\frac{\phi_K(t)}{\phi_u(t/h_n)}dt\\
&=\frac{1}{2\pi}\int e^{-it(x-X/h_n)}\phi_K(t)dt = K\{(x-X)/h_n\}.
\end{align*}
Again by change of variables $r=(x-X)/h_n$, we have $\int_0^1 E_{W|X}[K\{(x-W)/h_n\}/{h_n}]\,dx = \int|K(r)|dr=1$. There exists a constant $K_1$ such that $||k||_{\infty} \le K_1$ for any $k\in \m{K}$, then $\sup_{k\in\m{K}}\|k\|_{\infty}\lesssim \max\{1,\exp(\delta_n^2/2h_n^2)\}$. Next we estimate the term $\sigma^2_{\m{K}}$. For any $k\in\m{K}$ and $W\mid X \sim \mbox{N}(X,\delta^2_n)$,
\begin{align*}
k(W,X)^2 &=\frac{1}{h_n^2} \bigg(\int_0^1 \bigg[K_n\bigg(\frac{x-u}{h_n}\bigg)-E_{W|X} \bigg\{K_n\bigg(\frac{x-W}{h_n}\bigg)\bigg\}\bigg]f(X)h(x)dx\bigg)^2\\
&\lesssim \frac{1}{h_n^2} \bigg\{\int_0^1 K_n\bigg(\frac{x-u}{h_n}\bigg)dx\bigg\}^2+\frac{1}{h_n^2} \bigg\{\int_0^1E_{W|X} K_n\bigg(\frac{x-W}{h_n}\bigg)dx\bigg\}^2  \\
&\lesssim \max\{1,\exp(\delta_n^2 /h_n^2)\}.
\end{align*}
Therefore $\sup_{k\in\m{K}}E_{W|X}\{k(W,X)^2\}\lesssim \max\{1,\exp(\delta_n^2 /h_n^2)\}$.

Finally, we move to bound $E_{W|X}(\sup_{k\in\m{K}}|\sum_{j=1}^n k(W_j,X_j)|)$. By Cauchy-Schwarz inequality,
\begin{align*}
E_{W\mid X} &\bigg(\sup_{k\in\m{K}}\bigg|\sum_{j=1}^n k(W_j,X_j)\bigg|\bigg) \\
&\le\bigg[\,E_{W\mid X} \bigg\{\sup_{k\in\m{K}}\bigg|\sum_{j=1}^n k(W_j,X_j)\bigg|\bigg\}^2\,\bigg]^{1/2}\\
&\lesssim \bigg(\frac{1}{h_n^2} \sum_{j=1}^n E_{W\mid X}\bigg[\int \bigg|\,K_n\bigg(\frac{x-W_j}{h_n}\bigg)-E_{W\mid X} \bigg\{K_n\bigg(\frac{x-W_j}{h_n}\bigg)\bigg\}\bigg|\,dx\,\bigg]^2\,\bigg)^{1/2}\\
&\lesssim n^{1/2} \max\{1,\exp(\delta_n^2 /2h_n^2)\}.
\end{align*}
To apply the Lemma \ref{lem:Talagrand}, we choose $\delta_n\asymp h_n$ and same $\epsilon_n$ in Theorem \ref{thm:mthm}, we have $\exp(\delta^2_n / 2h_n^2)=O(1)$. By choosing $t=n\epsilon_n^2$, we have $n^{1/2}+\{2(n+n^{1/2})n\epsilon^2_n\}^{1/2}+ n\epsilon_n^2/3\lesssim n\epsilon_n.$

We now discuss bounding the probability $P(\|T_{1,n}\|_1>\epsilon_n/2)$ with $T_{1,n}$ defined in Equation \eqref{T_terms}. Recall that 
$$nT_{1,n}=\sum_{j=1}^nK_n\{(x-W_j)/h_n\}(Y_j-E_{Y\mid X}Y_j)/h_n=\sum_{j=1}^nK_n\{(x-W_j)/h_n\}\widetilde{Y}_j/h_n,$$ 
with $\widetilde{Y}_j\sim \mbox{N}(0,1)$ i.i.d. for $j=1,\ldots,n$, given $Y_j\mid X_j\sim \mbox{N}(f(X_j),1)$ for $j=1,\ldots,n$. Again by Hahn-Banach theorem, there exists a countable and dense subset $\m{T} \in L_{\infty}[0,1]$ and a class of bounded linear functionals on $L_{\infty}[0,1]$,
$$
\m{Q}=\bigg\{q = \sum_{j=1}^n \widetilde{q}(u_j),\ \widetilde{q}(u) = \int_0^1 \sum_{j=1}^nK_n\bigg(\frac{x-u}{h_n}\bigg)(Y_j-E_{Y\mid X}Y_j)\,t(x)\,dx, \ t\in \m{T}\bigg\},
$$
and $\|nT_{1,n}\|_1=\sup_{q\in \m{Q}}\|q\|_{\infty}$. 

We now proceed to estimate $\sigma^2_{\m{Q}}=\sup_{q\in\m{Q}}E_{Y\mid X}\{\sum_{j=1}^n \widetilde{q}(W_j)\}^2$ and $E_{Y\mid X}(\sup_{q\in \m{Q}}\|q\|_{\infty})$ in order to apply Lemma \ref{lem:Borell}. We first estimate $\sigma^2_{\m{Q}}$. Again, by change of variables and the fact $\|t\|_{\infty}\le 1$ we have
\begin{align*}
E_{Y\mid X}\bigg\{\sum_{j=1}^n\widetilde{q}(W_j)\bigg\}^2 &=\frac{1}{h_n^2}\sum_{j=1}^n\bigg\{\int_0^1 K_n\bigg(\frac{x-W_j}{h_n}\bigg)\,t(x)\,dx\,\bigg\}^2 \\
&\le \frac{1}{h_n^2}\sum_{j=1}^n\bigg\{\int_0^1 K_n\bigg(\frac{x-W_j}{h_n}\bigg)\,dx\,\bigg\}^2\\
&\le \sum_{j=1}^n\bigg(\int |K_n(u)|du\bigg)^2\lesssim n\max\{1,\exp(\delta_n^2 /h_n^2)\}.
\end{align*}
Next we estimate $E_{Y\mid X} (\sup_{q\in\m{Q}}\|q\|_{\infty})$, using the generalized Minkowski inequality, we obtain
\begin{align*}
E_{Y\mid X}\sup_{q\in\m{Q}}\|q\|_{\infty}&=E_{Y\mid X}(\|nT_{1,n}\|_1)\le \{E_{Y\mid X}(\|nT_{1,n}\|_1^2)\}^{1/2}\\
&\le \|[E_{Y\mid X}\{(nT_{1,n})^2\}]^{1/2}\,\|_1 \\
&= \int\bigg\{\,\frac{1}{h_n^2}\sum_{j=1}^nK_n\bigg(\frac{x-W_j}{h_n}\bigg)^2\,\bigg\}^{1/2}\,dx.
\end{align*}
The last equation in the preceding holds because $Y_j$'s are independent. By Jensen's inequality and change of variables it can be bound by
$\sum_{j=1}^n \int K_n\{(x-W_j)/h_n\}^2 dx\}^{1/2}/h_n=n^{1/2}\{\int K_n(u)^2 du\}^{1/2}/h_n$. Fixed any constant $a'>0$, one has $$\int K_n(u)^2du\le\int_{|u|>a'} (u^4/a'^4)\, K_n(u)^2\,du + \int_{|u|\le a'}K_n(u)^2du.$$ It has been shown in the proof of Proposition \ref{L1norm} that $\int u^4 K_n(u)^2du\lesssim \exp(\delta_n^2 /h_n^2)$, and it is easy to see that $\int K_n(u)^2 du\lesssim \max\{1,\exp(\delta_n^2 /h_n^2)\}$. Thus we have $E_{Y\mid X}(\sup_{q\in\m{Q}}\|q\|_{\infty})\lesssim n^{1/2}\,\max\{1,\exp(\delta_n^2 /h_n^2)\}/\sqrt{h_n}$.  Then applying Borell's inequality in Lemma \ref{lem:Borell} by choosing $x=n\epsilon_n$, $\delta_n\asymp h_n\asymp \epsilon _n^{1/\beta}$, where $\epsilon_n$ is defined in Theorem \ref{thm:mthm}, we have shown that $P(\|T_{1,n}\|_1>\epsilon_n/2)<e^{-n\epsilon_n^2/8}$.

We now estimate the probability $P(I_{2,n}>\epsilon_n/2)$, recall that $I_{2,n}$ is defined in Equation \eqref{I_terms}. By definition, $I_{2,n}=E_{Y,W\mid X} (\widehat{f}_n\widehat{p}_n)-E_{Y,W,X}(\widehat{f}_n\widehat{p}_n)$, then with simple calculation one can show that $E_{Y,W\mid X}(\widehat{f}_n\widehat{p}_n) = \sum_{j=1}^n K\{(x-X_j)/h_n\}f(X_j)/(n h_n)$. Similarly, by Hahn-Banach theorem, there exists a countable and dense set $\m{H}_1\in L_{\infty}[0,1]$ such that we can construct a class of bounded linear functionals
$$
\m{L}=\bigg\{l(u): u \mapsto\int \bigg[K\bigg(\frac{x-u}{h_n}\bigg)f(u)-E_X\,\bigg\{K\bigg(\frac{x-X}{h_n}\bigg)f(X)\bigg\}\bigg]\,h_1(x)\,dx,\quad h_1 \in\m{H}_1 \bigg\},
$$
and we have $\|nI_{2,n}\|_1=\sup_{l\in\m{L}}\|\sum_{j=1}^n l(X_j)\|_{\infty}$. To apply the Talagrand's inequality, we first bound
 $\sup_{l\in\m{L}}|\l(u)\|_{\infty}\le (\int |K\{(x-X_j)/h_n\}/ h_n|~dx)\|f\|_{\infty}$. Since $\int|K(u)|du\le K_3$ for some constant $K_3>0$, by change of variables and Assumption \ref{ass:A} one can show $\sup_{l\in\m{L}}\|\l(u)\|_{\infty} \le K_4$, for some constant $K_4>0$.

Second, we bound 
$\sup_{l\in\m{L}}E_X\{l(X)\}^2$. For any $l\in\m{L}$,
\begin{align*}
E_X\{l(X)\}^2 &\le 2E_X\bigg( \bigg\{\int \bigg|K\bigg(\frac{x-X}{h_n}\bigg)\bigg|dx\bigg\}^2+
 2\bigg[\int E_X \bigg\{K\bigg(\frac{x-X}{h_n}\bigg)\bigg\}dx\,\bigg]^2\bigg)\,\|f\|_{\infty}^2/h_n^2\\
  &\le K_5,
\end{align*}
for some constant $K_5>0$. Thus we show that $\sup_{l\in\m{L}}E_X\{l(X)^2\}\le K_5$.

At last, we have
\begin{align*}
&E_{X}\sup_{l\in\m{L}}\bigg|\sum_{j=1}^n l(X_j)\bigg| \\
&\le \bigg\{E_X\bigg(\sup_{l\in\m{L}}\bigg|\sum_{j=1}^n l(X_j)\bigg|\bigg)^2\bigg\}^{1/2}\\
&\le \frac{1}{h_n} \bigg(2n \bigg[E_X\bigg\{\int K\bigg(\frac{x-X}{h_n}\bigg)dx\bigg\}^2 + \bigg\{\int E_X K\bigg(\frac{x-X}{h_n}\bigg)dx\bigg\}^2 \bigg] \bigg)^{1/2} \|f\|_{\infty} \\
&\lesssim ({n}/{h_n})^{1/2}.
\end{align*}
Choosing $h_n\asymp \epsilon_n^{1/\beta}$ with $\epsilon_n$ defined in the Theorem \ref{thm:mthm}, then applying Talagrand's inequality yields the result $P(I_{2,n} >\epsilon_n/2) \le e^{-n\epsilon_n^2/8}$.

Finally, for $I_{3,n}$ defined in Equation \eqref{I_terms}, it is easy to see $I_{3,n} \le \epsilon_n$ by Proposition \ref{bias}. Combining the results of $I_{1,n}, I_{2,n}$ and $I_{3,n}$, we prove the inequality in Equation \eqref{test3}. \\

\noindent{\em Proof of Equation \eqref{test2}.} Inequality in Equation \eqref{test2} can be obtained directly from  Equation \eqref{test3}, as it can be seen as a special case of Equation \eqref{test3} by letting the regression function $f(x)\equiv c$ for some constant $c>0$.\\

\noindent{\em Proof of Equation \eqref{test1}.} The proof of inequality \eqref{test1} follows a same line of arguments in the proof of Equation \eqref{test3} and we omit some details. Let $P_{1,n}=\widehat{p}_n-E_{W\mid X}(\widehat{p}_n)$, $P_{2,n}=E_{W\mid X}(\widehat{p}_n)-E_{W,X}(\widehat{p}_n)$ and $P_{3,n}=E_{W,X}(\widehat{p}_n) - p$. First, we estimate $P(\|P_{1,n}\|_{\infty} > \epsilon_0/2)$. The difference is that we consider the empirical process directly in $\|\cdot\|_{\infty}$. Since the function $K_n(x)$ is continuous and bounded on $[0,1]$, by the separability of $C[0,1]$,  there exists a countable and dense set $T$ over $[0,1]$ and consider the class,
$$
\m{M}=\bigg\{m_x(u): u \mapsto\int e^{-itx/h_n} \bigg\{e^{itu/h_n}-E_{W|X} \bigg(e^{itW/h_n}\bigg)\bigg\} \frac{\phi_{K}(t)}{\phi_u (t/h_n)}dt,\ x\in T \bigg\},
$$
then $\|nP_{1,n}\|_{\infty}=\sup_{x\in T}|\sum_{j=1}^n m_x(W_j)|$. Also we  can show
\begin{align*}
&\sup_{x\in T}\|m_x\|_{\infty}\lesssim h_n^{-1}\exp( \delta_n^2/2h_n^2),\\
&\sup_{x\in T}E_{W|X}[m_x (W)]^2\lesssim h_n^{-2}\exp(\delta_n^2/h_n^2),\\
&E_{W|X} \sup_{x\in T} |\sum_{j=1}^n m_x (W_j)|\lesssim n^{1/2} h_n^{-1}\exp( \delta_n^2/h_n^2).
\end{align*}
Therefore choosing $\delta_n=o(h_n)$ and $h_n =o(\epsilon_n)$ with same $\epsilon_n$ in Theorem \ref{thm:mthm}. For any $\epsilon_0 >0$, take $t=\epsilon_0 n h_n^2$, one has
\begin{align*}
n^{1/2}h_n^{-1}&\exp(\delta_n^2/2h_n^2) \\
&+ \{2nh_n^{-2}\exp( \delta_n^2/h_n^2) +4{n}^{1/2}h_n^{-1}\exp( \delta^2/2h_n^2)\}^{1/2}\,(n\epsilon_0 h_n^2)^{1/2} + \epsilon_0 n h_n^2 \\ 
&< n\epsilon_0.
\end{align*}
By applying Lemma \ref{lem:Talagrand}, one can show $P(\|\widehat{p}_n-E_{W\mid X} (\widehat{p}_n)\|_{\infty}>\epsilon_0) \le e^{-\epsilon_0 nh_n^2}$. Similarly, for $P_{2,n}$ one can write $P_{2,n}=E_{W\mid X} (\widehat{p}_n) - E_{W,X}(\widehat{p}_n) = \sum_{j=1}^n \widetilde{g}_x (X_j)/(nh_n)$,  where $\widetilde{g}_x(u) =K\{(x-u)/h_n)\} - E_X [K\{(x-X)/h_n\}]$ for any $x\in T$. Construct the class $\m{G}=\{\widetilde{g}_x, x\in T\}$ with the countable and dense set $T$ over $[0,1]$, with same calculation by choosing $t=n\epsilon_0h_n^2$, $\delta_n=o(h_n)$ and $h_n =o(\epsilon_n)$, another application of Talagrand's inequality shows $P(P_{2,n}>\epsilon_0) \le e^{-\epsilon_0 nh_n^2}$. Combining the above results for $P_{1,n}$, $P_{2,n}$ and applying Proposition \ref{bias} to $P_{3,n}$ completes the proof of Equation \eqref{test1}.

\subsection{Proof of Lemma \ref{lem:priorthick}} 
The Kullback--Leibler neighborhood around $f_0$ has been studied extensively in Baysian literature. We give a brief argument mentioning the difference in our case, refer to \cite{shen2013adaptive} for extended proof. Under the Assumption \ref{ass:p0}, $p_0$ is compactly supported and lower-bounded. From Theorem 3 in \cite{shen2013adaptive}, there exists a density function $h_{\sigma}$ supported on $[-a_0, a_0]$ satisfying $H(p_0,\phi_{\sigma}*h_{\sigma})\lesssim\sigma^{\beta}$, for some constant $a_0>0$. Fix $\sigma^{\beta} = \widetilde{\epsilon}_n\{\log(1/\widetilde{\epsilon}_n)\}^{-1}$ and find $b'>\max{(1,1/(2\beta))}$ such that $\widetilde{\epsilon}_n^{b'}\{\log(1/\widetilde{\epsilon}_n)\}^{5/4} \le \widetilde{\epsilon}_n$. By Lemma 2 of \cite{ghosal2007posterior} there is a discrete probability measure $F'=\sum_{j=1}^N p_j \delta_{z_j}$ with at most $N\le D\sigma^{-1}\{\log(1/\sigma)\}^{-1}$ support points on $[-a_0, a_0]$, and $F'$ satisfies $H(\phi_{\sigma}*h_{\sigma}, \phi_{\sigma}*F')\le \widetilde{\epsilon}_n^{b'}\{\log(1/\widetilde{\epsilon}_n)\}^{1/4}$. We construct the partition $\{U_1, \ldots, U_M\}$ in the flavor of $c\sigma\widetilde{\epsilon}^{b'}_n\le\alpha(U_j)\le1$ for $j=1, \ldots, M$, where $M\lesssim\widetilde{\epsilon}_n^{1/\beta}\{\log(1/\widetilde{\epsilon}_n)\}^{1+1/\beta}$. Further denote the set $S_F$ of probability measure $F$ with $\sum_{j=1}^M |F(U_j)-p_j|\le 2\widetilde{\epsilon}_n^{2b'}$ and $\min_{1\le j \le M}{F(U_j)} \ge \widetilde{\epsilon}_n^{4b'}/2$ for sufficiently large $n$. Then $\Pi(S_F)\gtrsim \exp[-\widetilde{\epsilon}_n^{-1/\beta}\{\log(1/\widetilde{\epsilon}_n)\}^{2+1/\beta}]$. For each $F \in S_F$,
\begin{align*}
H(p_0,p_{F,\sigma})&\le H(p_0,\phi_{\sigma}*h_{\sigma}) + H(\phi_{\sigma}*h_{\sigma},\phi_{\sigma}*F')+H(\phi_{\sigma}*F',p_{F,\sigma})\\
&\lesssim \sigma^{\beta}+\widetilde{\epsilon}_n^{b'}\{\log(1/\widetilde{\epsilon}_n)\}^{1/4} +\widetilde{\epsilon}_n^{b'}\lesssim \sigma^{\beta}.
\end{align*}
Also we can show that for every $x\in[-a_0,a_0]$, $p_{F,\sigma}/p_0 \ge A_4 \widetilde{\epsilon}_n^{b'}/\sigma$ for some constant $A_4>0$, which leads to $\log\|p_0/p_{F,\sigma}\|_{\infty}\lesssim \log(1/\widetilde{\epsilon}_n) $.

\subsection{Proof of Lemma \ref{lem:gkl}} 
To prove Lemma \ref{lem:gkl}, by the definition of the Kullback--Leibler neighborhood defined in Equation \eqref{KL_fp}, it suffices to bound the Kullback--Leibler divergence and the second moment of Kullback--Leibler divergence between $g_{f_0,p_0}$ and $g_{f,p}$ from above, respectively. Based on Lemma 5.3 in \cite{van2009adaptive} and Lemma \ref{lem:sieve} in Appendix~\ref{sec:aux}, we have $\Pi \{KL(p_0,p) \le \epsilon_n^2\} \ge e^{-n\epsilon_n^2}$ and $ \Pi(\|f-f_0\|_{\infty} < \epsilon_n) \ge e^{-n\epsilon_n^2}$. Then using the convexity of the Kullback--Leibler divergence with respect to both arguments, we have
\begin{align*}
&KL(g_{f_0,p_0}, g_{f,p}) \\
&= KL\bigg(\frac{1}{2\pi \delta_n}\int e^{-\frac{1}{2} \{y-f_0(x)\}^2-\frac{1}{2\delta_n^2} (w-x)^2} dP_0,\, \frac{1}{2\pi \delta_n}\int e^{-\frac{1}{2} \{y-f(x)\}^2-\frac{1}{2\delta_n^2} (w-x)^2} \frac{p}{p_0} \,dP_0\bigg) \\
&\le \int KL \bigg( \frac{1}{2\pi \delta_n} \,e^{-\frac{1}{2} (y-f_0(x))^2-\frac{1}{2\delta_n^2} (w-x)^2}, \,\frac{1}{2\pi \delta_n}\, e^{-\frac{1}{2} \{y-f(x)\}^2-\frac{1}{2\delta_n^2} (w-x)^2}\, \frac{p}{p_0}\bigg) d P_0\\
&= \int \int \frac{1}{2\pi \delta_n}\, e^{-\frac{1}{2} \{y-f_0(x)\}^2-\frac{1}{2\delta_n^2} (w-x)^2} \log \bigg( \frac{e^{-\frac{1}{2} \{y-f_0(x)\}^2}}{e^{-\frac{1}{2} \{y-f(x)\}^2}} \,\frac{p_0}{p} \bigg)\, dy\,dw\, dP_0\\
&= \int [KL\{\mbox{N}(y; f_0,1), \mbox{N}((y; f,1)\} + \log(p_0/p )] \,dP_0\\
&\lesssim \|f_0 - f\|_{\infty}^2 + KL(p_0,p)\lesssim \epsilon_n^2,
\end{align*}
where $P_0$ denotes the distribution measure associated with $p_0$. Next, we decompose the second moment of the Kullback--Leibler divergence into,
\begin{align}\label{eq:denom}
 \int g_{f_0,p_0}\bigg(\log\frac{g_{f_0,p_0}}{g_{f,p}}\bigg)^2 &=  \int_{A_n} g_{f_0,p_0}\bigg(\log\frac{g_{f_0,p_0}}{g_{f,p}}\bigg)^2 +  \int_{A_n^c} g_{f_0,p_0}\bigg(\log\frac{g_{f_0,p_0}}{g_{f,p}}\bigg)^2\nonumber\\
 &=: I_1 + I_2,
\end{align}
where $A_n=\{y\in \mathbb{R}:|y|\le \gamma'/\epsilon_n\}$ for some constant $\gamma'>0$.\\
We first bound term $I_1$ in Equation \eqref{eq:denom}, apply the inequality
\begin{align*}
 \int_{A_n} g_{f_0,p_0}\bigg(\log\frac{g_{f_0,p_0}}{g_{f,p}}\bigg)^2 \le 2 H^2( g_{f_0,p_0}, g_{f,p}) (1 + \log\| (g_{f_0,p_0}/g_{f,p}) \mathds{1}_{A_n}\|_{\infty})^2.
\end{align*}
It is well known that $H^2( g_{f_0,p_0}, g_{f,p}) \le KL (g_{f_0,p_0}, g_{f,p})$, then to estimate $I_1$ it remains to estimate the term $\| (g_{f_0,p_0}/g_{f,p})\, \mathds{1}_{A_n}\|_{\infty}$. By definition,
\begin{align*}
&\bigg|\frac{g_{f_0,p_0}(y,w)}{g_{f,p}(y,w)}\bigg|\,\mathds{1}_{A_n} \\
&\le \bigg| \frac{\int_{A_n} e^{-\frac{1}{2} (y-f_0(x))^2}e^{-\frac{1}{2\delta_n^2} (w-x)^2}p(x)dx}{\int_{A_n} e^{-\frac{1}{2} (y-f_0(x))^2} [ e^{-\frac{1}{2} (y-f(x))^2} /e^{-\frac{1}{2} (y-f_0(x))^2} ]\, e^{-\frac{1}{2\delta_n^2} (w-x)^2}p(x)\,dx }\bigg|\cdot \bigg|\bigg|\,\frac{p_0}{p}\bigg|\bigg|_{\infty}\\
&\le \bigg|\bigg|\frac{e^{-\frac{1}{2} (y-f_0)^2}}{e^{-\frac{1}{2} (y-f)^2}} \mathds{1}_{A_n}\bigg|\bigg|_{\infty}  \,\bigg|\bigg|\frac{p_0}{p}\bigg|\bigg|_{\infty}.
\end{align*}
Based on the Assumption \ref{ass:f0}, $f_0$ is a $\beta$-smooth function supported on $[0,1]$ and hence there exists some constant $B'_0>0$ such that $\|f_0\|_{\infty}\le B'_0$. For $y\in A_n$, we have
\begin{align*}
\frac{e^{-\frac{1}{2} \{y-f(x)\}^2}}{e^{-\frac{1}{2} \{y-f_0(x)\}^2}} &= e^{\{f(x)-f_0(x)\}\{y-f_0(x)\} - \{f(x)-f_0(x)\}^2/2} \\
&\ge e^{-\|f-f_0\|_{\infty}(|y|+\|f_0\|_{\infty}) - \{f(x)-f_0(x)\}^2/2}  \\
&\ge  e^{-\epsilon_n(\gamma'/\epsilon_n+B'_0)-\epsilon_n^2/2}\ge e^{-2\gamma'}.
\end{align*} 
Thus $\|e^{- (y-f_0)^2/2}/e^{-(y-f)^2/2} \mathds{1}_{A_n}\|_{\infty} \le e^{2\gamma'}$. Based on Lemma \ref{lem:sieve}, for any $x\in [0,1]$ and $p\in \m{P}_n$, we have $\log \| p_0/p\|_{\infty}\lesssim \log(1/\epsilon_n)$. Therefore, we have shown 
\begin{align}\label{I_1}
I_1=\int_{A_n} g_{f_0,p_0} \{\log (g_{f_0,p_0}/g_{f,p})\}^2 \le 2\epsilon_n^2 \log ^2(1/\epsilon_n).
\end{align}
Next we estimate the term $I_2$ in Equation \eqref{eq:denom}. For all $y\in A_n^c$ and for any fixed $x\in[0,1]$, we choose $\gamma'>1$ such that $|y-f_0(x)|\ge |y|-\|f_0\|_{\infty}> \gamma'/\epsilon_n-B'_0 \ge 1/\epsilon_n$. By Fubini's theorem,
\begin{align*}
&\int_{|y|>1/\epsilon_n} g_{f_0,p_0} \bigg(\log\frac{g_{f_0,p_0}}{g_{f,p}}\bigg)^2\\
&\le\frac{1}{2\pi\delta_n} \int_0^1 \int_{|y-f_0(x)|>1/\epsilon_n} e^{-\frac{1}{2}\{y-f_0(x)\}^2} e^{-\frac{1}{2\delta_n^2}(w-x)^2} \\
&~~~~~~~~~~~~~~~~~~~~~~~~~~~~~~~~~~~\cdot\left(\log \frac{\int e^{-\frac{1}{2}(y-f_0)^2} e^{-\frac{1}{2\delta_n^2}(w-x)^2}\,p_0(x)\,dx}{\int e^{-\frac{1}{2}\{y-f(x)\}^2} e^{-\frac{1}{2\delta_n^2}(w-x)^2}\,p(x)\,dx}\right)^2 dy\,dw\, p_0(x)\,dx\\
&\le \frac{1}{\sqrt{2\pi}} \int_0^1 \int_{|y-f_0(x)|>1/\epsilon_n} e^{-\frac{1}{2}\{y-f_0(x)\}^2} \bigg(\log \bigg|\bigg|\frac{e^{-(y-f_0)^2/2}}{e^{-(y-f)^2/2}}\bigg|\bigg|_{\infty} + \log \bigg|\bigg|\frac{p_0}{p}\bigg|\bigg|_{\infty}\bigg)^2 dy\, p_0(x)\, dx.
\end{align*}
Let $z=y-f_0(x)$, we can show that for any $x\in[0,1]$, $e^{-\{y-f_0(x)\}^2/2+\{y-f(x)\}^2/2} \le e^{\epsilon_n |z| + \epsilon_n^2/2}$. Then
\begin{align*}
&\int_{A_n^c} g_{f_0,p_0} \bigg(\log\frac{g_{f_0,p_0}}{g_{f,p}}\bigg)^2\\
 &\le 4(2\pi)^{-1/2} \int_0^1 \bigg( \int_{|z|\ge 1/\epsilon_n} e^{-\frac{1}{2}z^2} (\epsilon_n z)^2\, dz+  \int_{|z|\ge 1/\epsilon_n} e^{-\frac{1}{2}z^2} \log^2(1/\epsilon_n)dz\bigg)\, p_0(x)\, dx\\
&\le 4(2\pi)^{-1/2} E_0 \bigg\{ \epsilon_n^2 \int_{t>1/\epsilon_n^2} e^{-t/2}t^{1/2}\,dt +\log^2(1/\epsilon_n) P(|Z|\ge 1/\epsilon_n) \bigg\} \\
& \le 4 (2\pi)^{-1/2} E_0\bigg\{ \epsilon_n^2 \int_{t>1/\epsilon_n^2} e^{-t/4}\, dt + \log^2(1/\epsilon_n) e^{-\epsilon_n^{-2}/8}\bigg\} \\
&\lesssim  e^{-\epsilon_n^{-2}/8 + \log\log(1/\epsilon_n)} < \epsilon_n^2,
\end{align*}
where $Z\sim \text{N}(0,1)$ and $E_0(\cdot)$ denotes taking expectation with respect to the measure associated with the density $p_0$. The third line in the preceding uses the change of variables letting $t=z^2$. 

Combining the above result for $I_2$ and the result in Equation \eqref{I_1} for $I_1$, we have shown $\int g_{f_0,p_0}(\log g_{f_0,p_0}/g_{f,p})^2\lesssim \epsilon_n^2$. And further we have
\begin{align*}
\bigg\{\int g_{f_0,p_0}\log\frac{g_{f_0,p_0}}{g_{f,p}}\lesssim \epsilon_n^2,\ &\int g_{f_0,p_0}\bigg(\log\frac{g_{f_0,p_0}}{g_{f,p}}\bigg)^2\lesssim \epsilon_n^2 \bigg\} \\
&\supset \{\|f-f_0\|_{\infty} \le \epsilon_n, \ KL(p_0,p)\le\epsilon_n^2 \},
\end{align*}
which yields the conclusion in Lemma \ref{lem:gkl}.

\section{Posterior Computation: A Gibbs Sampler}\label{sec:gibbcom}
In the following, we develop a Gibbs sampler to generate a Markov chain which will eventually  converge to the posterior distribution. We focus on the Gaussian process associated with a squared exponential kernel as an  illustration (in practice the algorithm can be applied to other kernels as long as they are symmetric). The squared exponential kernel is denoted by $c(x,x') = \exp\{-(x-x')^2/\lambda\}$ associated with a bandwidth parameter $\lambda$. Theorem \ref{thm:appgp} enforces the prior distributions $w_j\sim \mbox{N}(0, 2/\lambda)$, $s_j \sim \mbox{Unif}\,[0,2\pi]$ and $a_j\sim \mbox{N}(0,1)$ i.i.d. for $j=1,\ldots,N$. To ensure the conditional conjugacy, we place a gamma distribution $\text{Ga}(a_0,b_0)$ on the bandwidth $\lambda$ with a shape parameter $a_0$ and a scale parameter $b_0$. We place a Dirichlet process mixture of normals prior defined in Equation \eqref{eq:DPMM} over the covariate density, given more precisely by 
\begin{eqnarray}\label{eq:DPMMb}
X_i \sim \sum_{h=1}^{\infty} \pi_h \mbox{N}(\mu_h, \tau_h^{-1}), \quad (\mu_h, \tau_h) \sim  \mbox{N}(\mu_h;\mu_0, \kappa_0\tau_h^{-1})\mbox{Ga}(\tau_h; a_{\tau},b_{\tau}),
\end{eqnarray}
for $i=1,\ldots,n$. The prior on  $\pi_h$ is expressed as $\pi_h = \nu_h \prod_{l < h} (1-\nu_l)$ where $\nu_l \sim \mbox{Beta}(1, \alpha)$. Here we let $\alpha =1$. Denote the cluster label of $X_i$ by $S_i \in \{1,\dots, K\}$ indicating that $X_i$ is associated with $S_i {\rm th}$ component in the Dirichlet process Gaussian mixture prior for $i=1,\ldots,n$.  Then Equation \eqref{eq:DPMMb} can be also written as 
\begin{eqnarray*}
X_i \mid S_i, \mu, \tau \sim  \mbox{N}(\mu_{S_i}, \tau_{S_i}^{-1}), \quad (\mu_{S_i}, \tau_{S_i}) \sim  \mbox{N}(\mu_{S_i};\mu_0, \kappa_0\tau_{S_i}^{-1})\,\mbox{Ga}(\tau_{S_i}; a_{\tau},b_{\tau}),\\
i=1,\ldots,n.
\end{eqnarray*}
In both simulation studies and the real application, we set the hyperparameters $\mu_0=0, \kappa_0 = 1,a_{\tau} = 1, b_{\tau} =1$, and we choose $a_0 = 5, b_0 =1$ for the hyperprior $\text{Ga}(a_0,b_0)$. We remark that these hyperparameter choices are based on our preliminary numerical experiments. In addition, recall that we assume $\sigma = 0.2$ in simulation studies and we treat $\sigma^2$ as an unknown parameter endowed with an objective prior in real application. 

As below we provide a complete updating scheme of the Gibbs sampler. We use bold symbols to distinguish the vectors $\mathbf{a}, \mathbf{w}, \mathbf{s}, \boldsymbol{\mu},\boldsymbol{\tau},\boldsymbol{\pi}, \mathbf{S},\mathbf{X}, \mathbf{Y},\mathbf{W}$ accordingly.  Then the joint posterior distribution of $\{\mathbf{a}, \mathbf{w}, \mathbf{s}, \lambda, \mathbf{X}\}$ given observations $\{\mathbf{Y},\mathbf{W}\}$ can be factorized as
\begin{eqnarray*}
[\,\mathbf{a}, \mathbf{w}, \mathbf{s}, \lambda, \mathbf{X} \mid  \mathbf{Y},\mathbf{W}\,] \propto [\,\mathbf{Y} \mid \mathbf{X},\mathbf{a},\mathbf{w},\mathbf{s}, \lambda\,] \times [\mathbf{W} \mid \mathbf{X}]\times [\mathbf{w} \mid  \lambda]\times [\lambda]\times [\mathbf{a}] \times [\mathbf{s}] \times [\mathbf{X}].
\end{eqnarray*}
The updating scheme runs as follows: 
\begin{enumerate}
\item Update $[\,\mathbf{w} \mid -]$ in a block by sampling $[w_j \mid -] \propto  [\mathbf{Y}\mid \mathbf{X}, \mathbf{a}, \mathbf{w}, \mathbf{s}, \lambda]\, \mbox{N}(w_j; 0, 2/\lambda)$ independently using Metropolis-Hasting algorithm for $j=1,\ldots,N$.
\item Update $[\,\mathbf{s} \mid -]$ in a block by sampling $[s_j \mid -] \propto [\mathbf{Y}\mid \mathbf{X}, \mathbf{a}, \mathbf{w}, \mathbf{s}, \lambda]\, \mbox{Unif}\,[0,2\pi]$ independently using Metropolis-Hasting algorithm for $j=1,\ldots,N$.
\item Update $[\,\mathbf{a} \mid -] $ from a multivariate normal distribution $ \mbox{N}(\widetilde{\boldsymbol{\mu}},\widetilde{\mathbf{\Sigma}})$, with the mean vector $\widetilde{\boldsymbol{\mu}}= \widetilde{\Sigma}\,\mathbf{\Phi}^{\T}\,\mathbf{Y}/\sigma^2$, and the covariance matrix $\widetilde{\Sigma}=(\mathbf{\Phi}^{\T}\mathbf{\Phi}/\sigma^2 +\mathbf{I}_N)^{-1}$, where $\mathbf{\Phi}$ is a $n\times N$ Fourier basis matrix with $(i,j)$th element $\mathbf{\Phi}_{ij} = (2/N)^{1/2}\,\cos(w_j x_i +s_j)$ for $i=1,\ldots, n$, $j=1,\ldots,N$. And $\mathbf{I}_N$ denotes a $N\times N$ identity matrix.
\item Update the parameters $[\,\mathbf{S},\boldsymbol{\mu},\boldsymbol{\tau},\boldsymbol{\pi} \mid -]$ associated with the Dirichlet process Gaussian mixture prior as in \cite{ishwaran2001gibbs} with the number of mixture components truncated at $20$.
\item  Update $[\, \mathbf{X} \mid -]$ in a block by sampling 
$$[X_i| S_i,X_{-i},-]\propto  \mbox{N}(Y_i; \mathbf{\Phi}^\T_{i} \mathbf{a}, \sigma^2)\, \mbox{N}(W_i; X_i, \delta^2)\,  \mbox{N}(X_i; \mu_{S_i}, \tau_{S_i})$$ using Metropolis-Hasting algorithm for $i=1,\ldots,n$. Here $\mathbf{\Phi}^\T_{i}$ denotes the $i$th row of the matrix $\mathbf{\Phi}$ defined in Step 3. 
\item Update $[\,\lambda \mid -]$ from a gamma distribution $\mbox{Ga}(\widehat{a}, \widehat{b})$ with $\widehat{a}=a_0$ and $\widehat{b}=b_0/(1+b_0\sum_{j=1}^n w_j^2/4)$.
\item Update $[\,\sigma^2 \mid -]$ from a inverse-gamma distribution $\mbox{IG}(a_{\sigma}, b_{\sigma})$ with $a_{\sigma} = n/2$ and $b_{\sigma} = (\mathbf{Y} - \mathbf{\Phi}\ind_N)^{\T} (\mathbf{Y} - \mathbf{\Phi}\ind_N)/2$, where $\ind_N$ denotes a  $n\times 1$ vector of ones. (This step will be implemented only in the real example of Section \ref{sec:real}.)
\end{enumerate}
In particular, in Metropolis-Hasting algorithm used for updating $\{w_j\}$ in Step 1, we consider a random walk proposal $ w_j^{\mathrm{prop}} \sim \mbox{N}(w_j^{\mathrm{cur}} , 1/4)$ for $j=1,\ldots, N$, where $w_j^{\mathrm{cur}}$ denotes the current state and the proposal variance is tuned to obtain average pointwise acceptance rate around $0.7$. In Metropolis-Hasting algorithm used for updating $\{s_i\}$ in Step 2, we consider the independence proposal $s_i^{\mathrm{prop}} \sim \mbox{Unif}\,[0, 2\pi]$ for $i=1,\ldots,n$. We note that the averaged pointwise acceptance rate for $s_i$ is around $0.6$.  Finally, to update $\{x_i\}$ in Step 5, we use an adaptive proposal $x_i^{\mathrm{prop}} \sim \mbox{N}(W_i/\delta^2 + \mu_{S_i}\tau_{S_i},1/(1/\delta^2+\tau_{S_i}))$ for $i=1,\ldots,n$ with the averaged acceptance rate around $0.8$.  \\

\noindent{\bf Constructing the spontaneous credible bands.} We provide one example of constructing the spontaneous credible bands (CB) with $\gamma=0.95$ for out-of-sample prediction of some model $f(x,\theta)$ evaluated at a test data set $x_{test}$ of size $n_t$, based on $L$ number of posterior samples $\{\theta^{(l)}, l=1,\ldots L\}$ of parameter $\theta$ associated with the model $f$. Denote by $f^{(l)}(x) = f(x;\theta^{(l)})$ for $l=1,\ldots,L$ and let $\hat{f}(x) = (1/L)\sum_{l=1}^n f(x;\theta^{(l)})$ denote the posterior estimate of the function.  Then, for each $l=1\ldots,L$, we first calculate the maximum distance between the functions $\hat{f}(x)$ and $f^{(l)}(x)$ over the test data points, defined as $d_l = \max_{i=1,\ldots,n_t}|f^{(l)}(x_{{test},i}) - \hat{f}(x_{{text},i})|$. To find the simultaneous CB, we find the $95\%$ quantile of $\{d_l\}$ denoted by $d_{95\%}$ and take $d_{95\%}$ as the half range of the simultaneous CB. Then we define the spontaneous 95\% credible band as $[\hat{f}(x_{text})-d_{95\%}, \hat{f}(x_{text}) + d_{95\%}]$.

\section{Additional Numerical Results} \label{app:addres}

\begin{figure}[h!]
\centering


\begin{multicols}{3}
    \includegraphics[scale = 0.045]{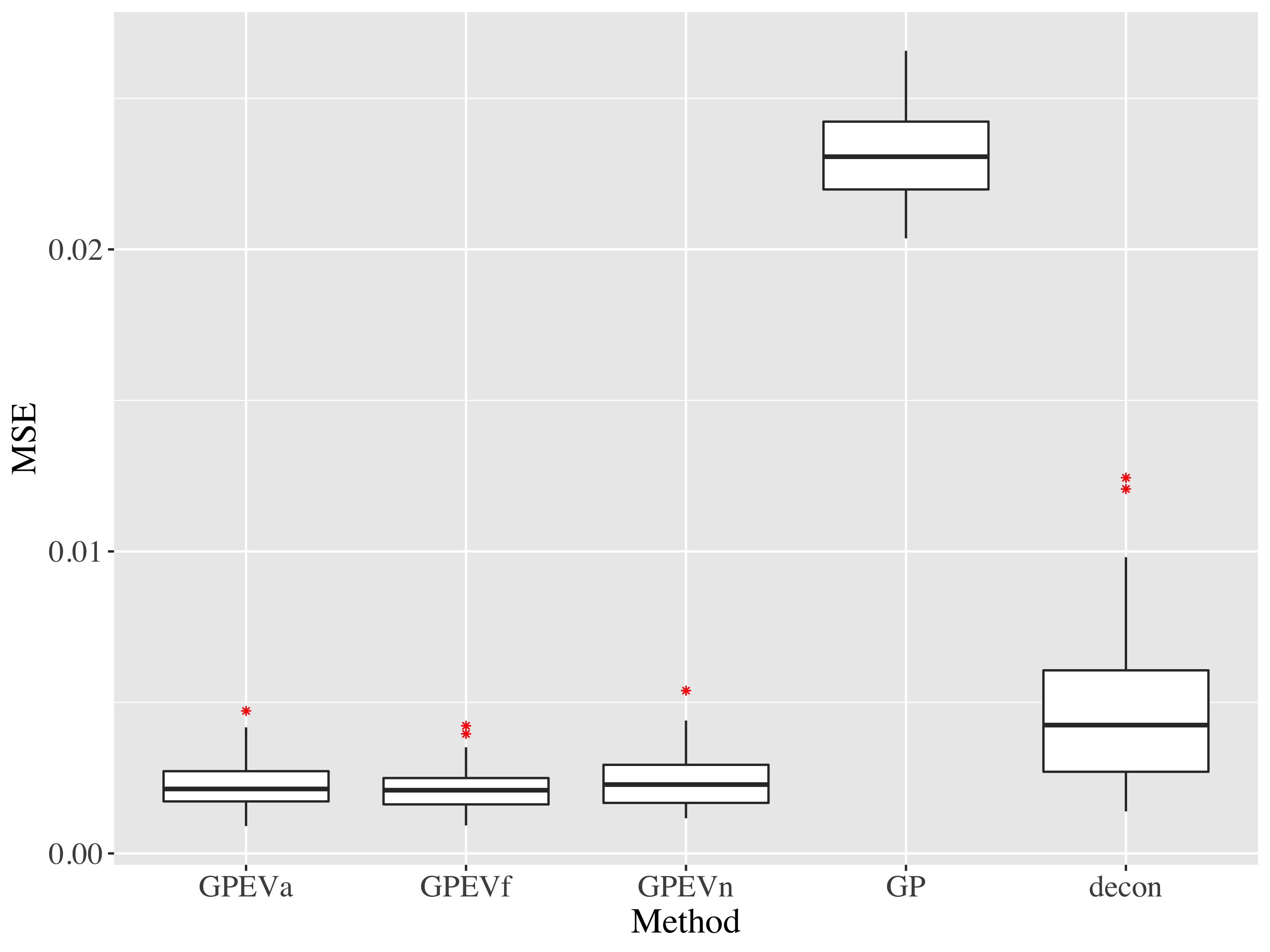}
    \includegraphics[scale = 0.045]{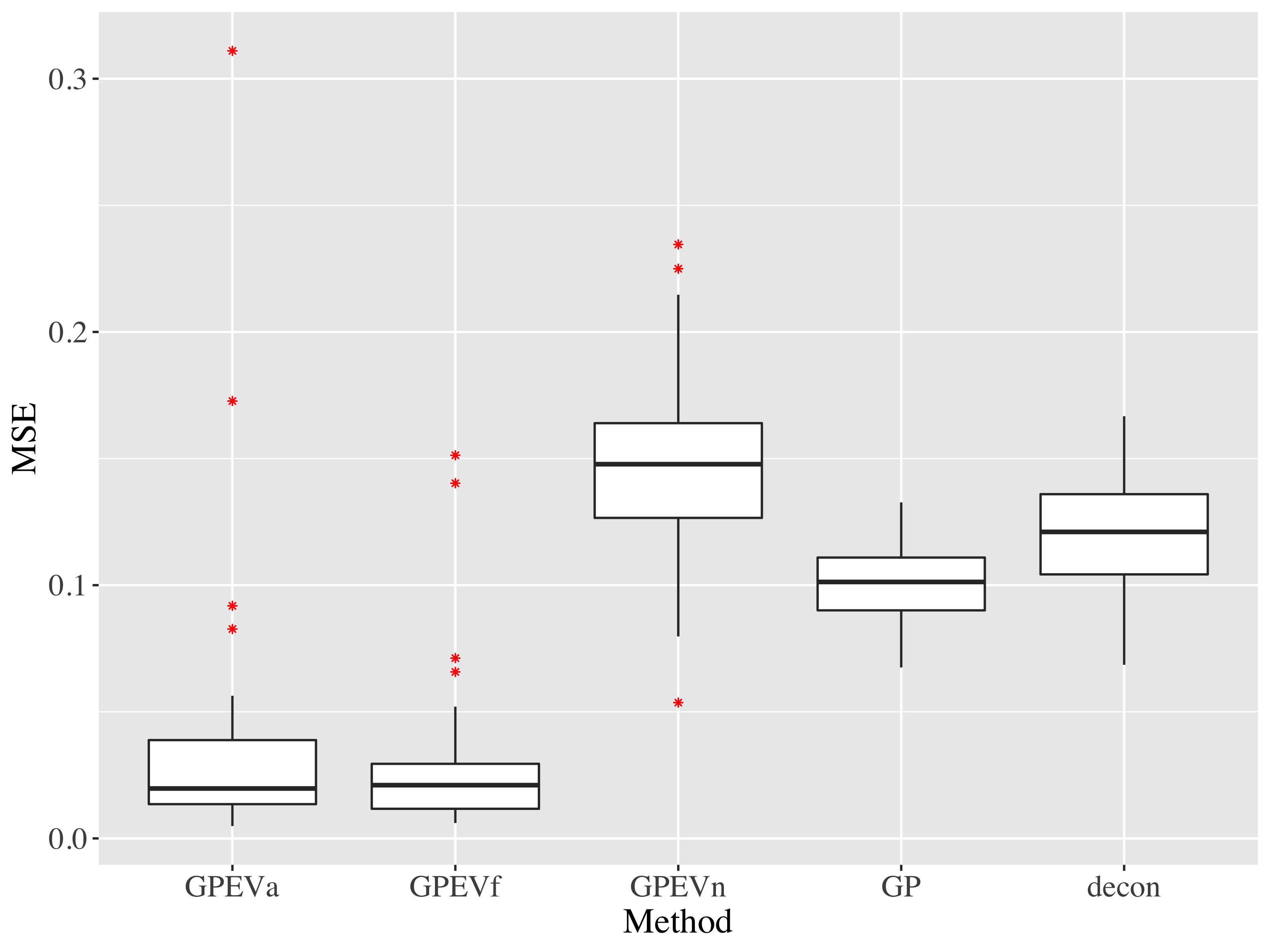}
    \includegraphics[scale = 0.045]{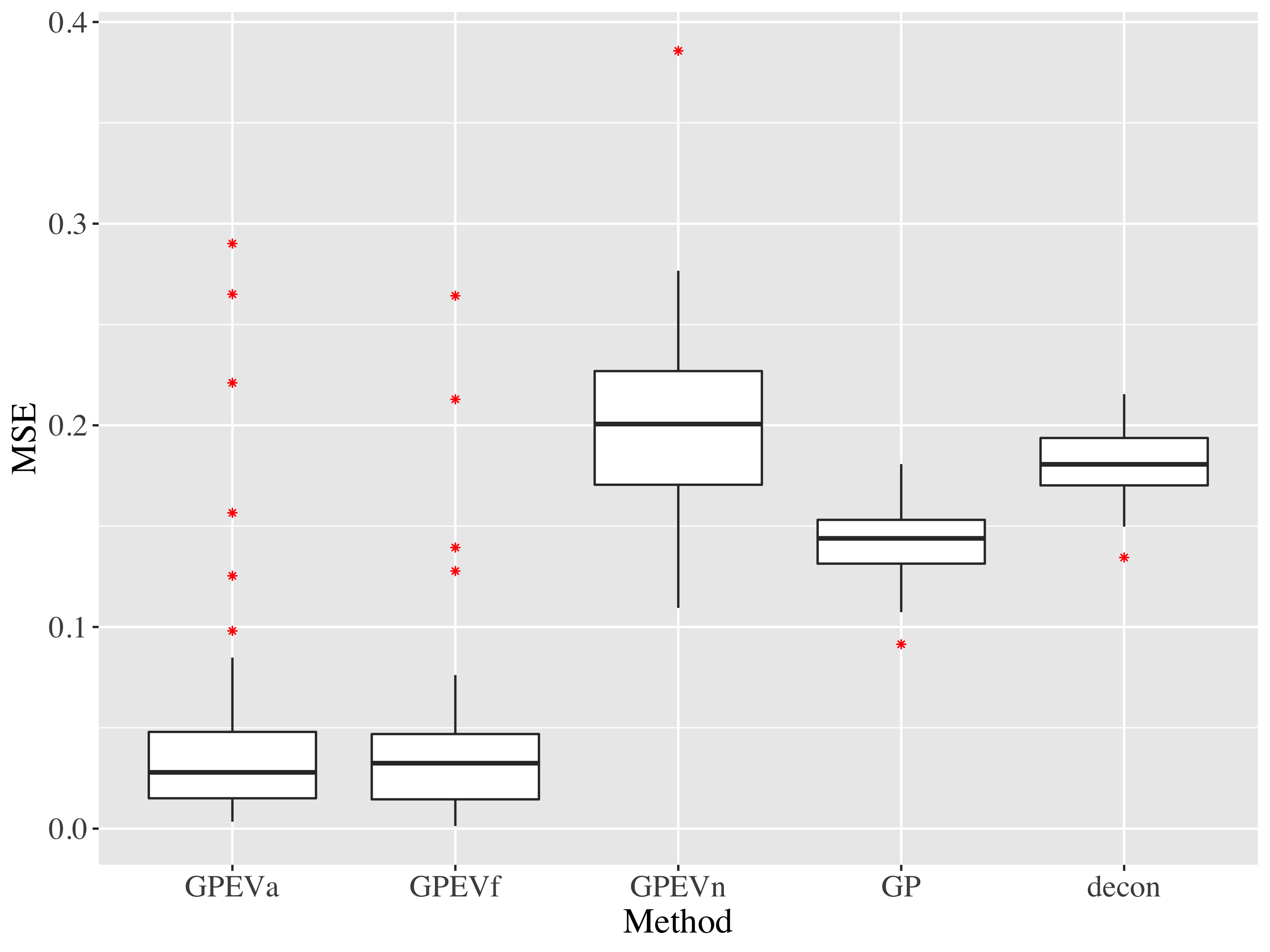}
\end{multicols}
\caption{Boxplots of mean squared errors for compared methods in Section \ref{sec:sims} over 50 replicated data sets of size $n = 250$ with $\delta^2 = 0.01$ (left panel), $\delta^2 = 0.6$ (middle panel) and  $\delta^2 = 1$ (right panel). In each panel the compared methods from left to right are \textsc{gpev}$_a$, \textsc{gpev}$_f$, \textsc{gpev}$_n$, \textsc{gp} and \text{decon}.} \small
\label{fig:boxplot_n250}

\end{figure}

\begin{figure}[h!]
\centering
\begin{multicols}{3}
    \includegraphics[scale = 0.045]{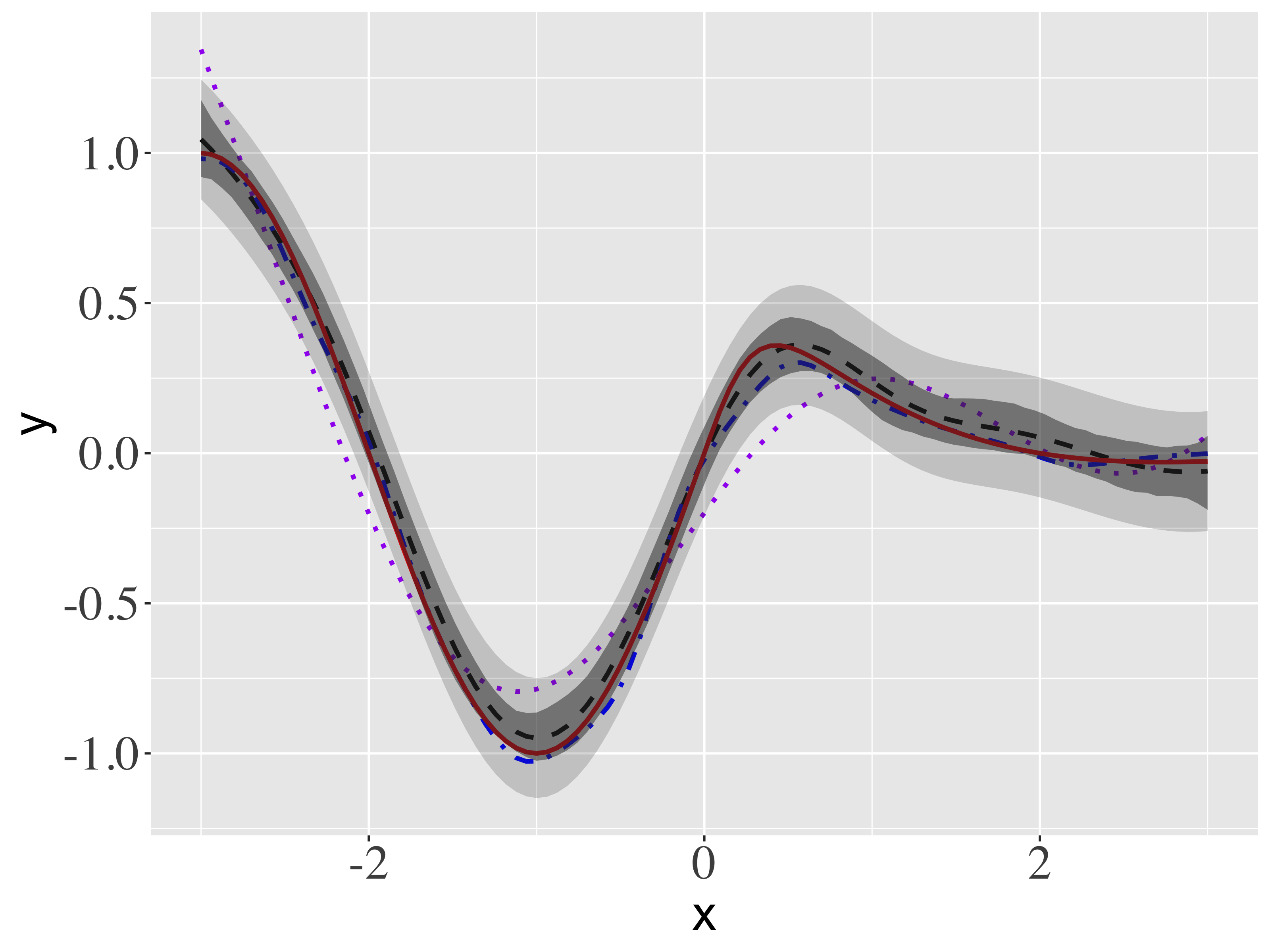}
    \includegraphics[scale = 0.045]{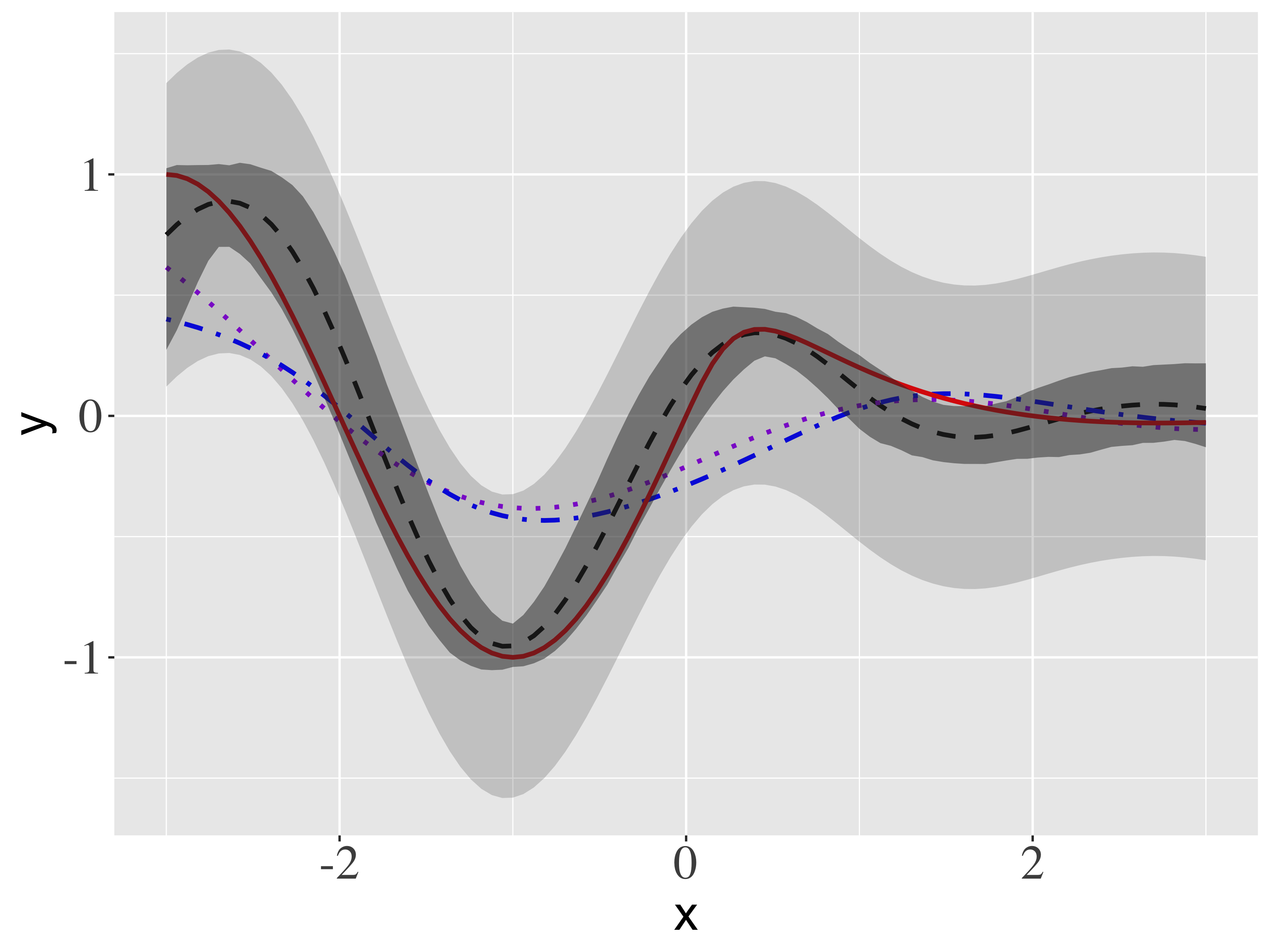}
    \includegraphics[scale = 0.045]{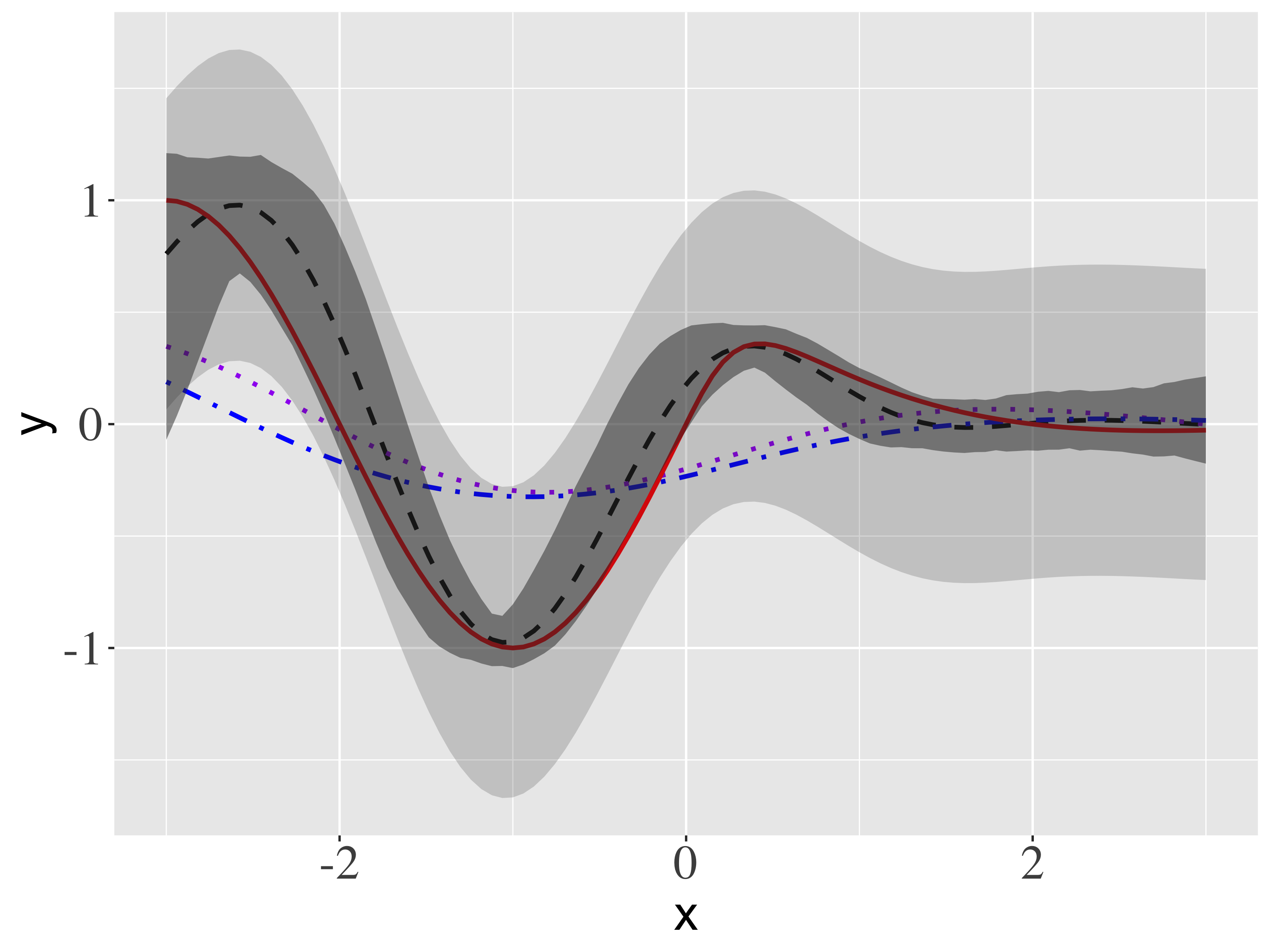}
\end{multicols}
\caption{Out-of-sample predictions of $f(x)$ for $\delta^2 = 0.01$ (left panel), $\delta^2 = 0.6$ (middle panel) and  $\delta^2 = 1$ (right panel) with sample size $n = 250$. The red solid line stands for the true function, the black dashed line stands for the predictive curve based on \textsc{gpev}$_a$, the blue dot-dashed line is based on \text{decon} and the purple dotted dashed line is based on \textsc{gp}. The darker and the lighter shades are the pointwise and simultaneous $95\%$ credible intervals of \textsc{gpev}$_a$, respectively.} \small

\label{fig:fit_n250}
\end{figure}

\begin{figure}[h!]
\centering 
\begin{multicols}{3}
   \includegraphics[scale = 0.045]{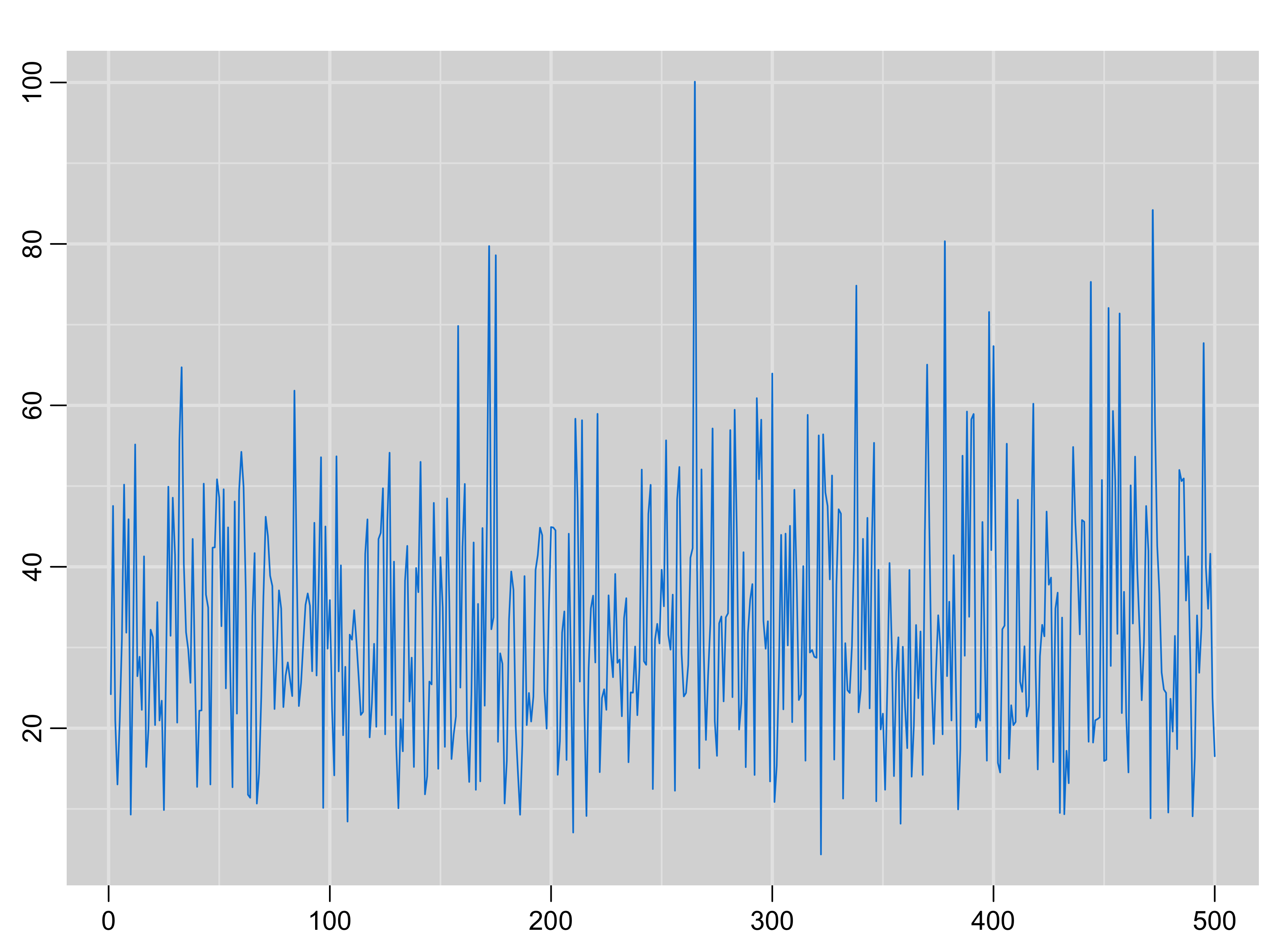}
    \includegraphics[scale = 0.045]{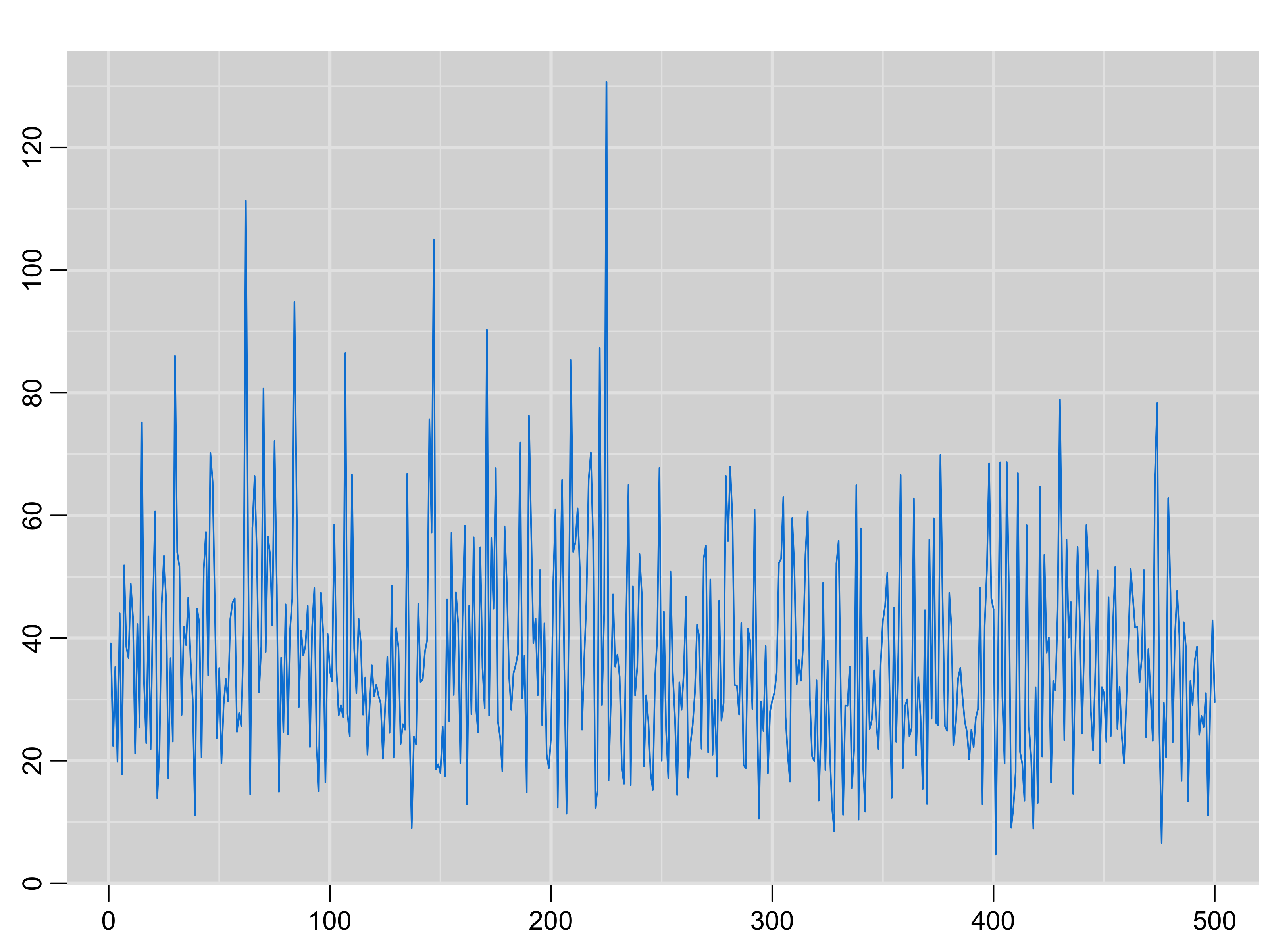}
    \includegraphics[scale = 0.045]{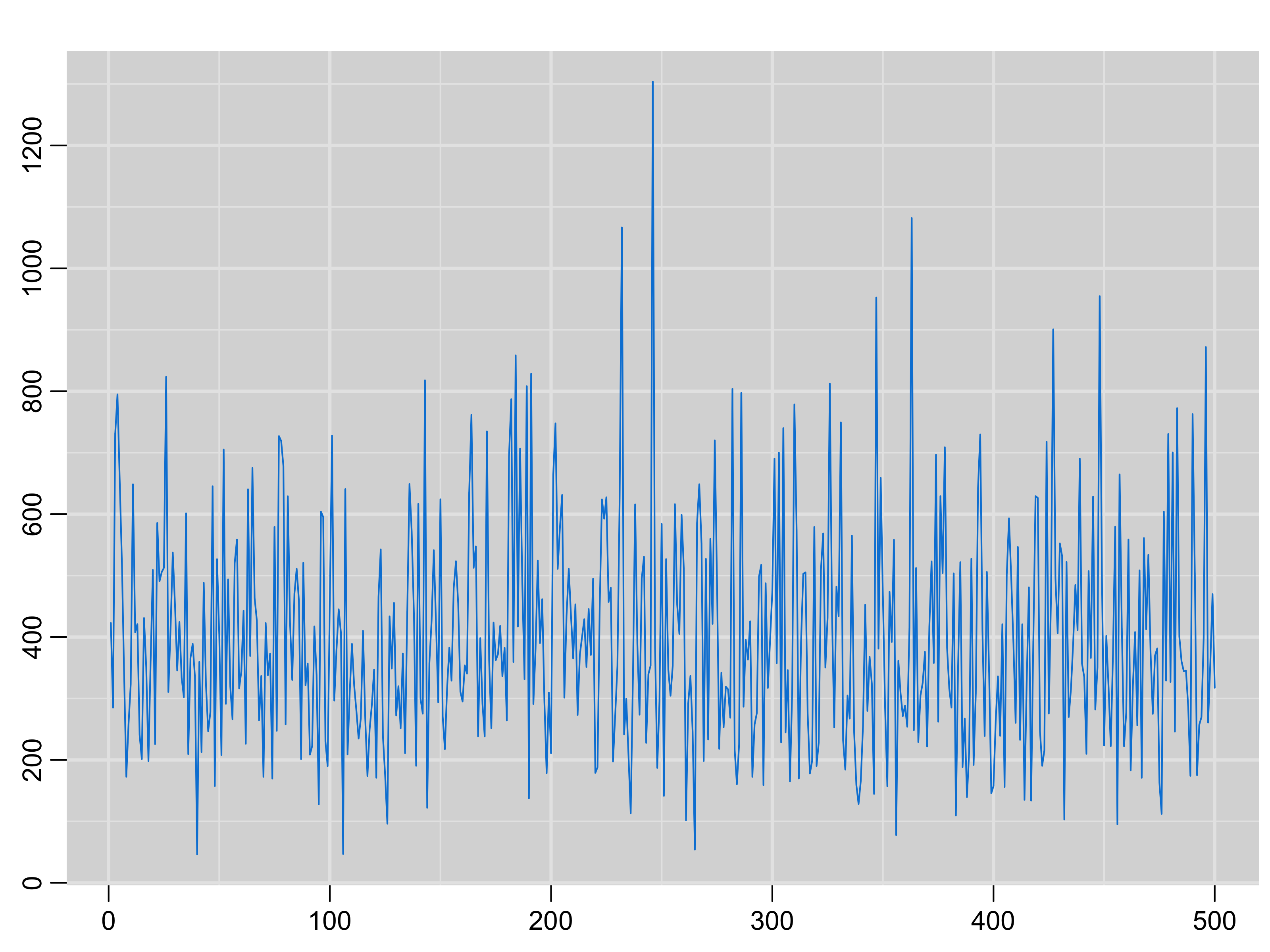}
\end{multicols}

\begin{multicols}{3}
    \includegraphics[scale = 0.045]{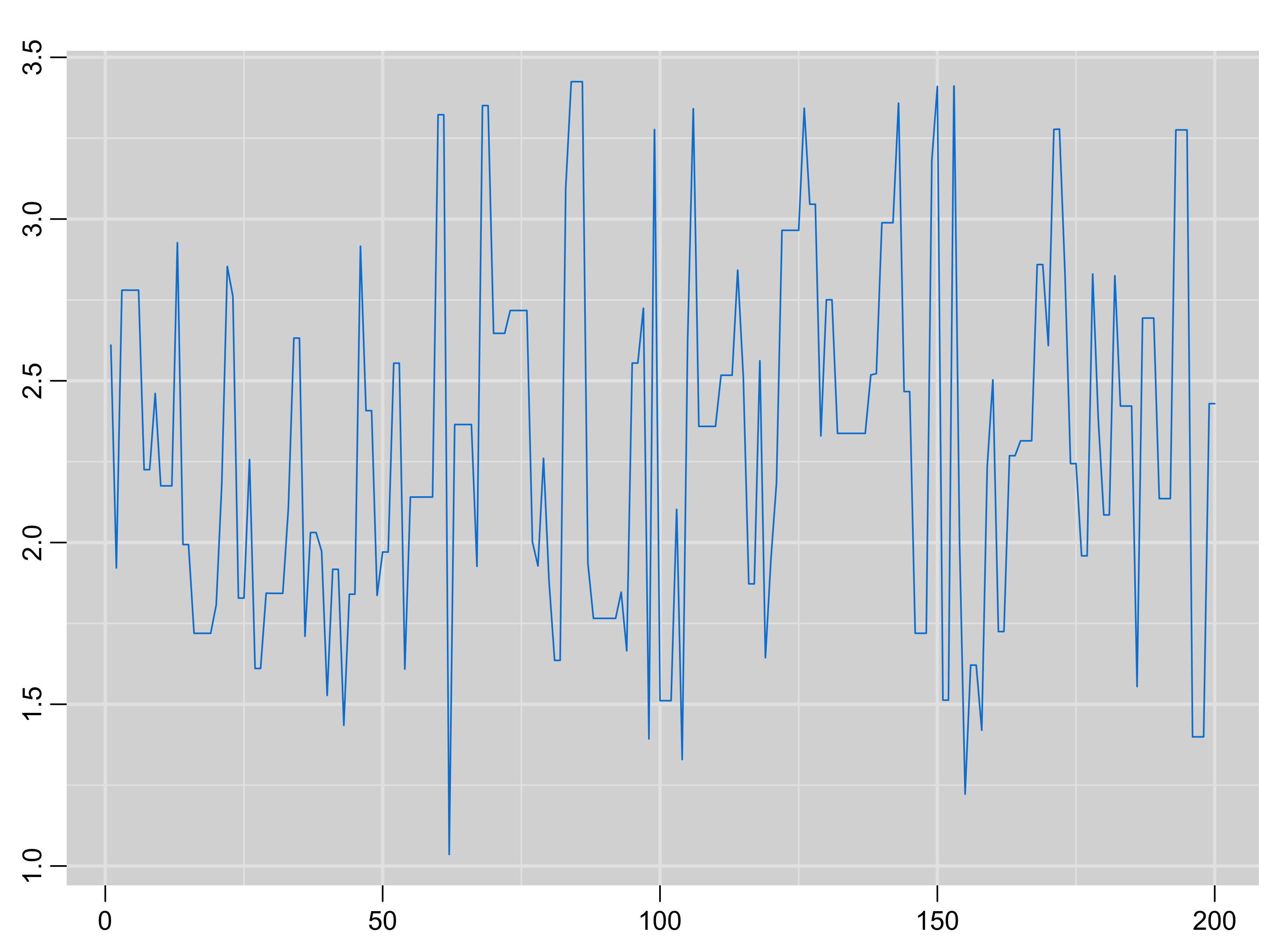}
    \includegraphics[scale = 0.045]{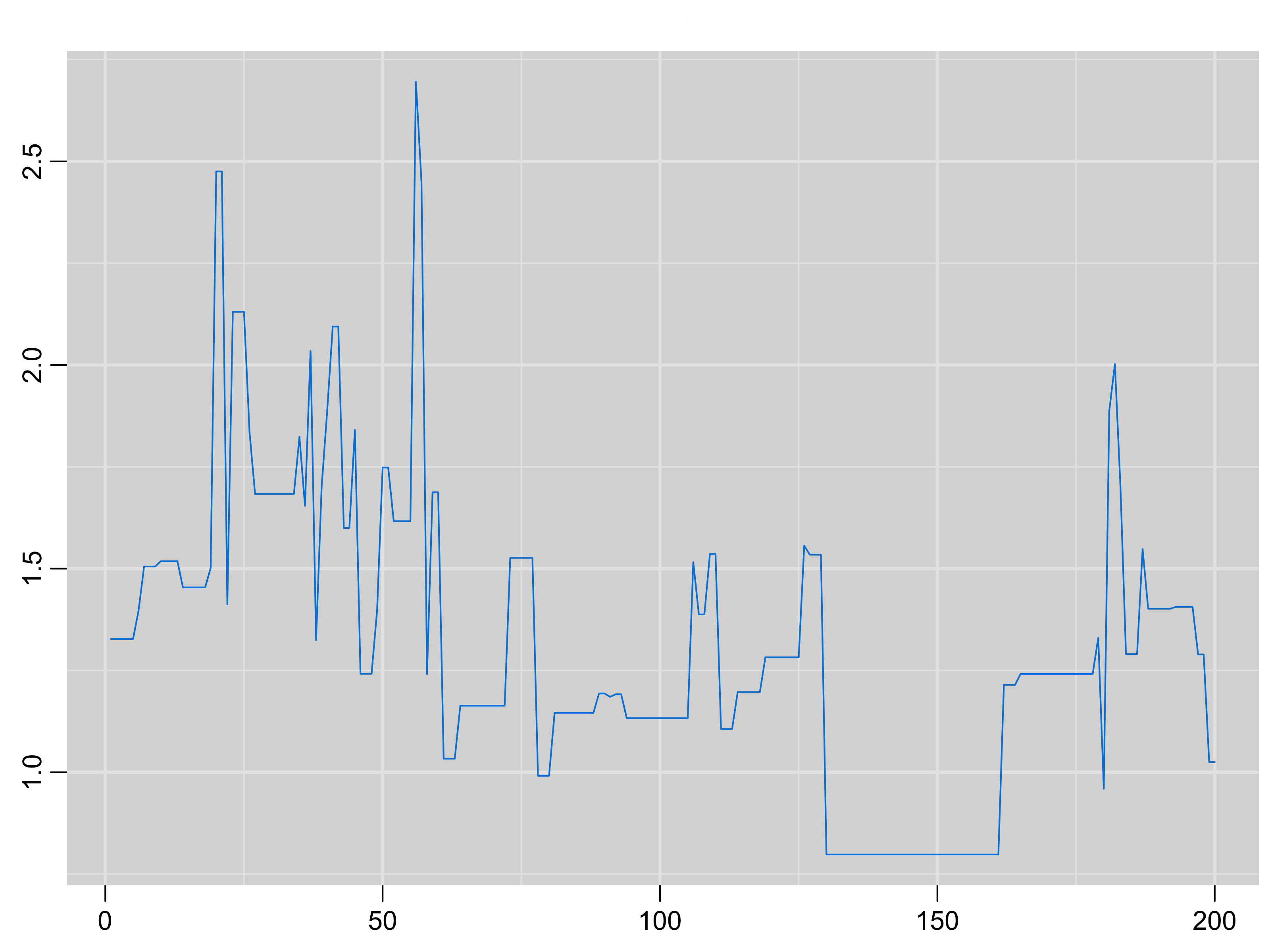}
    \includegraphics[scale = 0.045]{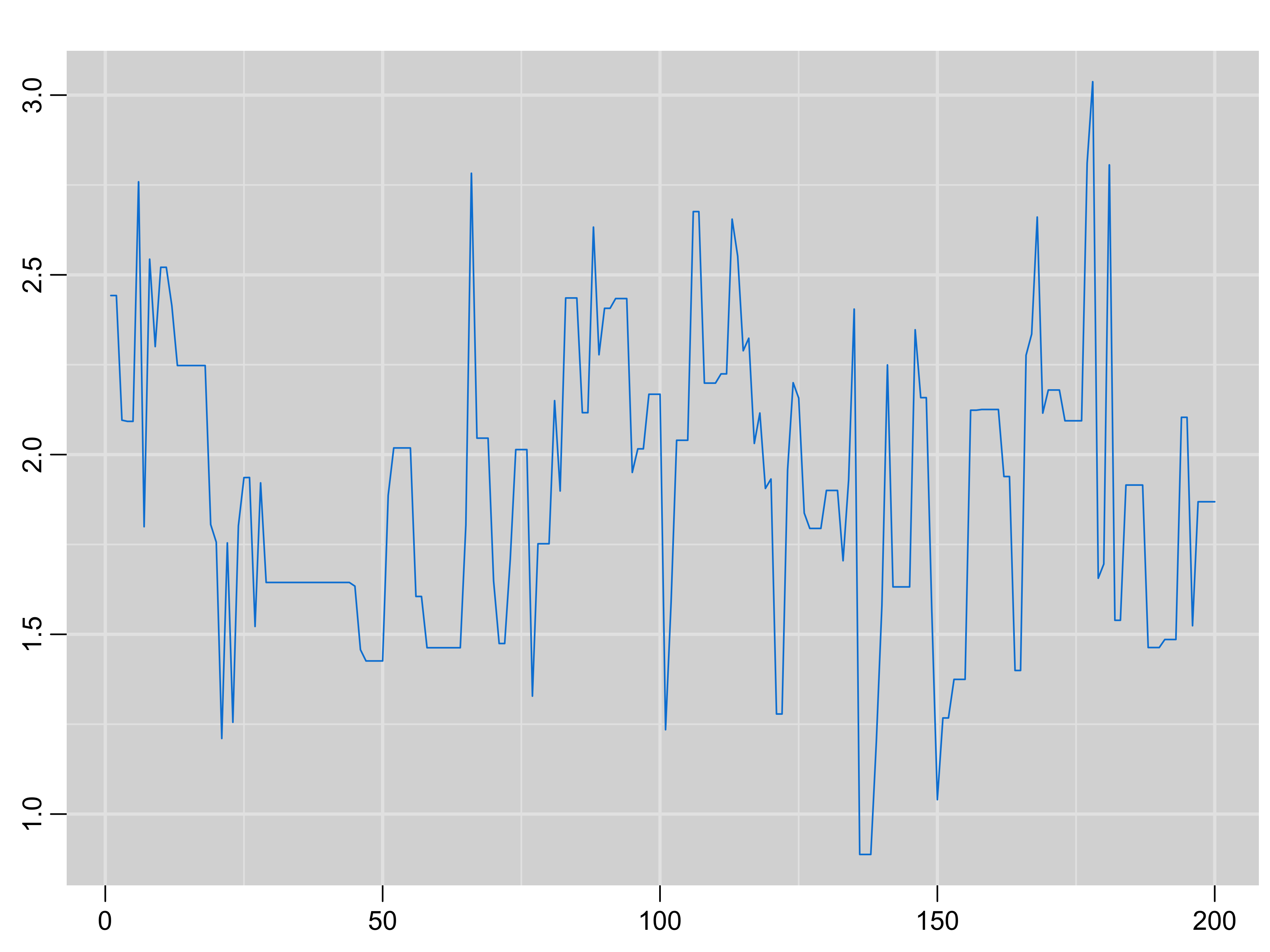}
\end{multicols}

\caption{Trace plots of posterior samples of $\lambda$ from \textsc{gpev}$_a$ (first row) and \textsc{gpev}$_f$ (second row) with sample size $n=100$. In each row, the values of $\delta^2$ are $0.01$ (left panel), $0.6$ (middle panel) and $1$ (right panel).}\small
\label{fig:mixing_lambda}

\end{figure}

\begin{figure}[h!]
\centering
\begin{multicols}{3}
    \includegraphics[scale = 0.045]{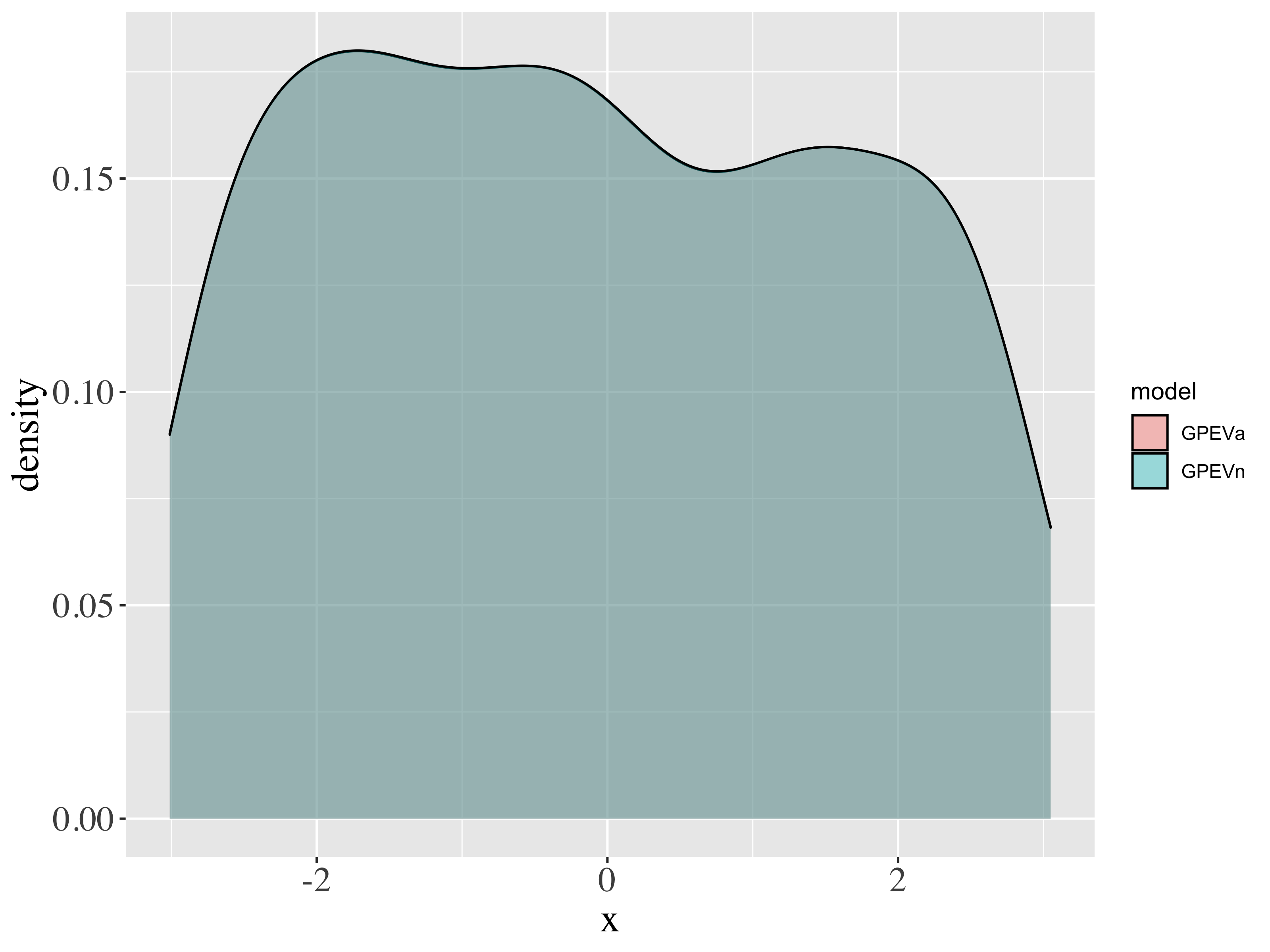}
    \includegraphics[scale = 0.045]{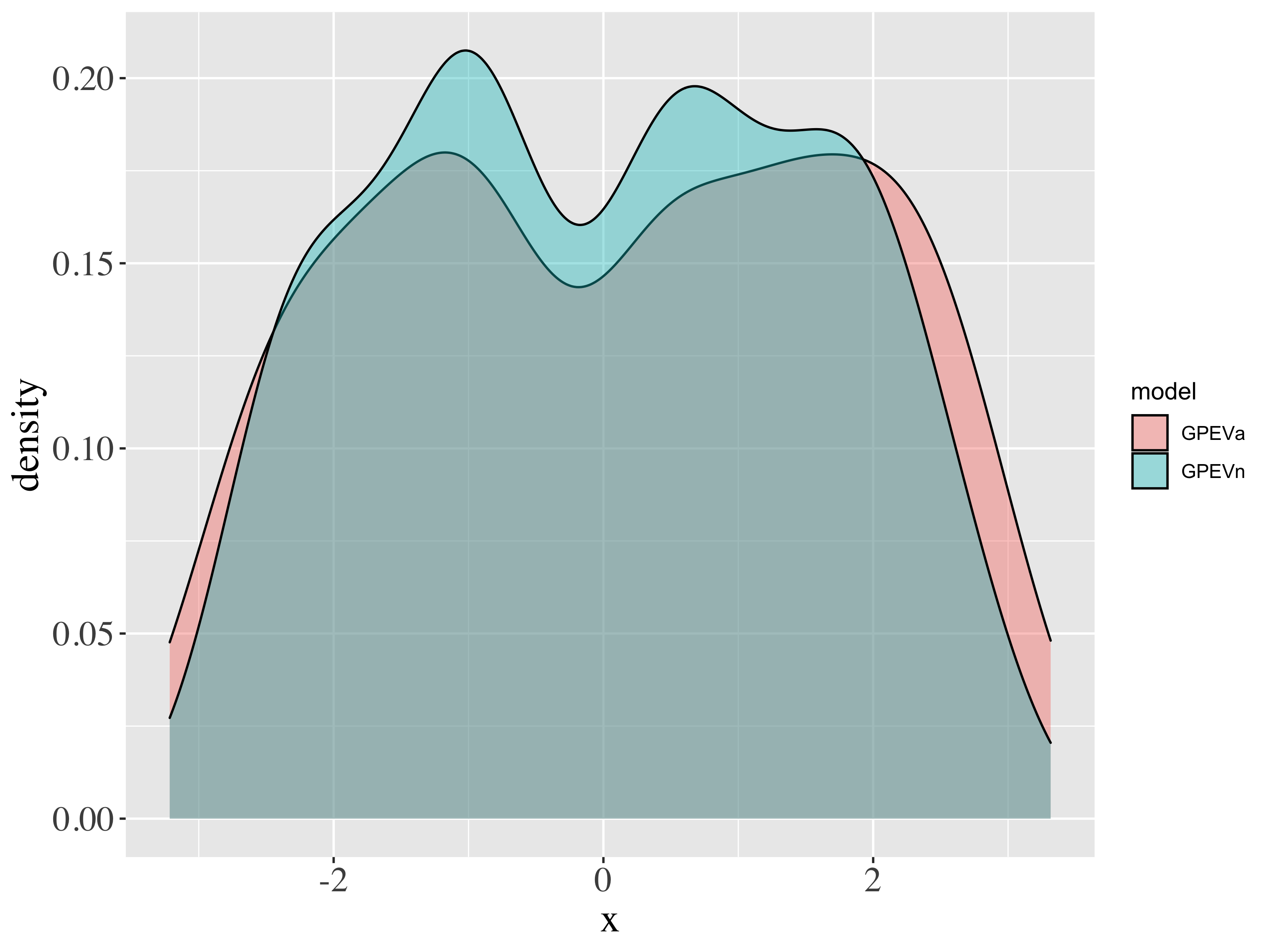}
    \includegraphics[scale = 0.045]{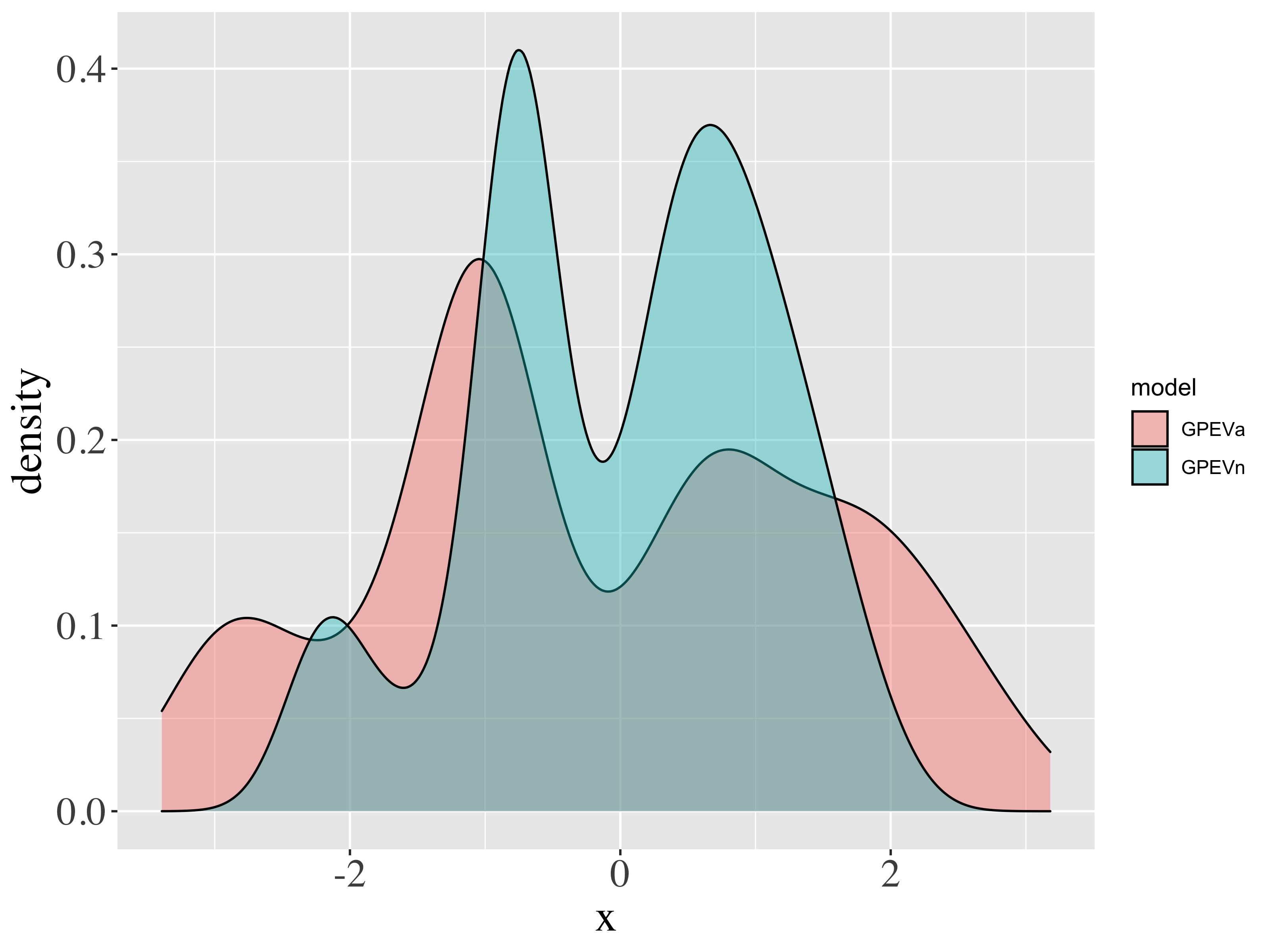}
\end{multicols}
\caption{Marginal posterior density plots of the covariate based on $\textsc{gpev}_a$ (red) and $\textsc{gpev}_n$ (green) with $n=500$. The values of $\delta^2$ are $0.001$ (left panel), $0.1$ (middle panel), $0.5$ (right panel). 
} \small
\label{fig:p0}
\end{figure}

\begin{figure}[h!]
\centering 
\begin{multicols}{3}
   \includegraphics[scale = 0.045]{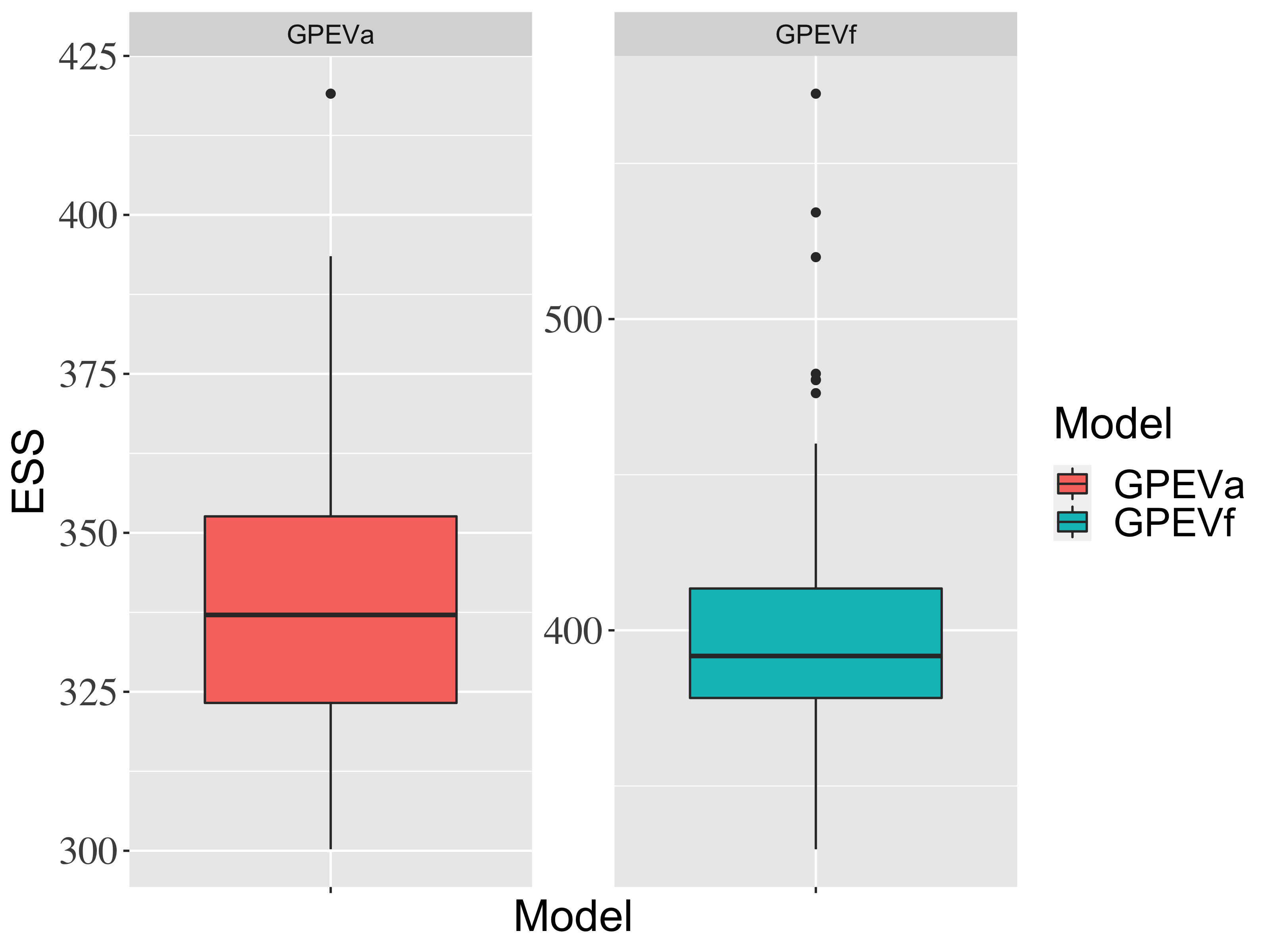}
    \includegraphics[scale = 0.045]{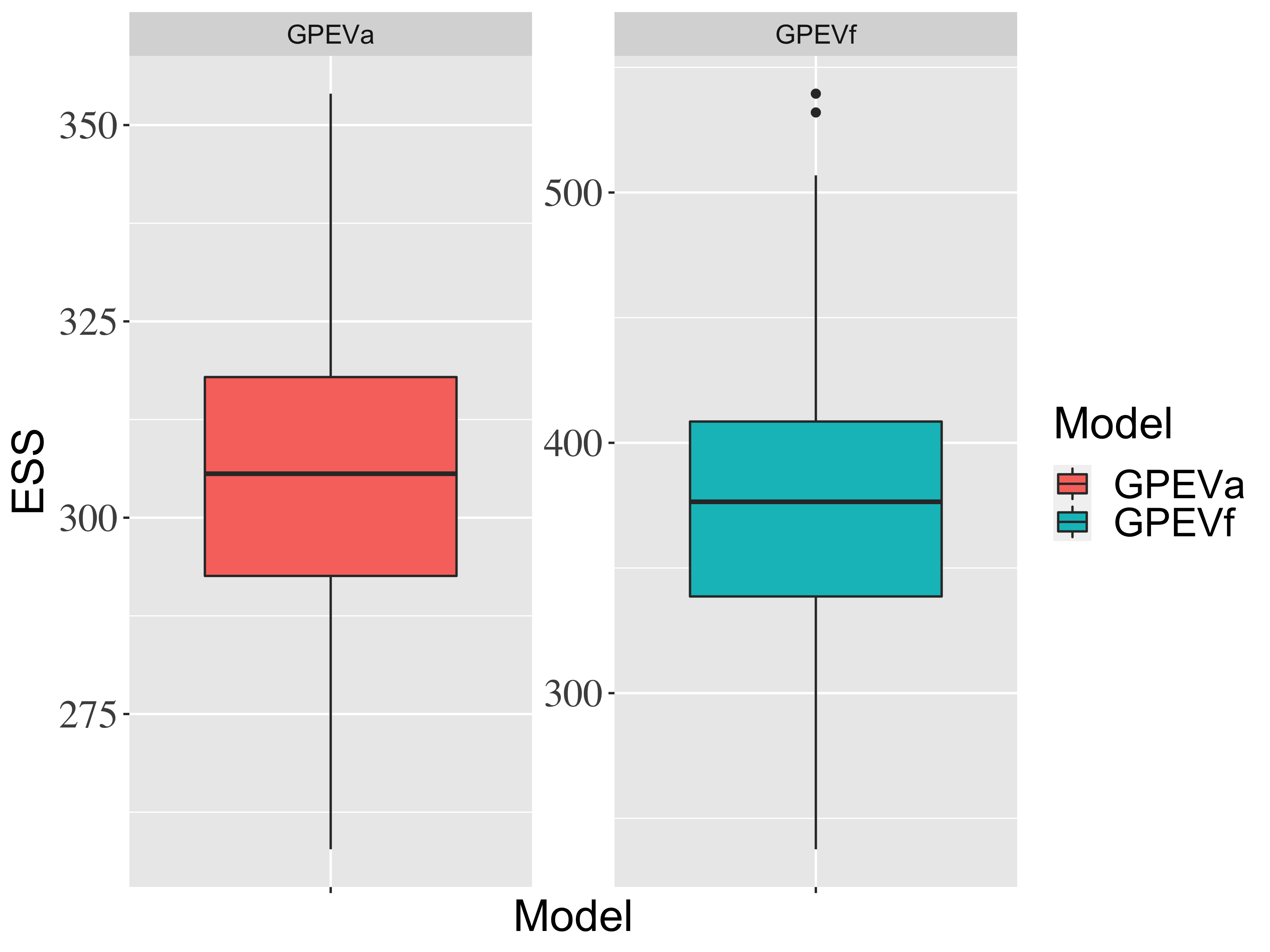}
    \includegraphics[scale = 0.045]{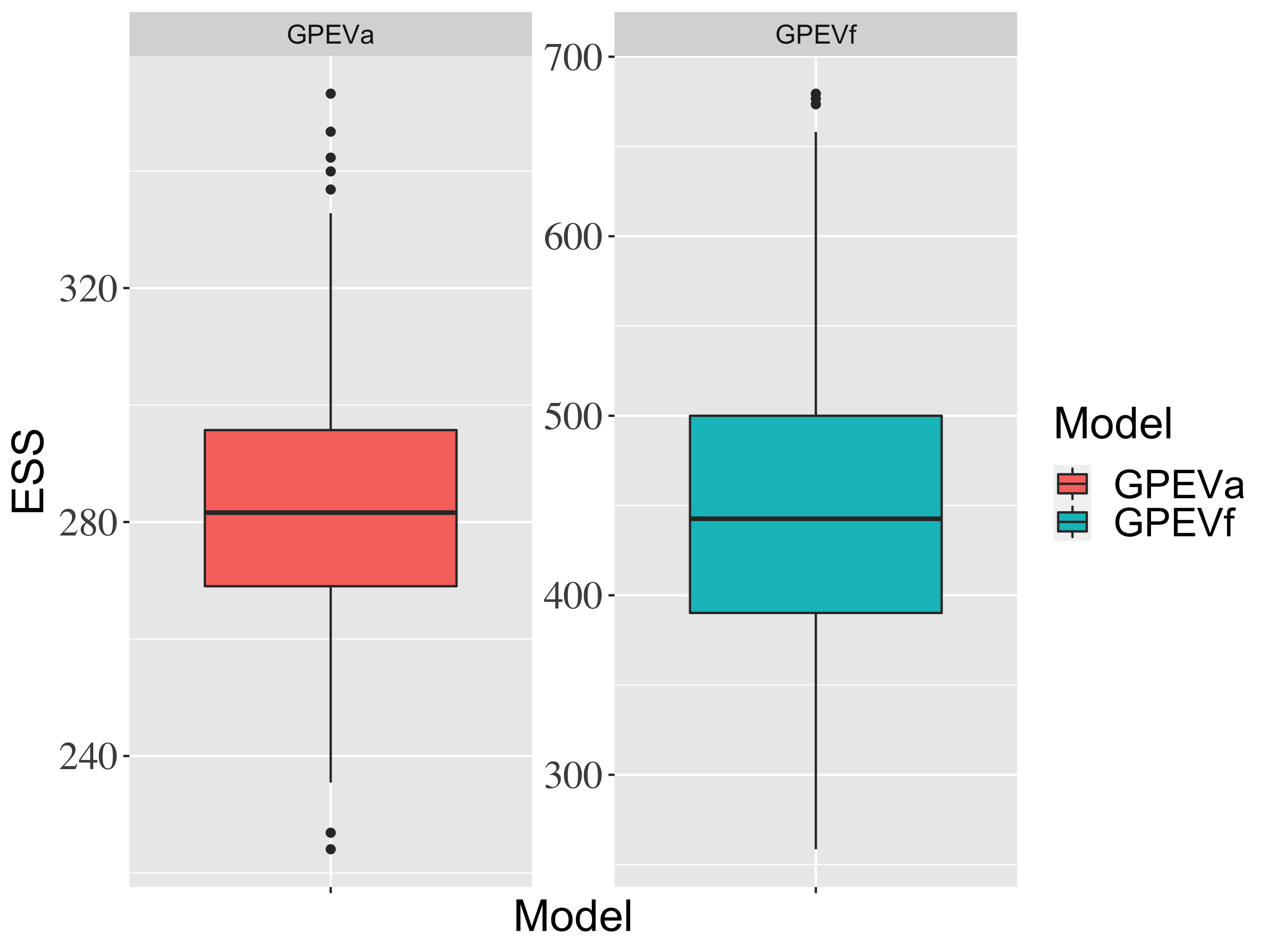}
\end{multicols}

\caption{Boxplots of effective sample sizes of function value estimated over training data points based on \textsc{gpev}$_a$ and \textsc{gpev}$_f$ over replicated data sets of sizes $n=100$ (left panel),  $n=250$ (middle panel), and $n=500$ (right panel). The effective sample sizes are averaged over 50 replicates with $\delta^2=0.01$ for $n=100,250$ and $\delta^2=0.001$ for $n=500$.}\small
\label{fig:ess-boxplot}

\end{figure}

\begin{figure}[h!]
\begin{center}
    \includegraphics[scale = 0.062]{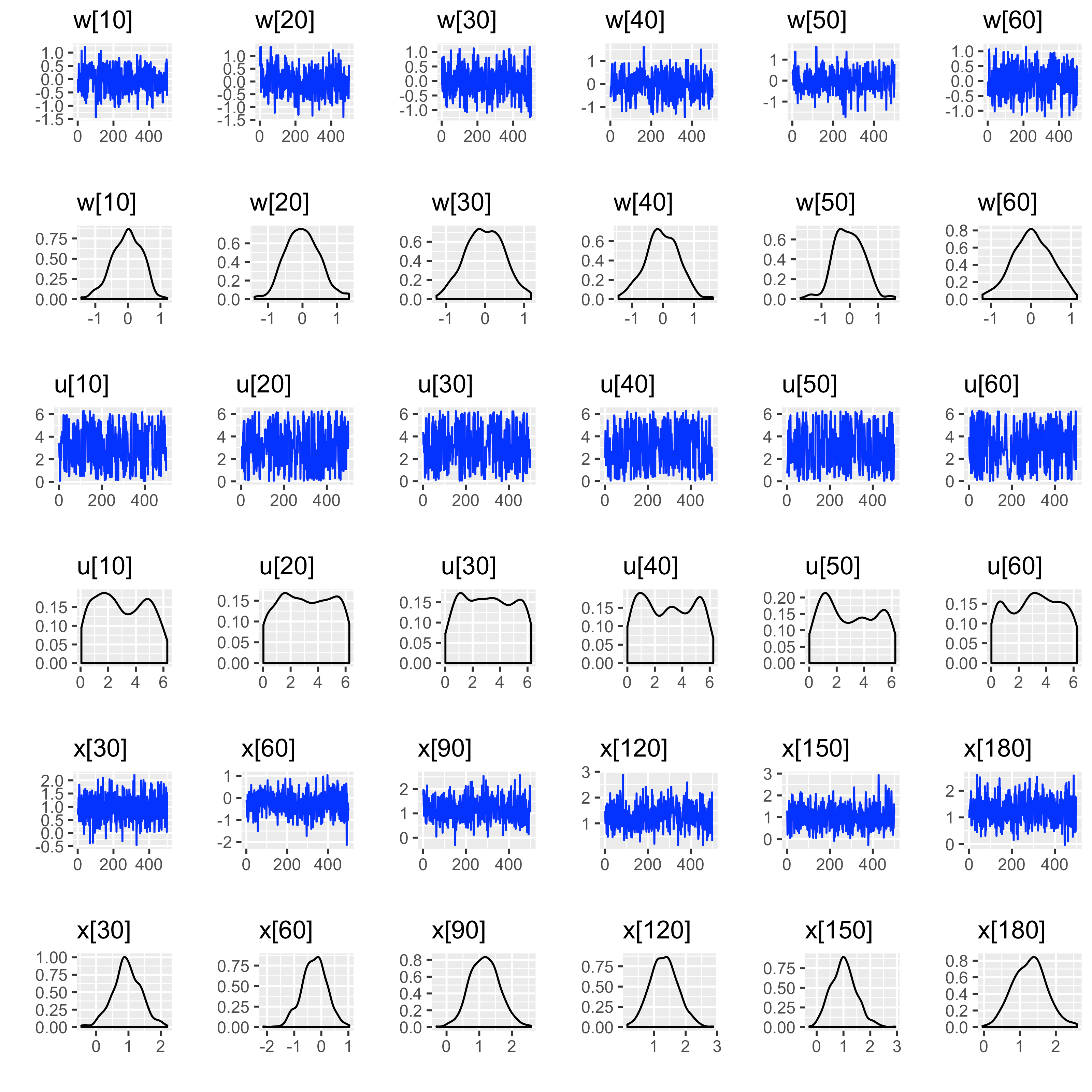}  
    \includegraphics[scale = 0.062]{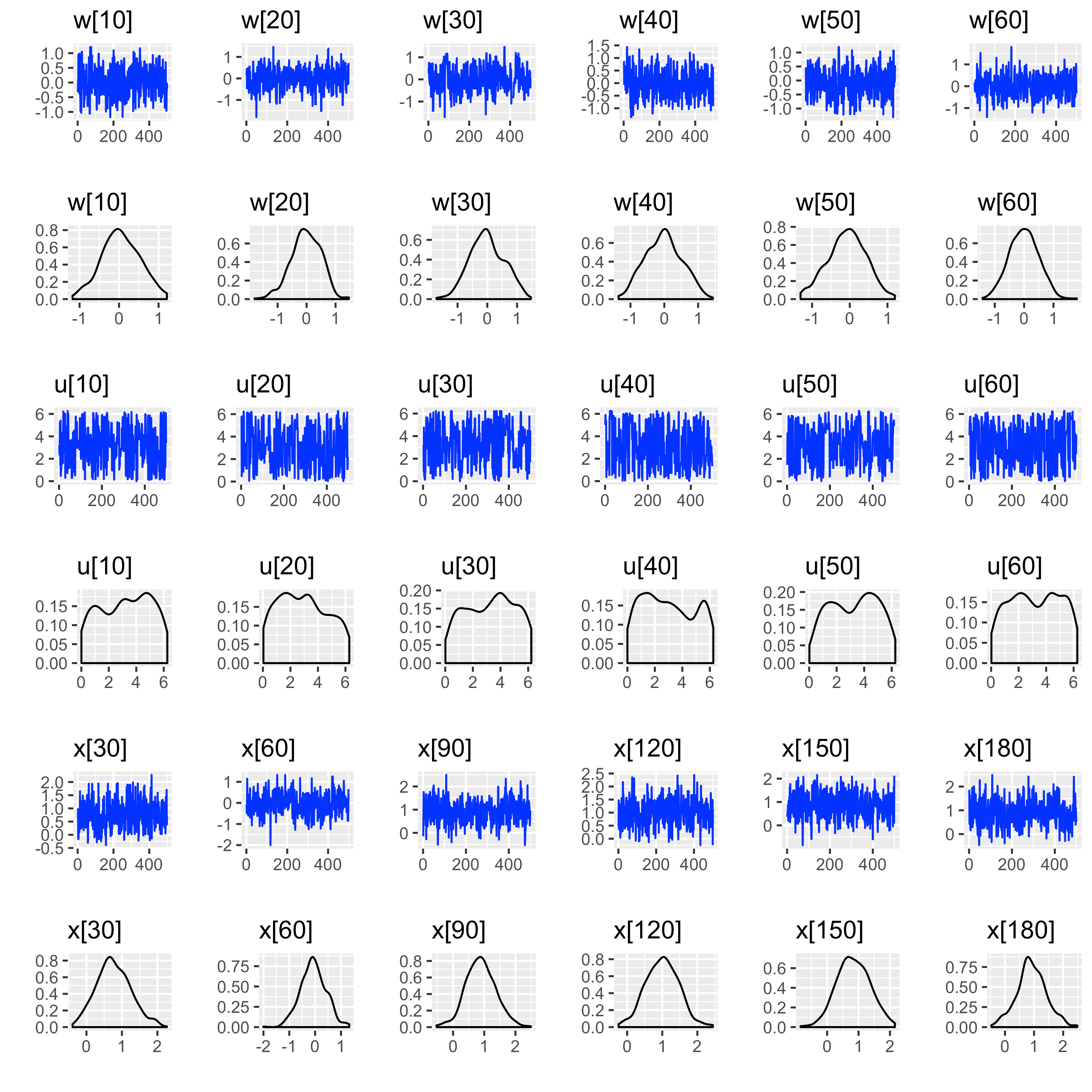}
 
\end{center}
\caption{Trace plots and density plots of the 500 posterior samples of a subset of $\{w_j, s_j, x_j\}$ from treatment group with $\delta^2 = 0.35$ (left panel) and with unknown $\delta^2$ (right panel) in the data example in Section~\ref{sec:real}.}
\label{fig:trace_treat}
\end{figure}

In this section, we provide additional numerical results for $n=250$ under the same setting in Section~\ref{sec:sims}, refer to Table \ref{tab:f_250}, Figure \ref{fig:boxplot_n250} and Figure \ref{fig:fit_n250}. {\tcr We include the \textsc{amse} values for estimating the true locations for \textsc{gpev}$_a$ and \textsc{gpev}$_n$ in Table \ref{x_mse} under all three settings of sample sizes.} We also collect diagnostic summaries under the settings in Section~\ref{sec:sims} including the mixing of the Markov chain of hyperparameter associated with covariance kernel in Figures \ref{fig:mixing_lambda}, marginal posterior density plot of covariate based on \textsc{gpev}$_a$ (Figure \ref{fig:p0}) and effective sample sizes for estimated function values over training data points  for \textsc{gpev}$_a$ and \textsc{gpev}$_f$ in Figure \ref{fig:ess-boxplot}. At last, we provide trace plots and density plots of parameters associated with \textsc{gpev}$_a$ (Figure \ref{fig:trace_treat}) for the real application in Section~\ref{sec:real}. 

\begin{table}[htbp]
\footnotesize
  \begin{center}
    \begin{tabular}{clccccccc} \toprule 
     &  &  \multicolumn{6}{c}{$\delta^2$} \\
         \cmidrule(lr){3-8}
   $n$ & Method  & 0$\cdot$01 & 0$\cdot$2 & 0$\cdot$4 & 0$\cdot$6 & 0$\cdot$8 & 1 \\
     \cmidrule(lr){1-8}
  \multirow{5}{*}{$250$}
   &\textsc{gpev}$_a$    &  0$\cdot$23 (0$\cdot$07)  &   0$\cdot$86 (0$\cdot$51)   &  2$\cdot$10 (2$\cdot$33)  &   3$\cdot$44 (4$\cdot$83)   &  3$\cdot$60 (4$\cdot$86)  &  4.$\cdot$86 (6$\cdot$21)  \\
    & \textsc{gpev}$_f$     &  0$\cdot$21 (0$\cdot$07)   &  0$\cdot$78 (0$\cdot$46)    &  1$\cdot$62 (0$\cdot$98)  &  2$\cdot$80 (2$\cdot$84)  &  2$\cdot$94 (3$\cdot$44)  &  4$\cdot$26 (4$\cdot$91)  \\
    & \textsc{gpev}$_n$       &  0$\cdot$24 (0$\cdot$09) &  4$\cdot$24 (1$\cdot$25)  &  10$\cdot$41 (2$\cdot$78)  &  14$\cdot$43 (3$\cdot$61)  &  18$\cdot$28 (4$\cdot$72)  &  20$\cdot$23 (4$\cdot$77)  \\
    & \textsc{gp}       &  2$\cdot$31 (0$\cdot$15)  &  4$\cdot$44 (0$\cdot$62)  &  7$\cdot$38 (1$\cdot$13)  &  10$\cdot$06 (1$\cdot$50) &  12$\cdot$28 (1$\cdot$72)  &  14$\cdot$29 (1$\cdot$85)  \\
 & \text{decon}   &  0$\cdot$48 (0$\cdot$27)   &  2$\cdot$99 (0$\cdot$94)  & 7$\cdot$45 (1$\cdot$62)  &  11$\cdot$91 (2$\cdot$01)  &  15$\cdot$57 (1$\cdot$99)  &  18$\cdot$17 (1$\cdot$77)  \\ \toprule
\end{tabular}%
\end{center}
\caption{Averaged mean squared errors (\textsc{amse}) defined as $\bbE\,[ K^{-1}\sum_{k=1}^K \{\,\widehat{f}(t_k) - f(t_k) \,\}^2]$ ($\widehat{f}(\cdot)$ denotes the proposed estimator of $f$, $\bbE(\cdot)$ denotes taking average over replicates) on an evenly spaced grid $(t_1, \ldots, t_K)$ of size $K=100$ over the interval $[-3, 3]$ and standard errors ($\times 10^{2}$) over 50 replicated data sets of size $n = 250$.}
\label{tab:f_250}%
\end{table}%

\begin{table}[htbp]
\footnotesize
  \centering
       \begin{tabular}{clccccccc} \toprule 
     &  &  \multicolumn{6}{c}{$\delta^2$} \\
         \cmidrule(lr){3-8}
   $n$ & Method  & 0$\cdot$01 & 0$\cdot$2 & 0$\cdot$4 & 0$\cdot$6 & 0$\cdot$8 & 1 \\
          \cmidrule(lr){1-8}
    \multirow{2}[0]{*}{100} & \textsc{gpev}$_a$ & 0$\cdot$92 (0$\cdot$13) & 13$\cdot$27 (2$\cdot$04) & 26$\cdot$41 (4$\cdot$85) & 37$\cdot$88 (6$\cdot$11) & 47$\cdot$68 (8$\cdot$93) & 57$\cdot$79 (10$\cdot$41) \\
          & \textsc{gpev}$_n$ & 0$\cdot$94 (0$\cdot$13) & 12$\cdot$86 (2$\cdot$12) & 34$\cdot$17 (4$\cdot$97) & 50$\cdot$73 (7$\cdot$76) & 64$\cdot$25 (11$\cdot$2) & 76$\cdot$86 (10$\cdot$57) \\
    \multirow{2}[0]{*}{250} & \textsc{gpev}$_a$ & 0$\cdot$89 (0$\cdot$09) & 12$\cdot$36 (1$\cdot$52) & 23$\cdot$70 (2$\cdot$95) & 33$\cdot$56 (4$\cdot$11) & 42$\cdot$98 (5$\cdot$16) & 52$\cdot$07 (6$\cdot$38) \\
          & \textsc{gpev}$_n$ & 0$\cdot$91 (0$\cdot$09) & 15$\cdot$95 (1$\cdot$59) & 33$\cdot$10 (3$\cdot$41) & 48$\cdot$49 (4$\cdot$47) & 62$\cdot$64 (5$\cdot$80) & 74$\cdot$40 (6$\cdot$54) \\
          &       & \multicolumn{6}{c}{$\delta^2$} \\
              \cmidrule(lr){3-8}
          $n$ & Method  & 0$\cdot$001 & 0$\cdot$005 & 0$\cdot$01 & 0$\cdot$1 & 0$\cdot$5 & 1\\
        \cmidrule(lr){1-8}          
    \multirow{2}[0]{*}{500} & \textsc{gpev}$_a$ & 0$\cdot$098 (0$\cdot$006) & 0$\cdot$46 (0$\cdot$03) & 0$\cdot$88 (0$\cdot$05) & 6$\cdot$54 (0$\cdot$53) & 28$\cdot$04 (2$\cdot$28) & 50$\cdot$63 (4$\cdot$36) \\
          & \textsc{gpev}$_n$ & 0$\cdot$098 (0$\cdot$006) & 0$\cdot$47 (0$\cdot$03) & 0$\cdot$89 (0$\cdot$05) & 7$\cdot$61 (0$\cdot$68) & 40$\cdot$26 (2$\cdot$90) & 74$\cdot$42 (5$\cdot$69) \\ \toprule
    \end{tabular}%
        \caption{Averaged mean squared errors  ($\times 10^{2}$) with standard errors ($\times 10^{2}$) in estimating the true locations defined as $\bbE\,[ n^{-1}\sum_{k=1}^n(\,\widehat{x}_k - x^\ast_k \,)^2]$, where $\{\widehat{x}_k\}$ denote the posterior estimate (mean) of covariates and $\{x^\ast_k\}$ denote the true locations, and $\bbE(\cdot)$ denotes taking average over replicates. The \textsc{amse} values and standard deviations are averaged over 50 replicated data sets of size $n =100, 250, 500$ separately.}
  \label{x_mse}%
\end{table}%

\newpage

\newpage 

\clearpage

\vskip 0.2in
\bibliography{gpev_1D}

\end{document}